\documentclass{article}
\usepackage{amsmath,latexsym,amssymb,amsthm,enumerate,amscd}
\newtheorem{theorem}{Theorem}[section]
\newtheorem{proposition}[theorem]{Proposition}

\newtheorem{definition}[theorem]{Definition}
\newtheorem{conjecture}{Conjecture}
\newtheorem{remark}[theorem]{Remark}
\newtheorem{lemma}[theorem]{Lemma}
\newtheorem{example}[theorem]{Example}

\newtheorem{thm}{Theorem}

\newtheorem{rem}[thm]{Remark}

\newtheorem*{uprop}{Proposition}
\newtheorem*{uthm}{Theorem}

\newtheorem*{ack}{Acknowledgment}

\numberwithin{equation}{section}
\numberwithin{table}{section} %TTT
\DeclareMathOperator{\Grass}{\rm{Grass}}

\def\BGrass(T){\mathrm{BGrass(T)}}
\def\GA{\mathrm{GA}}
\def\LA{\mathrm{LA}}

\def\cha{\mathrm{char}\ }

\def\In{\mathrm{In}}

\def\G_T{\mathrm{G_T}}

\def\Ann{\mathrm{Ann}\ }

\def\QRAM{\mathrm{QRAM}}

\def\Hilb{\mathrm{Hilb}}

\def\SGrass(T){\mathrm{SGrass(T)}}

\def\Grass{\mathrm{Grass}}

\def\Z_T{\mathrm{Z_T}}
\def\Z(E){{\bf{Z}}(E)}
\providecommand{\bysame}{\makebox[3em]{\hrulefill}\thinspace}
\def\cod{\mathrm{cod}\ }
\def\<{\left<}
\def\>{\rangle}

\def\ee{\end{equation}}

\def\Sec{\mathrm{Sec}}
\def\Sym{\mathrm{Sym}}

\def\ns{\footnotesize \it}

\def\Z{\mathfrak{Z}}

\def\RAM{\mathrm{RAM}}

\def\max{\mathrm{max}}

\title{The family $\G_T$ of graded Artinian quotients of
$k[x,y]$ of given Hilbert function}
\author{Anthony Iarrobino\\[.05in]
{\ns Department of Mathematics, Northeastern University, Boston, MA 02115, USA.
}\\[.2in]
Joachim Yam\'{e}ogo\\[.05in]
{\ns Lab. J.-A. Dieudonn\'{e}, UMR CNRS 6621,
Univ. de Nice-Sophia Antipolis, F-06108 Nice 06034 cedex 02, France
}\\[.2in]}

\date{August 26, 2002}

\begin{document}

\maketitle
%\subjclass{Primary:  14C20; Secondary 14C05}
%\footnote\ddag{AMS 2000 Subj. Class, Primary: 14C05; Secondary: 13D40} 
\begin{abstract}
 Let $R=k[x,y]$ be the polynomial ring over an
algebraically closed field $k$ of characteristic zero of
$\cha k>j$. Let $T=(1,2,\ldots ,\mu ,t_\mu ,t_{ \mu +1},\ldots ,t_j,0)
$ with $\mu \geq t_\mu
\geq t_{\mu +1} \geq \ldots \geq t_j \geq 0$ be a sequence of
nonnegative integers.  The nonsingular projective variety $\G_T$
parametrizes all graded ideals
$I$ of $R = k[x,y]$ for which the Hilbert function $H(R/I)=T$. Let $E$
be a monomial ideal of $R$, satisfying $H(R/E)=T$: the cell
$\mathbb{V}(E)$ is the family of graded ideals having initial monomials
$E$ in a suitable partial order. The graded ideal $I$ belongs to the
cell
$V(E)$ if the divisibility of a standard basis of each degree-i
homogeneous piece $I_i$ by powers of $x$ matches that of $E$; this
specifies the ramification $\QRAM (I_i,x)$ partition of the linear
system
$L(I_i)$ determined by $I_i$ at the point $p_x: x=0$ of
the projective line ${\check{\mathbb P}}^1$ (Proposition
\ref{prop4.5}). Likewise, letting $p$ be a point of
${\check{\mathbb P}}^1$ defined by $ax+by=0$, we define cells
${\mathbb V}(E,p)$, specifying the ramification at the point $p: ax+by =
0$.  When $j=
\mu$,
${\mathbb V}(E,p)$ is a Schubert cell of a Grassman variety
$\Grass(d,R_j), d= j+1-t_j$. In this case the intersection of
ramification conditions at different points is given by the Schubert
calculus (see \cite{I2,EH1,EH2}). Our goal is to generalize this
calculus to determine the homology ring $H^*(\G_T)$, in terms of the
classes of the cells
${\mathbb V}(E)$ giving the natural cellular decomposition of $\G_T$, a problem that
remains open in general.
\par
   We show that $\G_T$ is birational to a certain product $\SGrass(T)$
of small Grassman varieties (Proposition \ref{prop4.11}). We show (Theorem \ref{thm5.1})
\begin{uthm} Let ${k} = {\mathbb C}$. The birational map
from $\G_T$ to $\SGrass(T)$ induces an additive $\mathbb Z$ isomorphism
$\tau : H^{\ast}(\G_T) \longrightarrow H^{\ast}(\SGrass(T))$ of homology
groups.
\end{uthm}
  When $\mu<j$ this isomorphism is not usually an isomorphism of
rings.  We determine the ring $H^*(\G_T)$ when $T=T(\mu,j)=(1,2,\ldots
,\mu -1,\mu ,\ldots ,\mu,t_j=1)$, where $\G_T \subset \mathbb
P^\mu\times
\mathbb P^j$ (Theorem~\ref{thm7.3}). In this case $\G_T$ 
 is a desingularisation
of the $\mu$-secant bundle $\Sec (\mu,j)$ of the degree-$j$ rational
normal curve, or, equivalently, of the determinantal variety of
rank-$\mu$ 
$\mu +1$ by $j+1-\mu$ Hankel matrices. 
 We use this ring $H^*(\G_T)$ to determine the number of ideals
satisfying an intersection of ramification conditions at different
points (Example
\ref{ex7.4}). We also determine the classes in $H^*(\G_T)$ of the
pullback of the singular locus of
$\Sec  (\mu ,j)$ and of the pullbacks of the higher singular
loci -- the rank-$i$, with $i<\mu$ determinantal varieties
of the Hankel matrix (Theorem \ref{thm8.0}). \par
A main tool we use from Part I \cite{IY2} is the hook code
 or pruning code, a map $\mathcal D$ taking
the partition $P(E)$ determined by the monomials of $E^c$, the
complement of the monomial ideal $E$, to a sequence $\mathcal D(P(E))$ of smaller
partitions, each corresponding to a cell in a small Grassmann variety
of size determined by $T$. 
%\bigskip \noindent
%\hrule width4cm height0.5pt \smallskip \noindent
%AMS 1991 Subject Classification: Primary 14C20; 
%Secondary 13C40,14C05,14M12,14M15. \par \noindent
%Keywords:  Artin algebras, birational, cellular decomposition,
%desingularisation, determinantal variety, diagonal lengths, graded ideal,
%Grassmann, Grobner basis, Hankel matrix, Hilbert scheme, homology ring,
% hook code, initial ideal,linear system, parametrization, partitions,
% punctual Hilbert scheme,  ramification, rational curve, secant bundle,
% Schubert calculus, semismall resolution, 
%  Wronskian. 
\end{abstract}
%%%%%%%%%%%%%%%%%%%%%%%%%%%%%%%%%%%%%%%%%%%%%%%%%%%%%%%%%%
%%%%%%%%%%%%%%%%%%%%%%%%%%%%%%%%%%%%%%%%%%%%%%%%%%%%%%%%%% 

%%%%%%%%%%%%%%%%%%%%%%%%%%%%%%%%%%%%%%%%%%%%%%%%%%%%%%%%%%
\section{Introduction}\label{intro}
%%%%%%%%%%%%%%%%%%%%%%%%%%%%%%%%%%%%%%%%%%%%%%%%%%%%%%%%%%
The Hilbert function $H(A)$ of an Artin algebra $A$
over a field {k} specifies the dimension of each
degree-$i$ piece $A_{i}$ of $A$.  
%%%%%%%%%%%%%%%%%%%%%%%%%%%%%%%%%%%%%%%%%%%%%%%%%%%%%%%%%%
Let $T$ be a sequence of nonnegative integers possible for
a Hilbert function of an Artin quotient $A=R/I$ of the
polynomial ring $R = {k}[x,y]$. Then 
\begin{equation}\label{e0.1}
T = (1,2, \ldots , \mu,t_{\mu}, \ldots , t_{j},0), \  
where \ \mu \geq t_{\mu}\geq \cdots \geq t_{j}>0. 
\end{equation}
Here $\mu$ is the \emph{order} of $T$, the order or initial degree of
the ideal $I$. Letting $M=(x,y)$ be the maximal ideal of $R$, we have
$M^\mu\supset I\supset M^{j+1}$.
%%%%%%%%%%%%%%%%%%%%%%%%%%%%%%%%%%%%%%%%%%%%%%%%%%%%%%%%%
We denote by $\G_T$ the variety that parametrizes the
graded ideals $I$ of $R$ defining a quotient algebra    
$A = R/I$ of Hilbert function $T$.  
%%%%%%%%%%%%%%%%%%%%%%%%%%%%%%%%%%%%%%%%%%%%%%%%%%%%%%%%%
The first author had shown that  $\G_T$ is a
smooth projective variety [I1].
%%%%%%%%%%%%%%%%%%%%%%%%%%%%%%%%%%%%%%%%%%%%%%%%%%%%%%%%%
We consider these first a natural class of varieties
worthy of study in their own right.
The classical Grassmann varieties $\Grass(d,R_j)\cong \Grass(d,j+1)$
are the special case $j=\mu$.
%%%%%%%%%%%%%%%%%%%%%%%%%%%%%%%%%%%%%%%%%%%%%%%%%%%%%%%%%%
G.~Gotzmann has used the properties of $\G_T$ in his proof of
the simple connectivity of the Hilbert scheme of curves
[Gm1]; he has also used the varieties $\G_T$ to determine the
dimension of the postulation strata $\Hilb^H(\mathbb P^2)$ of the punctual Hilbert
scheme
$\Hilb^n
\mathbb P^2$ (see \cite{gotz0,gotz2},\cite[Theorem 5.51]{IK}). We will note several
other applications shortly.
%%%%%%%%%%%%%%%%%%%%%%%%%%%%%%%%%%%%%%%%%%%%%%%%%%%%%%%%%%

J.~Brian\c{c}on and the second author had studied a
certain ``vertical stratification"  of the family
$\Hilb^{n}R$ of all colength $n$ ideals in the local
ring ${k}\{ x,y\}$ (see \cite{Br,Y1,Y2}). 
%%%%%%%%%%%%%%%%%%%%%%%%%%%%%%%%%%%%%%%%%%%%%%%%%%%%%%%%%% 
L.~G\"{o}ttsche then used the methods of G. Ellingsrud and
S.~A.~Str\"{o}mme \cite{ES1,ES2} and of A.~Bialynicki -~Birula \cite{B}
to give a cellular decomposition of the variety $\mathrm{Z_T}$, parametrizing all ideals
$I\subset R$ with Hilbert function $H(R/I)=T$, and a similar decomposition of
$\G_T$ (graded ideals).
%%%%%%%%%%%%%%%%%%%%%%%%%%%%%%%%%%%%%%%%%%%%%%%%%%%%%%%%%%
Each cell $V(E)$ corresponds to a monomial ideal
$E$ of $R$ defining a quotient algebra of Hilbert
function $T$;  
%%%%%%%%%%%%%%%%%%%%%%%%%%%%%%%%%%%%%%%%%%%%%%%%%%%%%%%%%%
the cell consists of the ideals collapsing to $E$ under a
${\mathbb C}^{\ast}$-action.  The second author
showed that the Ellingsrud-Str\"omme-G\"ottsche cells
$\mathbb Z(E)$ on $\mathrm{Z_T}$ and also $\mathbb V(E)$ on $\G_T$ are in fact the
vertical cells of Brian\c{c}on \cite{Br}. In the
special case of graded ideals the cell $\mathbb V (E)$ is just the
family of graded ideals having initial monomial ideal $E$ (see
\cite{Y4}). J. Yam\'{e}ogo has studied the closure of these
cells of $\mathrm{Z_T}$ and $\G_T$ in \cite{Y3,Y4,Y5}).
%%%%%%%%%%%%%%%%%%%%%%%%%%%%%%%%%%%%%%%%%%%%%%%%%%%%%%%%%
L.~G\"{o}ttsche used the cellular decomposition to
determine the Betti numbers of $\mathrm{Z_T}$ and of $\G_T$, when ${k}=
{\mathbb C}$ \cite{G2,G4}.
%%%%%%%%%%%%%%%%%%%%%%%%%%%%%%%%%%%%%%%%%%%%%%%%%%%%%%%%%
G.~Gotzmann independently obtained a simple formula for
the Betti numbers in codimension one \cite{gotz2}.
%%%%%%%%%%%%%%%%%%%%%%%%%%%%%%%%%%%%%%%%%%%%%%%%%%%%%%%%%
Here we give a new description of the homology groups 
$H^{\ast}(\G_T)$, generalizing Gotzmann's.
%%%%%%%%%%%%%%%%%%%%%%%%%%%%%%%%%%%%%%%%%%%%%%%%%%%%%%%%%
Letting $\Delta =\Delta (T)$ denote the difference function of
$T$, so $\delta _i=t_{i-1}-t_i$, we denote by
\begin{equation}
\SGrass(T)=\prod_{\mu\le i\le j}\Grass(\delta_{i+1},
1+\delta_i+\delta_{i+1})
\end{equation}
 a product of ``small" Grassmann
varieties.
%%%%%%%%%%%%%%%%%%%%%%%%%%%%%%%%%%%%%%%%%%%%%%%%%%%%%%%%%
We first show that $\G_T$ is birational to $\SGrass(T)$
(Proposition \ref{prop4.11}). We denote by $H^{i}(\G_T)$ the
homology group in codimension $i$.  We then show (Theorem \ref{thm5.1})
%%%%%%%%%%%%%%%%%%%%%%%%%%%%%%%%%%%%%%%%%%%%%%%%%%%%%%%%%%
\begin{uthm}
Suppose that ${k} = {\mathbb C}$. The birational map
from $\G_T$ to $\SGrass(T)$, the product of small
Grassmannians induces an additive isomorphism $\tau :
H^{\ast}(\G_T) \longrightarrow H^{\ast}(\SGrass(T))$ of homology
groups: 
\begin{equation}
H^{\ast}(\G_T) \cong_{add}
\prod_{\mu(T)\leq i\leq j(T)}
H^{\ast}(\Grass(\delta_{i+1},\delta_i+\delta_{i+1}+1))
\end{equation}
 \end{uthm}
%%%%%%%%%%%%%%%%%%%%%%%%%%%%%%%%%%%%%%%%%%%%%%%%%%%%%%%%%% 
\noindent The simple nature of the isomorphism $\tau$ in terms of a ``hook code''
for partitions is striking.
%%%%%%%%%%%%%%%%%%%%%%%%%%%%%%%%%%%%%%%%%%%%%%%%%%%%%%%%%%
If $E$ is a monomial ideal of $R$, the quotient $R/E$ 
has a basis of monomials in the shape of a partition
${P}(E)$, having ``diagonal lengths" $T$ (Definition~\ref{def4.6}).
%%%%%%%%%%%%%%%%%%%%%%%%%%%%%%%%%%%%%%%%%%%%%%%%%%%%%%%%%
A hook in the shape ${ P}(E)$ consists of a corner
monomial $\omega$, an arm $(\omega,x\omega, \ldots ,\mu
= x^{a-1}\omega)$, and a foot $(\omega, y\omega, \ldots
, \mu = y^{b-1}\omega)$, with $x\mu$ and $y\mu \in E$,
but $\mu , \mu \not \in E$. Such a hook has {\em hand}
$\mu = x^{a-1}\omega$ and {\em arm minus leg difference}
$a-b$. We denote by ${\cal H}({ P})$ the set of ``difference-one'' hooks
of
$P$, those for which $a-b=1$. We show (Theorem
\ref{thm4.9})
%%%%%%%%%%%%%%%%%%%%%%%%%%%%%%%%%%%%%%%%%%%%%%%%%%%%%%%%%
\begin{uthm} The dimension of the cell
${\mathbb V}(E)$ satisfies $\dim({\mathbb V}(E))=\# {\cal H}({ P(E)})$, the
total number of difference-one hooks of the partition $P(E)$.
\end{uthm}
%%%%%%%%%%%%%%%%%%%%%%%%%%%%%%%%%%%%%%%%%%%%%%%%%%%%%%%%%
\noindent We then define a partition ${\cal Q}_{i}({
P})$ with $\delta_{i+1} = t_{i}-t_{i+1}$ parts. Each
part of ${\cal Q}_{i}({ P})$ is the number of
difference-one hooks having a given degree-$i$ hand
monomial $\mu$. 
%%%%%%%%%%%%%%%%%%%%%%%%%%%%%%%%%%%%%%%%%%%%%%%%%%%%%%%%% 
The hook code ${P} \rightarrow {\cal D}({P})$ is
the sequence  ${\cal D}({ P}) 
																				= 
							({\cal Q}_{\mu}({ P}), \ldots ,
								{\cal Q}_{j}({ P}))$.
%%%%%%%%%%%%%%%%%%%%%%%%%%%%%%%%%%%%%%%%%%%%%%%%%%%%%%%%%
The partition ${\cal Q}_{i}({ P})$ determines in the
usual way a Schubert variety and homology class in
$\Grass (\delta_{i+1}, 1+\delta_i+\delta_{i+1})$. By Theorem \ref{thm4.19} this class is
the degree-$i$ component of the isomorphism $\tau$ shown in
Theorem \ref{thm5.1}. The hook code is studied in
detail in Part I (see \cite{IY2}), and we give the results we need in Section
\ref{hookcode}.
%%%%%%%%%%%%%%%%%%%%%%%%%%%%%%%%%%%%%%%%%%%%%%%%%%%%%%%%%

A second reason for studying $\G_T$ is that it
parametrizes ideals of linear systems on the projective
line ${\check{\mathbb P}}^1$.
%%%%%%%%%%%%%%%%%%%%%%%%%%%%%%%%%%%%%%%%%%%%%%%%%%%%%%%%% 
A linear system on  ${\check{\mathbb P}}^1$ is determined by
a $d$-dimensional vector subspace $V$ of the space
$R_{j}$ of degree-$j$ forms, and consists of the zero-sets
$p_{f} \in {\check{\mathbb P}}^1$ of the elements $f\in V$
(Definition~\ref{linearsystem}).
%%%%%%%%%%%%%%%%%%%%%%%%%%%%%%%%%%%%%%%%%%%%%%%%%%%%%%%%%
The {\em degree sequence} $N_{p}(V)$ at a point $p:
ax+by=0$ of ${\check{\mathbb P}}^1$ is an increasing
sequence of integers specifying the order of vanishing at
$p$ of a ``good" basis of $V$ at $p$.
%%%%%%%%%%%%%%%%%%%%%%%%%%%%%%%%%%%%%%%%%%%%%%%%%%%%%%%%%
If $V_{u}(p)$ is the $u$-dimensional subspace of $V$
maximally divisible by $ax+by$, and $n_{i}(p,V)$ is the
power of $ax+by$ dividing $V_{d+1-i}(p)$, then $N_{p}(V) =
(n_{1}(p,V), \ldots n_d(p,V))$.
%%%%%%%%%%%%%%%%%%%%%%%%%%%%%%%%%%%%%%%%%%%%%%%%%%%%%%%%% 
The {\em ramification} of a linear system at $p$ is a
partition ${ Q}\RAM _{p}(V)$ whose parts are the
sequence $N_{p}(V)$ decreased by $(0,1, \ldots , d-1)$ (Definition \ref{def1.2}).
%%%%%%%%%%%%%%%%%%%%%%%%%%%%%%%%%%%%%%%%%%%%%%%%%%%%%%%%%
The problem of characterizing the possible ramification
sequences at different points for linear systems on a
curve was a classical topic of study.
%%%%%%%%%%%%%%%%%%%%%%%%%%%%%%%%%%%%%%%%%%%%%%%%%%%%%%%%%
When the linear system is the canonical linear system,
this was the classical study of Weierstrass points. The
first author had used Schubert calculus to study the
problem of ramification for linear systems over ${{\mathbb P}}^{1}$ 
(see the last appendix of \cite{I2}). \ 
%%%%%%%%%%%%%%%%%%%%%%%%%%%%%%%%%%%%%%%%%%%%%%%%%%%%%%%%%
D.~Eisenbud and J.~Harris elaborated this viewpoint to
study the ramification of linear systems on curves in a
series of papers (see \cite{EH1,EH2}). Our goal is to extend these results to
ideals of linear systems. \par
A third application of the varieties $\G_T$ is as a desingularization of the secant
bundles to a degree-$j$ rational curve. When
$T=T(\mu ,j)=(1,2,\ldots ,\mu ,\mu ,\ldots ,\mu ,1)$ and $2\mu <j+1$, then $\G_T$ is
a natural desingularisation of the determinantal variety $V(\mu ,j)$
of $(\mu +1)\times (\mu +1)$ minors of the generic $(\mu +1)\times (j+1-\mu )$
Hankel matrix (see section \ref{sec5C}). 
These determinantal varieties are also the secant bundles of 
the degree-$j$ rational normal curve, and as well have applications to
control theory. The homology of these determinantal varieties has 
been studied in non-projective context by A.~L.~Gorodentsov and
B.~Z.~Shapiro in \cite{GS}. If $V(\mu ,j)$ is given its
natural stratification by rank, then $\G_T $ is a semismall
resolution in the sense of intersection theory, and the homology
ring structure of $\G_T $ in this case is known (Theorem 
\ref{thm7.3}). As a first step to determining the
homology of $V(\mu ,j)$ we determine the class of the pullback of the
rank $i<\mu$ locus in Theorem~\ref{thm8.0}.
\par
For some other choices of $T$ the varieties $\G_T$ are also natural desingularizations
of interesting varieties ${\mathbf  Y}_T$: for example, take ${\mathbf  Y}_T$ to be
the closure of a Hilbert function stratum for any one of the three families of
algebras associated to the Grassmanian
$\Grass(d,R_j)$ parametrizing
$d-$dimensional vector space of forms of degree $j$. These are the family $\Grass
(d,j)$ of ancestor algebras, the family $\LA (d,j)$ of level algebras, or the family
$\GA (d,j)$ of algebras
$A=R/(V)$ determined by the ideal of $V$ \cite[Theorem 2.32]{I3}. The Hilbert function
strata in the last example
$\GA (d,j)$ are in fact the strata giving the decomposition of the restricted tangent
bundle to a rational curve in $\mathbb P^j$ as a direct sum of line bundles
\cite{GhISa}.\par A fourth reason for studying $\G_T$ is that the variety
${\mathbf  Z}_T $ parametrizing all ideals of Hilbert function $T$ in ${k}[[x,y]]$ is
fibred over $\G_T$ by an affine space of known dimension; and the union of the
$\mathrm{Z_T}$ for all $T$ of fixed length
$n$ form the punctual Hilbert scheme $\Hilb^{n}{k}[[x,y]]$, the fibre of
$\Hilb^n(\mathbb P^2)$ over a point $np_0$ of the symmetric product $\Sym^n(\mathbb
P^2)$. The punctual Hilbert scheme has shown itself important lately in several
ways. M.~Haiman has used 
$\Hilb^{n}{\Bbb C}[[x,y]]$ in solving an
important combinatorial problem, an $n!$-conjecture involving representation theory
\cite{haim}.  J.~Cheah has recently determined the homology groups of
nonsingular varieties $\Hilb^{n,n-1}{k}[[x,y]]$ parametrizing pairs $(I_{n}\subset
I_{n-1})$ where
$I_{n}$ is a colength-$n$ ideal and $I_{n-1}$ is a colength-$n-1$ ideal in ${k}[[x,y]]$. She also uses the difference-$a$ hooks in her description of
these groups (\cite{chea}). Her study and also her compact review of the known homology
results for $\Hilb^{n}{k}[x,y]$ -due to G.~Gotzmann, L.~G\"ottsche, et al -
is a complement to our study of the homology of the fine strata $\G_T$. \par
Motivated by physics, I. Grojnowski, H. Nakajmima, and many others have 
found deep connections among the homology for different
$n$ of the Hilbert scheme $\Hilb^n(X)$ for surfaces $X$ (see,
for example, \cite{nakaj2,nakaj3,LQW}). On the one hand, the projective varieties
$\G_T$ and the bundles
$\mathrm{Z_T}$ over them may be regarded as subschemes of the local punctual Hilbert
scheme, the fibre over the point $n\cdot p_0$ of the symmetric product. On the other
hand, the varieties $\G_T$ appear to be rather more complicated, since the Picard group
$\mathrm{Pic} (\Hilb^n
\mathbb P^2)\otimes \mathbb Q=\mathbb Q\oplus \mathbb Q$, whereas the Picard group of
$\G_T$ is by Theorem~\ref{thm5.1} in general of higher rank.
\par
This article is a revision of the preprint \cite{IY1}; we have
changed the organization, title and exposition throughout. It is the geometric portion,\
orginially entitled ``Part II'' of a work whose first avatar was a 1991/1992 preprint.
We refer here to the preprint \cite{IY2} on the combinatorial aspects of the hook code
for partitions as ``Part I'', even though this article is no longer titled ``Part II'',
and we plan to split ``Part I'' into two shorter articles.\par
%%%%%%%%%%%%%%%%%%%%%%%%%%%%%%%%%%%%%%%%%%%%%%%%%%%%%%%%% 
We now summarize the sections. In Section \ref{sec2} we
study the special case $\mu =j$ where there is a single linear system, rather than an
ideal. We define the ramification of a linear system at a point of
${\check{\mathbb P}}^1$ and we connect the ramification -- a partition -- with the cells
${\mathbb V}(E)$ parametrizing vector spaces having a given initial monomial vector
space $E$. We then use the Wronskian of a vector space of forms $V$ in
$R_j$ to determine the total ramification number of $V$ at all points of
${\check{\mathbb P}}^1$.
%%%%%%%%%%%%%%%%%%%%%%%%%%%%%%%%%%%%%%%%%%%%%%%%%%%%%%%%%
In Section \ref{Grasscells} we apply the Schubert calculus to this case, and
in Section \ref{classicalresults} we summarize the known results on the
intersection of ramification conditions on linear systems over
${\check{\mathbb P}}^1$.
%%%%%%%%%%%%%%%%%%%%%%%%%%%%%%%%%%%%%%%%%%%%%%%%%%%%%%%%%%

In Section \ref{sec4} we extend our study of ramification to an ideal of
linear systems, our main theme. Using our results from
Section \ref{sec2} we show that the ramification loci are the same as the
cells $\mathbb V(E)$ (Proposition
\ref{prop4.5})
%%%%%%%%%%%%%%%%%%%%%%%%%%%%%%%%%%%%%%%%%%%%%%%%%%%%%%%%%
\begin{uprop} 
			Let $E$ be a monomial ideal	with $H(R/E) = T$, 
			and $p \in {\check{\mathbb P}}^1$, and suppose that 
			$I$ is an ideal with $H(R/I) = T$. 
			The following are equivalent: 
						\begin{itemize}
													\item[i.] $I\in {\mathbb V}(E,p)$
														\item[ii.] $\QRAM (I,p) = 
																							\QRAM (E)$.
 					\end{itemize}
\end{uprop}
%%%%%%%%%%%%%%%%%%%%%%%%%%%%%%%%%%%%%%%%%%%%%%%%%%%%%%%%%
\noindent Thus, if $p: ax+by =0$ is a point of 
${\check{\mathbb P}}^1$, the condition $I\in {\mathbb V}(E,p)$,
specifies both the initial form ideal of $I$ in a certain
basis $(ax+by,C)$ for ${k}[x,y]$, and also the
ramification of $I$ at $p$. We also define the partition $P(E)$ determined by the
monomial ideal
$E$: its Ferrers graphs is that of a complementary basis $E^c$ to $E$ in $R$
(Definition \ref{def4.6}).
%%%%%%%%%%%%%%%%%%%%%%%%%%%%%%%%%%%%%%%%%%%%%%%%%%%%%%%%%
In Section \ref{sec4B} we show  we show that the cell 
$\mathbb V(E)$ is an affine space, with parameters given by certain coefficients
of the generators of $I\in \mathbb V(E)$, corresponding to the difference-one hooks
of the partition
$P(E)$ (Theorem
\ref{thm4.9}). We also show that there is a birational map: $\G_T\to \SGrass(T)$, the product of small
Grassmanians (Proposition \ref{prop4.11}).  In Section \ref{dimcell}, after giving a
formula  for the dimension of the affine space fibre $F_{\mathbb{Z/V}}$ of the
cell $\mathbb Z(E)$ (all ideals) over
$\mathbb V(E)$ (graded ideals) in Proposition \ref{prop4.13B}, we reconcile our
dimension formulas for
$\mathbb {Z}(E)$ and $\mathbb V(E)$ with formulas of L. G\"{o}ttshe.  In Section
\ref{hookcode} we define the hook code of a partition $P$ of diagonal lengths $T$, and
we show that the hook code gives an isomorphism between two distributive lattices, the
second related to the product of Schubert cells on the small Grassmanians studied
earlier (Theorem
\ref{thm4.19}).  In Section \ref{homgroups} we show that there is an additive
isomorphism given by the hook code between the homology  $H^\ast (\G_T)$, and $H^\ast
(\SGrass(T))$ (Theorem \ref{thm5.1}), and we determine the Poincar\'{e} polynomial of
$\G_T$ (Theorem \ref{thm5.2}).
%%%%%%%%%%%%%%%%%%%%%%%%%%%%%%%%%%%%%%%%%%%%%%%%%%%%%%%%%
In Section \ref{sec4F} we summarize what we know about the ramification loci for
ideals $I\in \G_T$. The main limitation in comparison with the similar theory for
$\Grass(d,R_j)$ is that the intersection of ramification conditions 
$\Z = \overline{{\mathbb V}(E,p)}\cap 
 \overline{{\mathbb V}(E',p')}$ at different points is not
necessarily dimensionally proper!. 
This was shown by the second author for $T =(1,2,3,2,1)$ \cite{Y1}.
%%%%%%%%%%%%%%%%%%%%%%%%%%%%%%%%%%%%%%%%%%%%%%%%%%%%%%%%%
When $Z$ has the right codimension, then its homology
class $[Z]$ can be read from the homology ring
$H^{\ast}(\G_T)$ (Theorem \ref{thm7.0}).
%%%%%%%%%%%%%%%%%%%%%%%%%%%%%%%%%%%%%%%%%%%%%%%%%%%%%%%%%
 In particular if $Z$ is a zero-dimensional set,  the
number of vector spaces satisfying given ramification
conditions can be calculated, if we know the homology
ring $H^{\ast}(\G_T)$ in terms of the classes of the
cells $\overline{{\mathbb V}(E)}$.
%%%%%%%%%%%%%%%%%%%%%%%%%%%%%%%%%%%%%%%%%%%%%%%%%%%%%%%%%

%%%%%%%%%%%%%%%%%%%%%%%%%%%%%%%%%%%%%%%%%%%%%%%%%%%%%%%%%

In Section \ref{sec5A} we review what is known about the homology
ring structure $H^{\ast}(\G_T)$. 
There is a natural immersion 
$\iota : \G_T \longrightarrow \BGrass(T)$ into a variety
$\BGrass(T)$ that is a product of big Grassmannians
$\Grass (i) = \Grass(i+1-t_{i},i+1)$: we take the ideal $I$
into its degree-$i$ pieces, for each $i$. Recently, A. King and C.
Walter have shown that $\iota ^{\ast} : H^{\ast}(\BGrass(T))
\longrightarrow H^{\ast}(\G_T)$ is a surjection \cite{KW}. 
%%%%%%%%%%%%%%%%%%%%%%%%%%%%%%%%%%%%%%%%%%%%%%%%%%%%%%%%%
A result of the second author exhibits $\G_T$ as the
zeroes of a section of a vector bundle on a variety
$\mathrm{G_{T'}}\times \Grass(j)$ where $T'$ is simpler than $T$
\cite{Y7}. These give hope for further progress in finding the
homology ring $H^{\ast}(\G_T)$.
%%%%%%%%%%%%%%%%%%%%%%%%%%%%%%%%%%%%%%%%%%%%%%%%%%%%%%%%%
In Section \ref{sec5B} we determine the homology ring 
$H^{\ast}(\G_T)$ for 
$T(\mu ,j) = (1,2,\ldots ,\mu ,\mu ,\ldots \mu ,1)$ (Theorem \ref{thm7.3}).
%%%%%%%%%%%%%%%%%%%%%%%%%%%%%%%%%%%%%%%%%%%%%%%%%%%%%%%%%
We illustrate our approach by determining the ideals
satisfying a certain intersection of ramification
conditions on $\G_T$, $T=T(\mu ,j)$ (Example~\ref{ex7.4}). In Section \ref{sec5C}
we study the $\mu$-secant variety $\Sec(\mu,j)$ to the degree-$j$ rational normal
curve, of which $\G_T , T=T(\mu,j)$ is a desingularisation. We determine the classes
in 
$H\ast (\G_T)$ of the pullbacks of the higher singular loci of $\Sec(\mu,j)$ (Theorem
\ref{thm8.0}).
%%%%%%%%%%%%%%%%%%%%%%%%%%%%%%%%%%%%%%%%%%%%%%%%%%%%%%%%%

As a result of Theorem \ref{thm4.19}, and Proposition \ref{prop4.5}, we
have coded the ramification conditions on ideals of
$\G_T$ at a fixed point $p$ by the hook code. 
It remains to read off the ramification 
$\QRAM (I_{i},p)$ of each piece $I_{i}$ given the
code.  
This problem is studied and solved in Algorithm~2.29 of
Part I \cite{IY2}. There we defined a ``strand map" and with it we
constructed a partition ${Q}_{i}(E,p)$ that
determines $E_{i}$ directly from ${\cal D}(E)$. 
Here we show that ${Q}_{i}(E,p)$ is the dual of the
complement $\QRAM (E_{i},p)$ (Lemma \ref{lemma1.9}).
%%%%%%%%%%%%%%%%%%%%%%%%%%%%%%%%%%%%%%%%%%%%%%%%%%%%%%%%%

The sum of the parts of $\QRAM (E_{i},p)$ is the
codimension of the condition on $\Grass (i)$ that a vector
space $V$ in $R_{i}$ satifies 
$\QRAM (V,p) = \QRAM (E_{i},p)$. Concerning
the ramification of an ideal $I$ at different points, we
show, denoting by $\ell (P)$ the sum of the parts
of $P$ (Lemma \ref{lemma1.4} and Proposition \ref{prop7.5}),
%%%%%%%%%%%%%%%%%%%%%%%%%%%%%%%%%%%%%%%%%%%%%%%%%%%%%%%%%
\begin{uthm} 
If $I\in \G_T$, then for each $i$,
\begin{equation}
 \sum_{p\in {\check{\mathbb P}}^1} \ell
(\QRAM (I_{i},p)) = t_{i}(i+1-t_{i}) . 
\end{equation}
\end{uthm}\noindent
Summing over $i$, we find that the sum of the lengths of
the ramification conditions that $I$ satisfies at all
points $p$ is $$\sum_{p\in {\check{\mathbb P}}^1} \ell (\QRAM (I,p))
 = \dim(\BGrass(T) ) .$$ 
%%%%%%%%%%%%%%%%%%%%%%%%%%%%%%%%%%%%%%%%%%%%%%%%%%%%%%%%%
This result is a consequence of a stronger result that for each $i$
the set of ramification partitions \\ $\QRAM(I_{i},p), \ p\in
{\check{\mathbb P}}^1\}$, or, equivalently, the set of Schubert
classes $\{ {Q}(I_{i},p), \ p\in {\check{\mathbb P}}^1\}$ must
be ``complementary": their intersection
 can be calculated using the Littlewood-Richardson rule: the intersection must be
nonempty and nonzero provided the codimensions are not too large, and the
calculated homology class is nonzero. (Proposition \ref{prop2.4}). There is a Wronskian
morphism 
\begin{equation*}
W: \G_T
\longrightarrow {{\mathbb P}}, \text{ where }{{\mathbb P}} 
= \prod_{i=\mu}^j {{\mathbb P}}^{N_{i}} , \ n_{i}= t_{i}(i+1-t_{i}).
\end{equation*}
 The morphism $W$  is
a product of finite covering maps $w_{i}: \Grass (i)
\longrightarrow {{\mathbb P}}^{N_{i}}$ studied in Section~\ref{sec2}
(Proposition \ref{prop2.2}). The map $W$ is a
finite cover of its image (Proposition \ref{prop7.7}), but in general its image has
large codimension in ${{\mathbb P}}$. What are the equations
describing the image of $W$? These equations constitute
mysterious, hidden relations among the ramification of
$I$ at different points of ${\check{\mathbb P}}^1$. 
%%%%%%%%%%%%%%%%%%%%%%%%%%%%%%%%%%%%%%%%%%%%%%%%%%%%%%%%%
%%%%%%%%%%%%%%%%%%%%%%%%%%%%%%%%%%%%%%%%%%%%%%%%%%%%%%%%%%

This article and Part I \cite{IY2} replace the preprint [I-Y1]. 
Part I is the combinatorial portion; this is the
algebraic-geometric portion. Our intent in following a
referee's suggestion to split the article is to make the
results more accessible to specialists in each area.
%%%%%%%%%%%%%%%%%%%%%%%%%%%%%%%%%%%%%%%%%%%%%%%%%%%%%%%%%%

%%%%%%%%%%%%%%%%%%%%%%%%%%%%%%%%%%%%%%%%%%%%%%%%%%%%%%%%%%
\begin{ack}It was a question of G.~Gotzmann [Gm-1] about the 
relation of his calculation of a simple formula for the
rank of  $\mathrm{Pic}(\G_T)$ (see Theorem \ref{thm5.1}) with the
results  of L.~G\"{o}ttsche \cite{G2} that led
to our work. 
We thank J.~Brian\c{c}on, S.~Diesel, D.~Eisenbud,
G.~Ellingsrud, J.~Emsalem, A.~Geramita, G.~Gotzmann,
A.~Hirschowitz, P.~LeBarz, K.~OÕHara, M.~Merle, R.~Stanley,
A.~Suciu, S.~Xamb\'{o}-Descamps, J.~Weyman, B.~Zaslov and
A.~Zelevinsky for their comments. We thank H.~Matsumura,
organizer of Commutative Algebra and Combinatorics in
Nagoya in 1990, A.~Galligo and other organizers of the
MEGAS seminar at Nice in April 1992, E.~Kunz,
H.J.~Nastold, and L.~Szpiro, organizers of a Commutative
Algebra meeting at Oberwohlfach in 1992,  and B.~Sturmfels and D. Cox, organizers of a
portion of the 1992 Regional Geometry Institute on Computational
Algebraic Geometry. Revisions and the addition of Section \ref{sec5C} 
were made during a visit of the first author to 
the Laboratoire J.-A. Dieudonn\'e, UMR CNRS 6621, in 1997, Further revisions were
made soon after S. Kleiman's 60th birthday conference. We
thank Steve Kleiman for informative discussions and encouragement.
\end{ack}
%%%%%%%%%%%%%%%%%%%%%%%%%%%%%%%%%%%%%%%%%%%%%%%%%%%%%%%%%%

\tableofcontents

 %%%%%%%%%%%%%%%%%%%%%%%%%%%%%%%%%%%%%%%%%%%%%%%%%%%%%%%%%%
% \newpage
%%%%%%%%%%%%%%%%%%%%%%%%%%%%%%%%%%%%%%%%%%%%%%%%%%%%%%%%%%
%%%%%%%%%%%%%%%%%%%%%%%%%%%%%%%%%%%%%%%%%%%%%%%%%%%%%%%%%%
\section{Vector spaces of forms, linear systems on $\mathbb P^1$, and
Wronskians}\label{sec2}
%%%%%%%%%%%%%%%%%%%%%%%%%%%%%%%%%%%%%%%%%%%%%%%%%%%%%%%%%%
%%%%%%%%%%%%%%%%%%%%%%%%%%%%%%%%%%%%%%%%%%%%%%%%%%%%%%%%%%
We let $R$ denote the polynomial ring ${k}[x,y]$ in 
two variables, where ${k}$ is an algebraically closed
field. 
%%%%%%%%%%%%%%%%%%%%%%%%%%%%%%%%%%%%%%%%%%%%%%%%%%%%%%%%%
We let $R_{j} = \langle x^{j},x^{j-1}y,  \ldots , y^{j}\rangle$ be
the subspace of degree-$j$ homogeneous polynomials.  In Section \ref{ramifWronsk} we
first describe how a $d$-dimensional vector space $V$ of
$R_{j}$ is a linear system on the projective line ${\check{\mathbb P}}^1$. We then define the ramified points 
$p\in {\check{\mathbb P}}^{1}$ of $V$, and a ramification
partition $\QRAM_p(V)$ of the vector space $V$ at
each point $p\in {\check{\mathbb P}}^1$ (Definition \ref{def1.2}). Recall that
$\Grass(d,R_j)$ is the Grassmanian parametrizing $d-$dimensional subspaces of $R_j$; it
has dimension
$N=d(j+1-d)$. We define a Wronskian morphism from
$\Grass(d,R_j)$ to $\mathbb P^N$, and use it to show that the total lengths
of the ramification partitions for a fixed $V$ is
$N =
\dim (V)\cdot \cod (V)$ (Definition \ref{def1.4},Lemma \ref{lemma1.4}). We will also
show that
$\QRAM_p(V)$ is determined by the initial monomials for $V$ in a particular
basis for $R$ (Lemma \ref{lemma1.9}). Given $E$ a set of monomials of $R_j$ we define
the ``cell''
$\mathbb V(E)$ as the subscheme of the Grassmanian having initial monomials $E$, and
$\mathbb V(E)$ an analagous cell at the point $p\in {\check{\mathbb P}}^1$
(Definition \ref{def1.11}.\par
 In
Section \ref{Grasscells} we
 show that if $E$ is a monomial vector space, the
subscheme ${\mathbb V}(E,p) \subset \Grass(d,R_{j})$ is a Schubert cell of the
Grassmannian (Lemmas
\ref{lemma2.1A},\ref{lemma2.1B}). We also show that the Wronkian morphism 
$w: \Grass(d,R_{j})
\longrightarrow {{\mathbb P}}^{N}$ is a finite cover
(Proposition \ref{prop2.2}). The intersection of ramification loci
at $s$ different points of ${\check{\mathbb P}}^1$ is nonempty and has a class given
by the Schubert calculus, provided their codimensions sum to less than $N$ and that
class is nonzero (Proposition \ref{prop2.4}). In  Section \ref{classicalresults} we
summarize in this context what is known classically for the intersection of
ramification loci of linear systems at distinct points of ${\check{\mathbb P}}^1$. We
will give a similar summary in Section \ref{sec4F} for the varieties
$\G_T$. 
%%%%%%%%%%%%%%%%%%%%%%%%%%%%%%%%%%%%%%%%%%%%%%%%%%%%%%%%%
\subsection{Ramification and the Wronskian}\label{ramifWronsk}
\begin{definition}\label{linearsystem} 
{\sc linear system of a vector space} $V\subset R_{j}$.
\   If $p = (\alpha_{p},\beta_{p})$ is a point of the
projective space ${\check{\mathbb P}}^1 = {{\mathbb P}}(R_{1}^{\vee})$, we let $L_{p}$ be any linear form
$ax+by \in R_{1}$ vanishing at $p: a\alpha_{p}+b\beta_{p}
= 0$;
the vector space $\langle L_{p}\rangle $ and the class 
$\{ L_{p} \  up \  to \  k^{\star}-multiple \}$ in 
${{\mathbb P}}^{1} = {{\mathbb P}}(R_{1})$ are uniquely determined
by $p$. 
Conversely, given $L \inÊR_{1}$ we let $p_{L} \in {\check{\mathbb P}}^1 = {{\mathbb P}}(R_{1}^{\vee})$ be the point
where it vanishes.
The vector space $\langle f\rangle $ spanned by a single degree-$j$ 
homogeneous polynomial $f$ in $R = {k}[x,y]$, $f =
L_{1}\cdots L_{j}$ corresponds to a zero cycle  ${\bf
p}{_f} = p_{L_{1}}+\cdots + p_{L_{j}}$ of $j$
points - counting multiplicities - on the projective line
${\check{\mathbb P}}^1$. 
A $d$\ -dimensional vector subspace $V$ of         
$R_{j}$ corresponds to a unique linear system on ${\check{\mathbb P}}^1$ :
\begin{equation}\label{e1.1}
{\cal L}(V) = \{ p_{f} \ | \ f\in V\} \subset
{\check{\mathbb P}^j} =
\Sym^{j}({\check{\mathbb P}}^1). 
\end{equation}
${\cal L}(V)$ is a linear subspace of ${\check{\mathbb P}^j}$ and has projective dimension $(d-1)$.
\end{definition}
%%%%%%%%%%%%%%%%%%%%%%%%%%%%%%%%%%%%%%%%%%%%%%%%%%%%%%%%%%
%xxx\newtheorem{exam}{Example}[section] \begin{exam}
\begin{example}\label{ex1.1}
 {\em Let}
\begin{equation*}
 \begin{array}{ccccccc}
 V_{a}  &=  &\langle y^{2}x-a^{2}x^{3}, &yx^{2}+ax^{3}\rangle  	&=
 &\langle y(xy+ax^{2}), &x(xy+ax^{2})\rangle .\\
  & &f &g   &       	&ag+f    	 &g\\
  \end{array}
\end{equation*}
The general element of the one-dimensional family
${\cal L}(V_{a})$  is the zero set $
p_{F_{t}}$ of
\begin{equation*}
 F_{t}=t(y^{2}x-a^{2}x^{3})+(1-t)(yx^{2}+ax^{3}).
\end{equation*} 
Here, ${\cal L}(V_{a})$  consists of all sets of
three  points in ${\check{\mathbb P}}^1$ including
both $P : x = 0$ and $P': y+ax=0$, the base
points of ${\cal L}(V_{a})$. 
 In particular, ${\cal L}(V_{a})$  contains
$2P+P'$  and $2P'+P$. (See Figure \ref{fig-linsyst}).
\end{example}
\begin{figure}[hbtp]
\begin{center}
\leavevmode
$$\begin{array}{c}
\mbox{\small{\mbox{$\left\{ \begin{array}{ccccc}
 & 2p & + & p' &  
\\
 -- & \mbox{\large{$\bullet$}} & --- & \bullet & -- 
\\
& x=0 &  & y+ax=0 &  
\end{array}
\right\}$}}} 
\\ 
\\ 
\mbox{\small{\mbox{$\left\{ \begin{array}{ccccccc}
& p'' & + & p & + & p' &  
\\
 -- & \bullet & --- & \bullet & --- & \bullet & -- 
\\
& y-ax=0 & & x=0 & & y+ax=0 &  
\end{array}
\right\}$}}}
\\ 
\\
\mbox{\small{\mbox{$\left\{ \begin{array}{ccccc}
 & p & + & 2p' &  
\\
 -- & \bullet & --- & \mbox{\large{$\bullet$}} & -- 
\\
 & x=0 &  & y+ax=0 & 
\end{array}
\right\}$}}}
\end{array}
$$
\end{center}
\caption{Three elements of the linear system ${\mathcal L}(V_{a}).$ See Example
\ref{ex1.1}.
\label{fig-linsyst}}
\end{figure}
%
%\vspace{15mm}
%\special{picture fig/exam(1) scaled 720}      
%
%\vspace{5mm}
%
%%%%%%%%%%%%%%%%%%%%%%%%%%%%%%%%%%%%%%%%%%%%%%%%%%%%%%%%%%%%
%\hspace{3mm}
%\underline{{\bf Figure 1.} {\em Three
%elements of the linear system} ${\cal L}(V_{a})$}.\\ 
%%%%%%%%%%%%%%%%%%%%%%%%%%%%%%%%%%%%%%%%%%%%%%%%%%%%%%%%%%%

We now define the ramification $\QRAM_p(V)$ of $V$ at a point $p\in
{\check{\mathbb P}}^1$. It will be a subpartition of the $\dim V\times \cod V$
rectangular partition $B(d,j+1-d)$ with $d$ parts of size $j+1-d$. We will consider
$\QRAM_p(V)$ to have
$d$ parts, although some of the parts may be zero.
%%%%%%%%%%%%%%%%%%%%%%%%%%%%%%%%%%%%%%%%%%%%%%%%%%%%%%%%%
%\newtheorem{def12}[exam]{Definition} 
\begin{definition}[Ramification of $V\subset R_j$ at
$p$]\label{def1.2} Suppose that
$L$  is a linear form
$L=ax+by$ corresponding to the point $p: ax+by=0$  of
${\check{\mathbb P}}^1$, that $V=\langle f_{1}, \ldots ,
f_{d}\rangle $  is a $d$-dimensional subspace of
$R_{j}$, and that $F=(f_{1}, \ldots , f_{d})$ 
satisfies
\begin{equation}\label{e1.2}
f_{i} = L^{n_{i}(p,V)}g_{i}(p,V), \ and \
n_{1}(p,V)< n_{2}(p,V)< \cdots <n_{d}(p,V). 
\end{equation}
The {\em{degree sequence}} $N_{p}(V)$ of $V$ at $p$ is 
\begin{equation}\label{edegseq}
N_{p}(V) = (n_{1}(p,V), \ldots ,n_{d}(p,V)).
\end{equation}
It is easy to see that there is a basis $F$ of $V$
satisfying \eqref{e1.2} and that $N_{p}(V)$ is independent of
the choice of $F$.
The {\em ramification sequence} $\QRAM_{p}(V)\subset B(d,j+1-d)$ of
$V$ at the point $p\in {\check{\mathbb P}}^1$ is the
partition constructed in the usual way from $N_{p}(V)$:
\begin{equation}\label{e1.3}
\QRAM _{p}(V) = (r_{1}(p,V), \ldots ,r_{d}(p,V)),
 \ \ r_{i}(p,V) = n_{i}(p,V)-(i-1). 
\end{equation} 
We say that $V$ is {\em unramified} at $p$
if $N_{p}(V) = (0,1,\ldots ,d-1)$, the minimum possible
sequence, so $\QRAM _{p}(V) = (0, \ldots
,0)$. Otherwise, $V$ is {\em ramified} at $p$.  
The {\em total  ramification} $r(p,V)$ of $V$ at
$p$ satisfies  
\begin{equation}\label{e1.4}
r(p,V)= \sum_{i}
			r_{i}(p,V)= \ell(\QRAM_p(V)),
\end{equation}
the length of $\QRAM_p(V)$. Thus, $V$ is ramified at $p$ iff the total
ramification
$r(p,V)>0$.
\end{definition}
\begin{example}\label{ex13}
Let $p$ be $x=0$, and let $E$ be the monomial vector space
\begin{equation*}
E=\{ x^{n_{1}}y^{j-n_{1}}, \ldots ,
x^{n_{d}}y^{j-n_{d}}\}, \ n_{1} < \ldots < n_{d}
\end{equation*}
then we denote by $\QRAM (E)$ the sequence $\QRAM _{x}(E) = (n_{1}-0, \ldots , n_{d}-(d-1))$, and by
$r(E)$ the total ramification 
\begin{equation}\label{etotalram}
r(E) = \ell(\QRAM(E))=(\sum_{1\leq
i\leq d}n_{i}) \ - \ d(d-1)/2
\end{equation}
 of $E$ at $p$. (Recall that $\ell(P)$ for a partition $P$ is the sum of its parts.)
\end{example}

%%%%%%%%%%%%%%%%%%%%%%%%%%%%%%%%%%%%%%%%%%%%%%%%%%%%%%%%%
\noindent {\bf Remark.} For most vector spaces $V$ of
dimension two in $R_{3}$,   there are  four distinct
 ramification points $P_{i}, \ i=1, \ldots ,4$  such
that $2P_{i}+P'_{i}$ is an element of the linear 
system ${\cal L}(V)$ for some $P'_{i}$. For the
space $V_{a}$  of Example \ref{ex1.1} there are only two
ramification points $P_1$ and $P_2$, at each of which
$\QRAM _{P_1}(V_{a}) = \QRAM _{P_2}(V_{a}) =
(1,1)$; thus, $V_{a}$ has total ramification
two at each of $P_1$, $P_2$. In either case the sum of
the total ramification of $V$ over all points of 
${\check{\mathbb P}}^1$ is \  $4 = \cod (V)\cdot \dim (V)$. We now
introduce a form, the Wronskian determinant $W(V)$ of
$V$. We will show in Lemma \ref{lemma1.4} \ that the multiplicity
of its roots at $p$ is the total ramification of $V$ at
$p$, and that its degree is $N = \cod (V)\cdot \dim (V)$.  
%%%%%%%%%%%%%%%%%%%%%%%%%%%%%%%%%%%%%%%%%%%%%%%%%%%%%%%%%

%%%%%%%%%%%%%%%%%%%%%%%%%%%%%%%%%%%%%%%%%%%%%%%%%%%%%%%%%
%\newtheorem{def14}[exam]{Definition}
\begin{definition}[Wronskian]\label{def1.4}
Suppose that $V$ is a $d$-dimensional
subspace of $R_{j}$.  We will define up to nonzero
constant multiple a degree-$N$ {\em Wronskian form} $R$, 
where $N=d(j+1-d)$. The Wronskian form determines a unique
element, the Wronskian determinant $W(V)\in {\mathbb P}(R_{N})$,
the projective space. We let $R' = {k}[x,y,dx,dy]$, the
polynomial ring and let ${\cal R} ={k}[dx,dy]$. We define
a derivation $D : R'
\longrightarrow R'$  by
\begin{align*}
D&=0 \text { on }{\cal R}\\
 D&: R\longrightarrow R' \text{ satisfies }Df=f_{x}dx+f_{y}dy.
\end{align*}
Thus if $f_{i}\in R$ and
$g_{i}\in {\cal R}$ we have $D(\sum f_{i}g_{i}) = \sum
g_{i}(f_{_{i}x}dx+f_{_{i}y}dy)$.\\
%%%%%%%%%%%%%%%%%%%%%%%%%%%%%%%%%%%%%%%%%%%%%%%%%%%%%%%%% 
If $V=\langle f_{1},\ldots ,f_{d}\rangle$  is a  $d$-dimensional
vector subspace of $R_{j}$, we let $N = d(+j+1-d) =
(\dim (V))\cdot (\cod (V))$, and we define the following
degree-$N$ homogeneous polynomial 
$W(f_{1},\ldots ,f_{d})$:
\begin{equation}\label{e1.5}
W(f_{1},\ldots ,f_{d})=\det\left(
\begin{array}{cccc}                  
           f_{1} &\cdots &\cdots &f_{d}\\
            Df_{1} &\cdots &\cdots &Df_{d}\\
           \vdots &\vdots &\vdots &\vdots\\
             D^{d-1}f_{1} &\cdots &\cdots &D^{d-1}f_{d}\\
                    \end{array}
                    \right)/ (xdy - ydx)^{d(d-1)/2} \in
											R_{N} .\  
\end{equation} 
The polynomial $W(f_{1},\ldots ,f_{d})$ is a
degree-$N$ homogeneous form so is an element of $R_{N}$.
%%%%%%%%%%%%%%%%%%%%%%%%%%%%%%%%%%%%%%%%%%%%%%%%%%%%%%%%%%
Its class $W(V)$ mod ${k}^{\ast}$-multiple is
independent  of the basis chosen for $V$,  and we define
\begin{equation}\label{e1.6}
W(V) = W(f_{1}, \ldots , f_{d})\  mod\ 
{k}^{\ast}-multiple \ in \  {\mathbb
P}(R_{N}). 
\end{equation}
We denote by $\Grass(d,R_{j})$, the Grassmannian
parametrizing $d$-dimensional vector subspaces of
$R_{j}$. The map $V \longmapsto W(V)$ defines a {\em
Wronskian morphism}
\begin{equation}\label{e1.7}
w: \Grass(d,R_{j}) \longrightarrow {{\mathbb P}}^{N}, \ \
N=d(j+1-d)
\end{equation} 
\end{definition}
%%%%%%%%%%%%%%%%%%%%%%%%%%%%%%%%%%%%%%%%%%%%%%%%%%%%%%%%%%
%%%%%%%%%%%%%%%%%%%%%%%%%%%%%%%%%%%%%%%%%%%%%%%%%%%%%%%%%%       
%\newtheorem{exam15}[exam]{Example} 
\begin{example}\label{ex15}
 When $V_{a} = \langle  y^{2}x-a^{2}x^{3},
yx^{2}+ax^{3}\rangle $ the Wronskian $W(V_{a}) =
x^{2}(y+ax)^{2}$ up to ${k}^{\ast}$-multiple.
\end{example}
%%%%%%%%%%%%%%%%%%%%%%%%%%%%%%%%%%%%%%%%%%%%%%%%%%%%%%%%%%
%%%%%%%%%%%%%%%%%%%%%%%%%%%%%%%%%%%%%%%%%%%%%%%%%%%%%%%%%%
%\newtheorem{lem16}[exam]{Theorem}
\begin{lemma}\label{lemma1.4}
   If 
$\cha {k}$ $>j$ or is $0$, and $V\subset R_{j}$,
then 
 \begin{equation}\label{e1.8}
W(V)=\prod_{p\in {\check{\mathbb
P}}^1}L_{_{p}}^{r(p,V)}
\ up
\ to \ {k}^{\ast}-multiple. 
\end{equation}
The product is over
points $p$ of ${\check{\mathbb P}}^1$. Any vector space $V$
is unramified at all but a finite number of points   of \ 
${\check{\mathbb P}}^1$. The linear form $L_p$ divides $W(V)$ iff $\exists f\in V$ such
that
$(L_p)^d|f$. We have 
\begin{equation}\label{e1.9}
\sum_{p\in {\check{\mathbb P}}^1} r(p,V) =N.
\end{equation}
\end{lemma}
\begin{proof} For \eqref{e1.8}, it suffices, considering the
$PGL(1)$ action on $W(V)$, to take $L=x$. By
substituting $y=1$, then $W(V)(x,1)$ is a usual Wronskian
determinant, and has the value \  $cx^{r(x,V)}$ \ as a
consequence of a Van der Monde style calculation. From
Definition \ref{def1.4}, degree $W(V) = d(j+1-d)$, provided
$(xdy-ydx)^{d(d-1)/2}$ divides the numerator of \eqref{e1.5},
which is easily verified. The formula \eqref{e1.9} then follows
from \eqref{e1.8}. That $L_p|W(V)\Leftrightarrow (L_p)^d$ divides some $f\in V$
follows from \eqref{e1.8}.
\end{proof}\par
We now define a partition ${ Q}(V,p)\subset B(\cod V, \dim V)$, determined by
a standard basis of $V$ in the direction $p$. 
If $p\in {\check{\mathbb P}}^1$ is the point $L=0$, and $C\in
R_{1}$ satisfies $\langle C\rangle \neq \langle L\rangle $ we use the reverse
alphabetic order on the set of degree$-j$ monomials of
$R$, $C^{j}< C^{j-1}L< \ldots <L^{j}$. If $f\in R_{j}$
we let $\In_{p}(f) =$ the initial monomial of $f$ in this
basis. If $V\subset R_{j}$, we let $\In_{p}(V) =
\langle \In_{p}(f), f\in V\rangle $, a space spanned by monomials in
$C, L$. If $E_{L}$ is spanned by monomials in $C, L$ we
let $E_{L}^{c}$ be the complementary set of monomials in
$C, L$.
%%%%%%%%%%%%%%%%%%%%%%%%%%%%%%%%%%%%%%%%%%%%%%%%%%%%%%%%%%
%\newtheorem{def17}[exam]{Definition}
\begin{definition}[Partition ${ Q}(V,p)$.]\label{def1.7}
If $V \subset R_{j}$, if $p\in {\check{\mathbb
P}}^1$ is the point $L=0$, and if $V$ has dimension $d$ and
codimension $t=j+1-d$ in $R_{j}$, we let $E_{L} =
\In_{p}(V)$, and 
\begin{equation*}
E_{L}^{c} = \{L^{a_{1}}C^{j-a_{1}},
\ldots , L^{a_{i}}C^{j-a_{i}}, \ldots
,L^{a_{t}}C^{j-a_{t}}\},\ a_{1}< \ldots <a_{t}.
\end{equation*} 
We define the partition $Q(V,p)\subset B(t,d), t=j+1-d$ by
\begin{equation}\label{e1.10}
{ Q}(V,p) = (a_{1}, a_{2}-1, \ldots ,
a_{i}-(i-1),
\ldots ,a_{t}-(t-1)). 
\end{equation}
 \end{definition}
%%%%%%%%%%%%%%%%%%%%%%%%%%%%%%%%%%%%%%%%%%%%%%%%%%%%%%%%%%
The partition ${ Q}(V,p)$ has $t=\cod (V)$ parts, some
of which may be zero. Each part is no greater then
$j+1-t = \dim (V)$ (Lemma \ref{lemma1.8}).

We now give an equivalent definition of this partition.
If $E_{L}$ is a set of $d$ degree$-j$ monomials in
$(L,C)$, and $E_{L}^{c}$ is its complement, we let
${\cal S}(E_{L})$ denote the following set of ordered pairs
of degree$-j$ monomials in $L, C$:
$${\cal S}(E_{L}) = \{(\mu,\nu), \mu \in E_{L}, \nu \in
E_{L}^{c}, \ and \ \mu <\nu \}. \eqno{(1.11)}$$
Given a monomial $\mu \in E_{L}^{c}$, we define a subset
of degree$-j$ monomials 
$$ {\cal S}_{E_{L}}(\nu) = \{ \mu \in E_{L} \ | \ \mu <
\nu\}. \eqno{(1.12)}$$
If $\mathrm P$ is a partition, we let $\ell({\mathrm P})$ be
the number it partitions.
%%%%%%%%%%%%%%%%%%%%%%%%%%%%%%%%%%%%%%%%%%%%%%%%%%%%%%%%%%
%%%%%%%%%%%%%%%%%%%%%%%%%%%%%%%%%%%%%%%%%%%%%%%%%%%%%%%%%
%\newtheorem{lemm18}[exam]{Lemma}
\begin{lemma}\label{lemma1.8}
  Let
$\In_{p}(V) = E_{L}$. There is one part $q_{\nu}$ of the
partition ${ Q}(V,p)$ corresponding to each of the
$t$ cobasis monomials $\nu \in E_{L}^{c}$. The part
$q_{\nu}$ is \  $\#{\cal S}(E_{L}(\nu))$ and satisfies $q_\nu\le j+1-t$, thus we have
$Q(V,p)\subset B(\cod V, \dim V)$. Each partition
$ Q\subset B(t,j+1-t)$ occurs as the
partition ${ Q}(E_{L},p)$ for a unique monomial
vector space $E_{L}$.
\end{lemma}
\begin{proof}
If $\nu = L^{a_{i}}C^{j-a_{i}} \in
E_{L}^{c}$, then $i-1 \leq a_{i} \leq j-(t-1)$, to leave
room for the other monomials $\nu' \in E_{L}^{c}$ before
and after $\nu$. Also,
$$ {\cal S}_{E_{L}}(\nu) = \{ C^{j}, LC^{j-1}, \ldots ,
L^{a_{i}}C^{j-a_{i}}\} - \{ L^{a_{1}}C^{j-a_{1}},
\ldots , L^{a_{i}}C^{j-a_{i}}\}, $$
has $(a_{i}+1)-(i) = a_{i}-(i-1)$ elements. This is the
$i-$th part $q_{i}$ of ${ Q}(E,p)$, and $q_{i}$ thus
satisfies $0\leq q_{i} \leq j+1-t$. Conversely, given
$ Q$ let $ a_{i} = q_{i} + (i+1)$. Choosing
monomials $\nu_{i} = L^{a_{i}}C^{j-a_{i}}$ we obtain a
set $E_{L}({ Q})$ of $t$ distinct degree$-j$ monomials
such that ${ Q}(E_{L},p) = { Q}$. The space
$E_{L}$ is easily seen to be independent of the choice
of $C$, so is unique.
\end{proof}
\par\smallskip\noindent
{\em Notation}
\ Let
${\mathrm P}$ be a subpartition of the retangular partition $\mathrm{B}(t,j+1-t)$ with
$t$ parts of size $j+1-t$, so $P$ has $t$ parts (some possibly zero), each of size no
greater than
$(j+1-t)$. We denote by
${\mathrm P}^{c}$ its complement in the rectangular partition
$\mathrm{B}(t,j+1-t)$, and we
denote by ${\mathrm P}^{\wedge}\subset B(j+1-t,t)$ its dual partition
(switch rows and columns in the Ferrers diagram).
%%%%%%%%%%%%%%%%%%%%%%%%%%%%%%%%%%%%%%%%%%%%%%%%%%%%%%%%%
%\newtheorem{lemm19}[exam]{Lemma}
\begin{lemma}[Ramification partition $\QRAM_p(V)$ and  ${
Q}(V,p)$]\label{lemma1.9} If $p$ is any point of \ 
${\check{\mathbb P}}^1$, the partition ${ Q}(V,p)\subset B(t,j+1-t)$ is related
to $\QRAM _{p}(V)\subset B(j+1-t,t)$ by 
\begin{equation}\label{e1.12}
\QRAM _{p}(V) = ({ Q}(V,p)^{c})^{\wedge}.
\end{equation}
Each partition having $\dim (V)$ parts, each no
greater than $cod(V)$ occurs as $\QRAM _{p}(E_{L})$
for a unique monomial vector space $E_{L}$.
\end{lemma}
\begin{proof} If $\In_{p}(V) = E_{L} = \langle 
L^{n_{1}}C^{j-n_{1}}, \ldots , L^{n_{i}}C^{j-n_{i}},
\ldots , L^{n_{d}}C^{j-n_{d}}\rangle $, then by Definition \ref{def1.2},
${ Q}(V) = (n_{1}, \ldots , n_{i}-(i-1), \ldots ,
n_{d}-(d-1))$. The partition ${ Q}(V,p)$ in
Definition \ref{def1.7} is constructed similarly from the
complementary set of monomials $E_{L}^{c}$. The formula
\eqref{e1.12} thus reduces to an easily shown combinatorial
identity. The last statement is a consequence of
Lemma \ref{lemma1.8}. 
\end{proof} 
%%%%%%%%%%%%%%%%%%%%%%%%%%%%%%%%%%%%%%%%%%%%%%%%%%%%%%%%%
%\newtheorem{exam110}[exam]{Example}
\begin{example}\label{example110}
Let $V=(x^{4}, x^{3}y, x(x+y)^{3})$ and the point
$p: x=0$; then $n_{p}(V) = (4,3,1)$ so $\QRAM _{p}(V) = (4,3,1) -
(2,1,0) = (2,2,1)$ . We have $E = \In (V) = (xy^{3}, x^{3}y, x^{4})$ so
$E^{c} = (y^{4}, y^{2}x^{2})$ and ${ Q}(V,p) = (1,0) \subset B(3,3)$.
The complementary partition ${ Q}(V,p)^{c}$ \  in
$B(3,3)$ is $(3,2)$, whose dual partition is $(2,2,1)$;
thus $\QRAM (V,p) = ({ Q}(V,p)^{c})^{\wedge}$ (Lemma \ref{lemma1.9}).
\end{example}
%%%%%%%%%%%%%%%%%%%%%%%%%%%%%%%%%%%%%%%%%%%%%%%%%%%%%%%%%
%\newtheorem{def111}[exam]{Definition}
\begin{definition}[The cells ${\mathbb V}(E,p)$ in $\Grass
(d,R_{j})$.]\label{def1.11}
 If $E$ is a $d$-dimensional vector space
spanned by monomials, and $p\in {\check{\mathbb P}}^1$ is the point
$p: x=0$, then we denote by
${\mathbb V}(E)$ or ${\mathbb V}(E,p)$ the subfamily of the
Grassmannian $\Grass (d,R_{j})$
\begin{equation}\label{e1.14}
 {\mathbb V}(E) = \{ V \ | 
\ \In (V) = \In (E)\}. 
\end{equation}
\end{definition}
%%%%%%%%%%%%%%%%%%%%%%%%%%%%%%%%%%%%%%%%%%%%%%%%%%%%%%%%%
Let $p\in {\check{\mathbb P}}^1$ be arbitrary, and let
$E=\{x^{n_{1}}y^{j-n_{1}}, \ldots ,
x^{n_{d}}y^{j-n_{d}}\}, \ n_{1}< \ldots <n_{d}$ be a set
of $d$ monomials in $x, y$. We denote by  $E_{L} =
\{L^{n_{1}}C^{j-n_{1}}, \ldots , L^{n_{d}}C^{j-n_{d}}\}$,
the corresponding set of monomials in any basis $(L,C)$ for $R$,
$L=L_{p}$, $C\in R_{1}$ with $\langle C\rangle \neq \langle L\rangle $. We denote by
${\mathbb V}(E,p)$ the corresponding family to ${\mathbb V}(E)$ in
the basis $(L,C)$ for $R$ : ${\mathbb V}(E,p)$ is the set of
vector spaces whose initial forms in the basis $C^{j}<
C^{j-1}L< \ldots < L^{j}$ for $R_{j}$ is $E_{L}$. We have 
\begin{equation}\label{1.15}
{\mathbb V}(E,p) = \{ V\subset R_{j} \ | \ \dim_{k}V
=d,
\ and \ \QRAM (V,p) = \QRAM (E) \}.
\end{equation} 
We give ${\mathbb V}(E,p)$ the reduced subscheme structure
inherited from $\Grass(d,R_{j})$.
We will show that the closure $\overline{{\mathbb V}(E,p)}$ is a Schubert cell of
the Grassmannian $\Grass(d,R_{j})$ below in Lemma \ref{lemma2.1B}.
%%%%%%%%%%%%%%%%%%%%%%%%%%%%%%%%%%%%%%%%%%%%%%%%%%%%%%%%%
%%%%%%%%%%%%%%%%%%%%%%%%%%%%%%%%%%%%%%%%%%%%%%%%%%%%%%%%%%
%%%%%%%%%%%%%%%%%% Sections 2-3 %%%%%%%%%%%%%%%%%%%%%%%%%%
%%%%%%%%%%%%%%%%%%%%%%%%%%%%%%%%%%%%%%%%%%%%%%%%%%%%%%%%%%
%%%%%%%%%%%%%%%%%%%%%%%%%%%%%%%%%%%%%%%%%%%%%%%%%%%%%%%%%%

\subsection{Ramfication and Schubert cells of the Grassmanian
$\Grass(d,R_j)$}\label{Grasscells}
%%%%%%%%%%%%%%%%%%%%%%%%%%%%%%%%%%%%%%%%%%%%%%%%%%%%%%%%%
Let $E$ be a subspace of $R_{j}$ spanned by monomials in
$x, y$; we let $E^{c}$ denote the degree$-j$ monomials not
in $E$, which we term {\em cobasis} monomials.
%%%%%%%%%%%%%%%%%%%%%%%%%%%%%%%%%%%%%%%%%%%%%%%%%%%%%%%%%
For $p\in {\check{\mathbb P}}^1$ we denote by ${\cal F}_{j}(p)$
the flag 
\begin{equation}
{\cal F}_{j}(p):\, 0=F_{0}(p)\subset F_{1}(p)\subset \cdots \subset
F_{i}(p)\subset \cdots \subset F_{j+1}(p)=R_{j}
\end{equation}  
of subspaces of $R_{j}$,
where  $C\in R_{1}$ satisfies $\langle C\rangle \neq
\langle L_{p}\rangle $ and where 
\begin{equation}\label{e2.1}
F_{i}(p)=\langle L_{p}^{j},
L_{p}^{j-1}C,\ldots , L_{p}^{j+1-i}C^{i-1}\rangle .
\end{equation}
\begin{lemma}\label{lemma2.1A} Let $E=\{x^{n_{1}}y^{j-n_{1}},\ldots
		,x^{n_{d}}y^{j-n_{d}}\},\ \  n_{1}<n_{2}<\cdots
		<n_{d}$, and consider a point $p$ of ${\check{\mathbb P}}^1$. The following are equivalent:
		\begin{flalign}
		i.\, &  V\in {\mathbb V}(E,p), \text { or,
equivalently, } Q(V,p)=Q(E).\notag\\
		ii.\, & \QRAM(V,p) = \QRAM (E).&\label{e2.2}\\ 
		iii.\, &\dim_{k}(V\cap F_{j+1-n_{i}}(p))=d+1-i, \ \ 
		i=1,\ldots , d. &\notag
		\end{flalign}
\end{lemma}
%%%%%%%%%%%%%%%%%%%%%%%%%%%%%%%%%%%%%%%%%%%%%%%%%%%%%%%%%%
\begin{lemma}\label{lemma2.1B} For any $p\in {\check{\mathbb P}}^1$,
the subvariety ${\mathbb V}(E,p)$ is an open dense subset
of a Schubert cell $\overline{{\mathbb V}(E,p)}$ on
$\Grass(d,R_{j})$ and it is an affine space. The
dimension of ${\mathbb V}(E,p)$ is
\begin{equation}\label{e2.3}
\dim ({\mathbb V}(E,p)) = \#{\cal S}(E) = N-r(E) =
\ell({ Q}(E)).  
\end{equation}
	where $r(E)=\ell(\QRAM (E))$ of \eqref{etotalram} and it has  codimension  in
$\Grass(d,R_{j})$
\begin{equation}\label{e2.4}
\cod ({\mathbb V}(E,p))=r(E) = \ell (\QRAM(E))=(\sum_{1\leq i \leq d}n_{i}) - d(d-1)/2.
\end{equation}
\end{lemma}                               
%%%%%%%%%%%%%%%%%%%%%%%%%%%%%%%%%%%%%%%%%%%%%%%%%%%%%%%%%%
\begin{proof}[Proof of Lemma \ref{lemma2.1A}] It suffices to
consider $p: x=0$.  We have $V\in {\mathbb V}(E)$ iff $\In (V)
= E$. Let $(f_{1},\ldots ,f_{d})$ be a basis of $V$ such
that $\In (f_{i})=x^{n_{i}}y^{j-n_{i}}\in E$, for each $i$.
%%%%%%%%%%%%%%%%%%%%%%%%%%%%%%%%%%%%%%%%%%%%%%%%%%%%%%%%%%
Then  $n_{1}<n_{2}<\cdots < n_{d}$, $x^{n_{i}}|f_{i}$, and
it is easy to see that there is no basis $F'$ of $V$ with
larger $x$\ -powers than $n_{i}$.
%%%%%%%%%%%%%%%%%%%%%%%%%%%%%%%%%%%%%%%%%%%%%%%%%%%%%%%%%%
By Definition \ref{def1.2},  $\QRAM(V,p)=\QRAM(E)$.
From $\In (V) = E$ we see that \eqref{e2.2} is satisfied. Thus
$i. \Longrightarrow ii. \Longrightarrow iv.$, and the
converse implications are an easy consequence. Lemma \ref{lemma1.9}
shows $ii. \Longleftrightarrow iii.$ 
\end{proof}
\begin{proof}[Proof of Lemma \ref{lemma2.1B}]
The identification of ${\mathbb
V}(E)$  with a Schubert  cell is immediate from \eqref{e2.2}. That ${\mathbb
V}(E)$ is an affine space, and has the codimension $r(E)$ in
$\Grass(d,R_{j})$ is a standard result (see \S I.5 of
[GrH], or see the proof below of Theorem \ref{thm4.9} in the special case $T=(1,2,\ldots
,j,t,0)$). The formula
\eqref{e2.3} follows from
\eqref{e2.4} and
 Lemma \ref{lemma1.9}. 
\end{proof}\par
\noindent If $G$ is a complex algebraic variety, we grade
the homology ring $H^{\ast}(G,{\mathbb Z})$ by codimension,
and consider the intersection product 
$H^{u}(G,{\mathbb Z})\times H^{v}(G,{\mathbb Z})
\longrightarrow H^{u+v}(G,{\mathbb Z})$. Since our
varieties are connected, if $N= \dim(G)=d(j+1-d)$, $H^{N}(G) \cong
{\mathbb Z}\cdot\zeta_{0}$ and has basis the class $\zeta_{0}$
of a point. If \  $a\zeta_{0} \in H^{N}(G)$ we will denote
by $[a\zeta_{0}]$ the integer $a$.\\
%%%%%%%%%%%%%%%%%%%%%%%%%%%%%%%%%%%%%%%%%%%%%%%%%%%%%%%%%%
%%%%%%%%%%%%%%%%%%%%%%%%%%%%%%%%%%%%%%%%%%%%%%%%%%%%%%%%%%
\begin{proposition}\label{prop2.2} Assume ${k} = {\mathbb
C}$. The Wronskian morphism $w: \Grass(d,R_{j})
\longrightarrow {{\mathbb P}}^{N}$ is a finite cover of degree 
$[c_{1}^{N}] = N!/\binom{j}{d}$ where $c_{1}$ is the class
of a hyperplane in the homology
$H^{\ast}(\Grass(d,R_{j}))$. The degree of $w$ is the same
as that of $\Grass (d,R_{j})$ in the Pl\"ucker embedding.
\end{proposition}
\begin{proof} Clearly $w$ is a proper algebraic
morphism. We let $q\in {{\mathbb P}}^{N} = \Sym^{N}({\check{\mathbb
P}}^1)$ correspond to a set $\{q\} = p_{1}, \ldots ,p_{N}$ of $N$
distinct points of ${\check{\mathbb P}}^1$. For each $i$, The
condition that $\QRAM (V,p_{i}) = (0, \ldots ,0,1), 1\leq i \leq N$
is by Lemma \ref{lemma2.1A} a simplest Schubert condition on $\Grass(d,R_{j})$,
whose class is $c_{1}$ in the homology ring $H^{\ast}(\G_T)$.
Since the points $p_{1}, \ldots , p_{N}$ are distinct the
Schubert conditions are distinct. The intersection of all 
$N$ conditions will consist of $m=[c_{1}^{N}]$ points of
$\Grass(d,R_{j})$, provided that the intersection is
proper. Should the intersection not be proper, and there
is a one dimensional family of vector spaces in the
intersection, then we would impose an additional
condition $\QRAM (V,p_{N+1}) =0$ at any further
point $p_{N+1} \in {\check{\mathbb P}}^1$, and find a
vector space having $N+1$ ramification points,
contradicting Lemma \ref{lemma1.4}. Thus, there will be
$[c_{1}^{N}] = N!/\binom{j}{d}$ vector spaces $V$ having
total ramification one at each of the points $p_{i}$. By
\eqref{e1.9} any such vector space $V$ can
have no further ramification, and $w(V) = q$. We have
shown that $w: w^{-1}({{\mathbb P}}^{N} \setminus \Delta)
\longrightarrow {{\mathbb P}}^{N} \setminus \Delta$ is an
$m$-to-one cover, where $\Delta$ is the large diagonal;
since $\Grass(d,R_{j})$ is irreducible, $w$ is a finite
cover as claimed. That the integer $[c_{1}^{N}]$ is the
the degree of $\Grass(d,R_{j})$ in its Pl\"ucker
embedding is well known.  
\end{proof}
\begin{definition}\label{Definition 2.3} Recall that the homology
class of the intersection 
${\mathbb V}_{{ Q}_{1}}\cap {\mathbb V}_{{ Q}_{2}}$ of
two Schubert subvarieties of $\Grass (d,R_{j})$ contains a
term corresponding to ${Q}_{a}$ iff there are
permutations $\tau$ and $\sigma$ of $(1, \ldots , j+1-d)$
such that we can write each part $q_{a_{i}}$ of  
${Q}_{a}$ as $q_{a_{i}} = q_{1\sigma(i)}
+q_{2\tau(i)}$, the sum of parts in ${ Q}_{1}$ and 
${ Q}_{2}$. We say ${ Q}_{a}\subset { Q}_{1}
+{ Q}_{2}$, in this case, and extend this definition
to sums of $s$ partitions. 
\end{definition}
\begin{proposition}\label{prop2.4} Let $p_{1}, \ldots ,p_{s}$
be distinct points of ${\check{\mathbb P}}^1$, and let
$E_{1}, \ldots , E_{s}$ are monomial vector spaces of
dimension $d$ in $R_{j}$. Assume that $\sum \ell(Q(E_i))\le N=d(j+1-d)$. Then the
homology class of the intersection 
$\overline{{\mathbb V}(E_{1},p_{1})}\cap \ldots \cap 
\overline{{\mathbb V}(E_{s},p_{s})}$ of cells satisfies 
\begin{equation}\label{e2.5}
[{\mathbb V}(E_{1},p_{1})\cap \ldots \cap 
{\mathbb V}(E_{s},p_{s})] = \sum n_{a}[{ Q}_{a}], \ 
{Q}_{a} \subset { Q}(E_{1}) + \ldots + 
{Q}(E_{s}). 
\end{equation} 
Here the coefficients $n_{a}$ are given by the Schubert
calculus. If the above class in \eqref{e2.5} is nonzero, the codimension of the
intersection is 
\begin{equation}\label{e2.6}
\sum r(E_{i}) = N - \sum \ell({Q}(E_{i})).
\end{equation}
\end{proposition}
\begin{proof} Immediate from Lemmas \ref{lemma2.1A} and \ref{lemma2.1B}
and the Schubert calculus, provided that the intersection
$Y = \overline{{\mathbb V}(E_{1},p_{1})}\cap \ldots \cap 
\overline{{\mathbb V}(E_{s},p_{s})}$ is proper and nonempty.
The condition \eqref{e2.5} for nonemptiness is a consequence of
the Schubert calculus. It remains to show that the
intersection is proper if nonempty. {\sc wolog} we may add
codimension 1 conditions $\QRAM (V,p_{i}) = (0,
\ldots , 0,1)$ at a finite number of points $p_{s+1},
\ldots , p_{s'}$ so that $s'-s = N - \sum r(E_{i})$,
obtaining $Y' = Y\cap H_{s+1}\cap \ldots \cap H_{s'}$. Then the Wronskian of $V\in
Y'$ satisfies: 
$W(V)=\prod_{_{1\leq i \leq s}}
L_{p_{i}}^{r(E_{i})}\prod_{_{s\leq i\leq s'}} L_{P{i}}$.
Since $w: \Grass(d,R_{j}) \longrightarrow {{\mathbb P}}^{N}$ is
a finite cover, there are only a finite number of vector spaces
$V$ having a given Wronskian, so $Y'$ is a finite set of
points. It follows that $Y$ is proper if nonempty.
\end{proof}
\begin{figure}[hbtp]
\begin{center}
\leavevmode %(80,100)(96,-100);(110,130)(85,-130)
\begin{picture}(70,90)(40,-110)
\setlength{\unitlength}{1.4mm}
\put(30,-4){{\small \mbox{${\Bbb V}(E)$ Cell}}}
\put(84,-4){{\small \mbox{${ Q}(E)$ Code}}}
\put(39,-16){\makebox(0,0){{\small \mbox{$A$}}}}
\put(39,-40){\makebox(0,0){{\small \mbox{$B$}}}}
\put(39,-64){\makebox(0,0){{\small \mbox{${ Q}(C)$}}}}
\put(39,-88){\makebox(0,0){{\small \mbox{$F$}}}}
\put(39,-112){\makebox(0,0){{\small \mbox{$G$}}}}
\put(9,-64){\makebox(0,0){{\small \mbox{$C$}}}}
\put(63,-64){\makebox(0,0){{\small \mbox{$D$}}}}
\put(90,-124){\makebox(0,0){{\small \mbox{$\emptyset$}}}}
\put(43,-24){{\scriptsize \mbox{$y^{2}+a_{1}yx^{2}+b_{1}x^{3}$}}}
\put(37,-28){{\scriptsize \mbox{$y^{3}+a_{2}yx^{2}+b_{2}x^{3}$}}}
\put(49,-44){{\scriptsize \mbox{$yx^{2}+b_{1}x^{3}$}}}
\put(37,-52){{\scriptsize \mbox{$y^{3}+a_{1}y^{2}x+b_{2}x^{3}$}}}
\put(25,-64){{\scriptsize \mbox{$x^{3}$}}}
\put(7,-76){{\scriptsize \mbox{$y^{3}+a_{1}y^{2}x+b_{1}yx^{2}$}}}
\put(73,-68){{\scriptsize \mbox{$yx^{2}+a_{1}x^{3}$}}}
\put(67,-72){{\scriptsize \mbox{$y^{2}x+a_{2}x^{3}$}}}
\put(55,-88){{\scriptsize \mbox{$x^{3}$}}}
\put(43,-96){{\scriptsize \mbox{$y^{2}x+a_{1}yx^{2}$}}}
\put(55,-112){{\scriptsize \mbox{$x^{3}$}}}
\put(49,-116){{\scriptsize \mbox{$yx^{2}$}}}
\thicklines
\multiput(0,-32)(0,-24){4}{\line(1,0){96}}
\put(48,-80){\line(0,1){24}}
\thinlines
\multiput(30,-28)(6,4){2}{\line(1,0){6}}
\multiput(30,-24)(6,4){2}{\line(1,0){6}}
\multiput(30,-28)(6,4){2}{\line(0,1){4}}
\multiput(36,-28)(6,4){2}{\line(0,1){4}}
\multiput(60,-72)(6,4){2}{\line(1,0){6}}
\multiput(60,-68)(6,4){2}{\line(1,0){6}}
\multiput(60,-72)(6,4){2}{\line(0,1){4}}
\multiput(66,-72)(6,4){2}{\line(0,1){4}}
\multiput(42,-116)(6,4){2}{\line(1,0){6}}
\multiput(42,-112)(6,4){2}{\line(1,0){6}}
\multiput(42,-116)(6,4){2}{\line(0,1){4}}
\multiput(48,-116)(6,4){2}{\line(0,1){4}}
\multiput(30,-52)(12,8){2}{\line(1,0){6}}
\multiput(30,-48)(12,8){2}{\line(1,0){6}}
\multiput(30,-52)(12,8){2}{\line(0,1){4}}
\multiput(36,-52)(12,8){2}{\line(0,1){4}}
\multiput(36,-96)(12,8){2}{\line(1,0){6}}
\multiput(36,-92)(12,8){2}{\line(1,0){6}}
\multiput(36,-96)(12,8){2}{\line(0,1){4}}
\multiput(42,-96)(12,8){2}{\line(0,1){4}}
\multiput(0,-76)(18,12){2}{\line(1,0){6}}
\multiput(0,-72)(18,12){2}{\line(1,0){6}}
\multiput(0,-76)(18,12){2}{\line(0,1){4}}
\multiput(6,-76)(18,12){2}{\line(0,1){4}}
\multiput(42,-20)(6,4){2}{\rule{6mm}{4mm}}
\multiput(6,-72)(6,4){2}{\rule{6mm}{4mm}}
\multiput(30,-124)(6,4){2}{\rule{6mm}{4mm}}
\multiput(36,-48)(12,8){2}{\rule{6mm}{4mm}}
\multiput(30,-100)(12,8){2}{\rule{6mm}{4mm}}
\multiput(54,-76)(18,12){2}{\rule{6mm}{4mm}}
\multiput(84.5,-24)(.5,0){11}{\line(0,1){3.5}}
\multiput(84.5,-28)(.5,0){11}{\line(0,1){3.5}}
\multiput(90.5,-24)(.5,0){11}{\line(0,1){3.5}}
\multiput(90.5,-28)(.5,0){11}{\line(0,1){3.5}}
\multiput(36.5,-72)(.5,0){11}{\line(0,1){3.5}}
\multiput(36.5,-76)(.5,0){11}{\line(0,1){3.5}}
\multiput(84.5,-76)(.5,0){11}{\line(0,1){3.5}}
\multiput(90.5,-76)(.5,0){11}{\line(0,1){3.5}}
\multiput(84.5,-100)(.5,0){11}{\line(0,1){3.5}}
\multiput(84.5,-48)(.5,0){11}{\line(0,1){3.5}}
\multiput(90.5,-48)(.5,0){11}{\line(0,1){3.5}}
\multiput(84.5,-52)(.5,0){11}{\line(0,1){3.5}}
\thicklines
\multiput(84,-20)(0,-4){3}{\line(1,0){12}}
\multiput(84,-44)(0,-4){2}{\line(1,0){12}}
\multiput(84,-72)(0,-4){2}{\line(1,0){12}}
\multiput(84,-96)(0,-4){2}{\line(1,0){6}}
\multiput(36,-68)(0,-4){3}{\line(1,0){6}}
\put(84,-52){\line(1,0){6}}
\multiput(84,-28)(6,0){3}{\line(0,1){8}}
\multiput(84,-52)(6,0){2}{\line(0,1){8}}
\multiput(36,-76)(6,0){2}{\line(0,1){8}}
\put(96,-48){\line(0,1){4}}
\multiput(84,-76)(6,0){3}{\line(0,1){4}}
\multiput(84,-100)(6,0){2}{\line(0,1){4}}
\end{picture}
\end{center}
\vspace{8cm}
\renewcommand{\figurename}{Table}
\protect\caption{The cells ${\protect\Bbb V}(E)$ of $\Grass(2,4)$
and their codings ${Q}(E)$. \label{tab-cell}}
\end{figure}
%\newpage
%%%%%%%%%%%%%%%%%%%%%%%%%%%%%%%%%%%%%%%%%%%%%%%%%%%%%%%%%%%
%\begin{center} \underline{{\bf Table 1.}} \ \  {\bf
% The cells} ${\Bbb V}(E)$ {\bf of} $Grass(2,4)$ {\bf and
%their codings} ${Q}(E)$. \end{center}
%%%%%%%%%%%%%%%%%%%%%%%%%%%%%%%%%%%%%%%%%%%%%%%%%%%%%%%%%%%
%
%\vspace{12cm}
%
%\special{picture Gr2/table(1) scaled 800}
%
%\vspace{5mm}
%
%%%%%%%%%%%%%%%%%%%%%%%%%%%%%%%%%%%%%%%%%%%%%%%%%%%%%%%%%%%
%%%%%%%%%%%%%%%%%%%%%%%%%%%%%%%%%%%%%%%%%%%%%%%%%%%%%%%%%%%
%\hspace{.7cm}\\
%%%%%%%%%%%%%%%%%%%%%%%%%%%%%%%%%%%%%%%%%%%%%%%%%%%%%%%%%%
\par\noindent{\bf Table \ref{tab-cell}} gives the six cells of $\Grass (2,R_{3})$, for
the  $p: x=0$ direction. At left the cobasis monomials 
for the cell ${\mathbb V}(E)$ are shaded; the initial
monomials of generators of $V$ in the standard basis
are unshaded. We give at right the  partition
${Q}(E)$ (see Definition \ref{def1.7}).
%%%%%%%%%%%%%%%%%%%%%%%%%%%%%%%%%%%%%%%%%%%%%%%%%%%%%%%%%%

The codimension of a cell $\mathbb{V}(E)$ is its depth
below  the top of the table.  Thus, cells $\mathbb{V}(E)$
with higher powers of $x$  dividing the cobasis $E^{c}$
have a {\em higher} shaded portion $E^{c}$ in Table \ref{tab-cell}:
these cells correspond to codings ${Q}(E)$ with
{\em larger} area $\ell({Q}(E))$ giving the
dimension of the cell, or equivalently, to vector spaces
$V\in {\mathbb V}$ with {\em smaller} ramification at $p$.
%%%%%%%%%%%%%%%%%%%%%%%%%%%%%%%%%%%%%%%%%%%%%%%%%%%%%%%%%%

%%%%%%%%%%%%%%%%%%%%%%%%%%%%%%%%%%%%%%%%%%%%%%%%%%%%%%%%%%
\begin{example}[Coding of a cell]\label{ex2.5} The
partition ${Q}(B)$ of Table \ref{tab-cell} for the second cell 
$B$ has two rows. The length, two, of the first row
$Q_{\nu_{1}}$, $\nu_{1} =x ^{3}$ counts the two pairs of
monomials $(yx^{2},x^{3})$ and $(y^{3},x^{3})$;
equivalently, the length counts the two coefficients
$b_{1}, b_{2}$ of the basis of $V$ on $\nu = x^{3}$. The
length, one, of the second row $Q_{\nu_{2}}$,
$\nu_{2}=y^{2}x$ counts the single pair $(y^{3},y^{2}x)$;
equivalently, it counts the single coefficient $a_{2}$ of
$y^{3}+a_{1}y^{2}x+b_{2}$ on $y^{2}x$.
\end{example}
A given vector space $V$ has a coding ${Q}(V,p)$ for each point $p\in {\check{\mathbb P}}^1$. By
Lemmas \ref{lemma1.4}, and \ref{lemma2.1A}, for all but a
finite number of points $p$ the partition ${Q}(V,p)$ is
$B(j+1-d,d)$ and $\QRAM (V,p) = 0 $: we say that $V$ is in the
{\em generic cell} for such points $p$ in which the
partition ${Q}(V,p)$ 
is $B(j+1-d,d)$. We set $\cod({Q}(V,p)) = N
- \ell({Q}(V,p))$. By Lemma \ref{lemma1.9} $\cod({Q}(V,p))=\ell(\QRAM(V,p))$.
Equations \eqref{e1.4} and \eqref{e1.9} imply 
\begin{lemma}\label{lemma2.6} We have
\begin{equation}\label{e2.7}
\sum_{p\in {\check{\mathbb P}}^1} \cod({Q}(V,p)) = N =
\cod(V)\cdot \dim(V). 
\end{equation}
\end{lemma}
\begin{example}\label{ex2.7} In Table \ref{tab-cell} the vector
space 
$V_1= \langle yx^{2},x^{3}\rangle $
is the smallest cell in the
$p:x=0$ direction and has codimension four. By the sum formula
\eqref{e2.7},
$V_1$ lies in the generic cell in any other direction. 
 The vector space $V_2 = \langle y^{2}x+x^{3}, yx^{2}\rangle $
lies in the cell of codimension two at the right of
			the table, so must lie either in one cell of codimension
			two in some other direction , or in
			two cells of codimension one in two different
			directions (here, the latter, in directions $y-x, \ y+x$).
\end{example}
\subsection{Classical results for ramification conditions in $\Grass
(d,R_{j})$}\label{classicalresults}
%%%%%%%%%%%%%%%%%%%%%%%%%%%%%%%%%%%%%%%%%%%%%%%%%%%%%%%%%%
We state known results for the Schubert cells ${\mathbb
V}(E,p)$ of the Grassmannian variety $\Grass (d,R_{j})$.
Each cell consist of all vector spaces $V$ having given
ramification $\QRAM _{p}(V) = \QRAM (E)$ at a
point $p$ of ${\check{\mathbb P}}^1$; the cell also consists
of all the vector spaces $V$ having given initial space
$\In (V) = E_{L}$ in a basis $(L_{p},C)$ of $R$ (Lemma \ref{lemma1.9}). This
approach to ramification of linear systems over $\mathbb P^1$ by Schubert calculus
has been extended to families of linear systems on curves in \cite{EH1,EH2}). Our
approach here emphasizes the connections with combinatorics, in preparation for
considering ideals of linear systems in Section \ref{sec4}.\\
%%%%%%%%%%%%%%%%%%%%%%%%%%%%%%%%%%%%%%%%%%%%%%%%%%%%%%%%%%
%%%%%%%%%%%%%%%%%%%%%%%%%%%%%%%%%%%%%%%%%%%%%%%%%%%%%%%%%% 

We recall our notation.
%%%%%%%%%%%%%%%%%%%%%%%%%%%%%%%%%%%%%%%%%%%%%%%%%%%%%%%%%%
We fix \ $d=\dim (V)$, $j=\mathrm{degree}(V)$, so $\cod(V)=j+1-d$,
and
%%%%%%%%%%%%%%%%%%%%%%%%%%%%%%%%%%%%%%%%%%%%%%%%%%%%%%%%%%
we let $N=d(j+1-d)$, the dimension of $\Grass (d,R_{j})$.
%%%%%%%%%%%%%%%%%%%%%%%%%%%%%%%%%%%%%%%%%%%%%%%%%%%%%%%%%%
We let $\cod (\mathbb{V}(E)) = N-\dim(\mathbb{V}(E))$ be the
codimension of the cell $\mathbb{V}(E)$ in the Grassmannian 
$\Grass (d,R_{j})$. 
%%%%%%%%%%%%%%%%%%%%%%%%%%%%%%%%%%%%%%%%%%%%%%%%%%%%%%%%%%
The partitions ${Q}$ we consider will have $j+1-d$
parts  $0\leq q_{1}\leq \cdots \leq q_{j+1-d}\leq d$.
%%%%%%%%%%%%%%%%%%%%%%%%%%%%%%%%%%%%%%%%%%%%%%%%%%%%%%%%%%
We partially order the partitions by ${Q}\leq
{Q}'$ iff $q_{i}\leq q'_{i}$ for each $i$. 
%%%%%%%%%%%%%%%%%%%%%%%%%%%%%%%%%%%%%%%%%%%%%%%%%%%%%%%%%% 
We denote by ${Q}^{\wedge}$ the dual partition
obtained by exchanging rows and columns in the shape of
${Q}$:
%%%%%%%%%%%%%%%%%%%%%%%%%%%%%%%%%%%%%%%%%%%%%%%%%%%%%%%%%%
it has $d$ parts of sizes between zero and $j+1-d$.
%%%%%%%%%%%%%%%%%%%%%%%%%%%%%%%%%%%%%%%%%%%%%%%%%%%%%%%%%% 
The shape of ${Q}$ is included in a $\cod (V)\times
\dim (V)$ rectangle ${B}$;       
%%%%%%%%%%%%%%%%%%%%%%%%%%%%%%%%%%%%%%%%%%%%%%%%%%%%%%%%%%     
we denote by ${ Q}^{c}$ its complement in       
${B}$.
%%%%%%%%%%%%%%%%%%%%%%%%%%%%%%%%%%%%%%%%%%%%%%%%%%%%%%%%%%
If  $E=\{x^{n_{1}}y^{j-n_{1}},\ldots
		,x^{n_{d}}y^{j-n_{d}}\}$ we denote by $E^{\wedge}$    
		the space of monomials\\
		$E^{\wedge} = \{y^{n_{1}}x^{j-n_{1}},\ldots
		, y^{n_{d}}x^{j-n_{d}}\}$. 
%%%%%%%%%%%%%%%%%%%%%%%%%%%%%%%%%%%%%%%%%%%%%%%%%%%%%%%%%%
The ``length"  $\ell({Q})$  of a partition
$Q$ is  $\ell({Q}) = \sum q_{i}$, the        
``area" of its Ferrers graph, or shape.\\
%%%%%%%%%%%%%%%%%%%%%%%%%%%%%%%%%%%%%%%%%%%%%%%%%%%%%%%%%%

The result Cii. below follows from D.~Eisenbud and
J.~Harris in [E-H2]; the rest are consequences of Lemmas
\ref{lemma2.1A},\ref{lemma2.1B}, Propositions \ref{prop2.2} and \ref{prop2.4},
and the Schubert calculus (see \cite[Appendix 3]{I2}, or \cite{ste,EH1}).\par\medskip
%%%%%%%%%%%%%%%%%%%%%%%%%%%%%%%%%%%%%%%%%%%%%%%%%%%%%%%%%

\par {\bf A.} The cells ${\mathbb V}(E)$ are affine
spaces.\\
%%%%%%%%%%%%%%%%%%%%%%%%%%%%%%%%%%%%%%%%%%%%%%%%%%%%%%%%%

\par {\bf B.}	The set of Schubert cells
 ${\mathbb V}(E)$ correspond 1-1 with the set of 
			ordinary partitions ${Q}(E)$ of integers into no
			more than $\cod(V)$ parts, each of size no greater
			than $\dim (V)$.\\ 

			${\mathbb V}(E) \longrightarrow {Q}(E)$, a
			subpartition of a  $\cod (V)\times \dim (V)$ 
			rectangle.\\
		\indent ${Q}(E) = \{ q_{\nu}=\#S(\nu), \nu \in
		E^{c}\}, \ S(\nu) = \{ monomials \  | \  \mu \in E,
\mu \leq \nu\}$  \\
The partition
${Q}(E)$ is related to the ramification
$\QRAM(E)$ of $E$ at $x$ by 
\begin{equation}\label{e3.2}
{Q}(E) = \QRAM(E)^{c}. 
\end{equation}
If
$V\in {\mathbb V}(E)$,  then the ramification  of
$V$ at $p$  satisfies 
$\QRAM(V,p)=\QRAM(E)$ 
 by Lemma \ref{lemma1.9}. The total ramification of $V$ at
$p$  is the length of $\QRAM(V,p)$ (area of its Ferrers graph)
 and satisfies 
\begin{equation}\label{e3.3}
r(V,p)=N\, - \ell({Q}(E)).
\end{equation}
%%%%%%%%%%%%%%%%%%%%%%%%%%%%%%%%%%%%%%%%%%%%%%%%%%%%%%%%%%
%%%%%%%%%%%%%%%%%%%%%%%%%%%%%%%%%%%%%%%%%%%%%%%%%%%%%%%%%%

\par {\bf C.} The dimension of 
${\mathbb V}(E) =\ell({Q}(E))$.  And ${Q}(E^{\wedge}) = {Q}(E)^{c}.$ \\
%%%%%%%%%%%%%%%%%%%%%%%%%%%%%%%%%%%%%%%%%%%%%%%%%%%%%%%%%%
\par {\bf Ci.} The cells have
dimensionally proper intersections. The codimension of 
the intersection of ramification cells at two different
points
$x$, $L$  is the sum of their
codimensions,
\begin{equation}\label{e3.5}                    
\cod (\overline{{\mathbb V}(E,p)})\,  \cap \, 
\overline{{\mathbb V}(E',p)}) = \cod ({\mathbb V}(E,p)) + 
\cod ({\mathbb V}(E',p')),
\end{equation}
provided that the intersection of the Schubert cycles
corresponding to $E$, $E'$  are nonzero. Otherwise,
the intersection is empty.\\ 
%%%%%%%%%%%%%%%%%%%%%%%%%%%%%%%%%%%%%%%%%%%%%%%%%%%%%%%%%%

\par {\bf Cii.} {\bf Ramification behaves well
under specialization:}  If $V_{t} \  | \  t\neq 0$ 
and $W_{t} \  | \  t \neq 0, t \in {\bf A}^{1}$  are
algebraic families of vector spaces satisfying
$V_{t} \in {\mathbb V}(E,p_{t}), W_{t} \in V(E',p'_{t})$  
tending to a common limit 
$V_{0}=W_{0}$, 
if ${\mathbb V}(E,p_{t})$ and ${\mathbb V}(E',p'_{t})$ 
have dimensionally proper intersection, and
$p_{0}=p'_{0}$, then $V_{0}$  belongs to some
cell ${\mathbb V}(E'',p_{0})$ 
 where the class of $E''$
 occurs with positive coefficient in the intersection
of Schubert cycles corresponding to 
$\overline{{\mathbb V}(E,p)})\,  \cap \, 
\overline{{\mathbb V}(E',p)}$. 
(See \cite{EH2}).\\ 

%%%%%%%%%%%%%%%%%%%%%%%%%%%%%%%%%%%%%%%%%%%%%%%%%%%%%%%%%%
%%%%%%%%%%%%%%%%%%%%%%%%%%%%%%%%%%%%%%%%%%%%%%%%%%%%%%%%%%

\par {\bf D.} {\bf Exact duality:} 
if ${\mathbb V}(E)$  and ${\mathbb V}(E')$ 
have complementary dimension, then the intersection
$\overline{{\mathbb V}(E)}\, \cap \, \overline{{\mathbb
V}(E')} = {\mathbb V}(E_{0})$ 
is the class of a point iff
$E'= E^{\wedge}$; the intersection is empty if $E'\neq
E^{\wedge}$.\\
%%%%%%%%%%%%%%%%%%%%%%%%%%%%%%%%%%%%%%%%%%%%%%%%%%%%%%%%%%
%%%%%%%%%%%%%%%%%%%%%%%%%%%%%%%%%%%%%%%%%%%%%%%%%%%%%%%%%%

\par {\bf E.} {\bf Frontier property:} 
The closure $\overline{{\mathbb V}(E)}$ 
is the union of those cells  
${\mathbb V}(E')$  such that \linebreak ${Q}(E) \geq {Q}(E')$ in the sense of inclusion
of Ferrers diagrams.\\  
%%%%%%%%%%%%%%%%%%%%%%%%%%%%%%%%%%%%%%%%%%%%%%%%%%%%%%%%%%
%%%%%%%%%%%%%%%%%%%%%%%%%%%%%%%%%%%%%%%%%%%%%%%%%%%%%%%%%%

\noindent Our goal in this article is to see which of these
remarkably good properties for the ramification/Schubert
cells of $\Grass(d,R_j)$ parametrizing $d$-dimensional vector spaces
$V\subset R_j, R=k[x,y]$ extends to the analogous ramification cells
for the variety $\G_T$ parametrizing ideals $I$ of $R$ having Hilbert function
$H(R/I)=T$. This will lead us to study the homology ring of
$\G_T$, and to determine it in some special cases. \medskip
%%%%%%%%%%%%%%%%%%%%%%%%%%%%%%%%%%%%%%%%%%%%%%%%%%%%%%%%%%
%%%%%%%%%%%%%%%%%%% Section 5 %%%%%%%%%%%%%%%%%%%%%%%%%%%
%%%%%%%%%%%%%%%%%%%%%%%%%%%%%%%%%%%%%%%%%%%%%%%%%%%%%%%%%%
\section{Ramification cells in the family of graded
ideals of $R$}\label{sec4}
%%%%%%%%%%%%%%%%%%%%%%%%%%%%%%%%%%%%%%%%%%%%%%%%%%%%%%%%%%
We now show our main results, determining the homology
groups of $\G_T$, by connecting the dimension of the cells $\mathbb V(E)$ with the hook
code of the partition $P(E)$ having diagonal lengths $T$. In Section \ref{sec4A} we
extend our definitions from Section \ref{sec2} to ideals, and we show the equivalence
between suitable ramification $\QRAM(I,p)=E$ of the graded ideal $I$ at $p\in
{\check{\mathbb P}}^1$, and the ideal having initial monomial ideal $E$ (Proposition
\ref{prop4.5}). We also define the partition
$P(E)$ (Definition \ref{def4.6}). In Section~\ref{sec4B} we show that the cell 
$\mathbb V(E)$ is an affine space, with parameters given by certain coefficients
of generators, determined by the difference-one hooks of $P(E)$ (Theorem \ref{thm4.9}).
We also show that there is a birational map $\mathcal G: \G_T\to \SGrass(T)$, the
product of small Grassmanians (Proposition \ref{prop4.11}). In Section \ref{dimcell}
we determine the fibre dimension of the cell $\mathbb Z(E)$ (all ideals) over
$\mathbb V(E)$ (graded ideals) in Proposition \ref{prop4.13B}, and we reconcile our
formulas with those of L.~ G\"{o}ttsche, using properties of the difference-zero
hooks of
$P(E)$ (Lemmas
\ref{lemma4.16A},\ref{lemma4.16B}, Remark \ref{remark4.17}). In Section \ref{hookcode}
we define the hook code of a partition $P$ of diagonal lengths $T$, and we show that
the hook code gives an isomorphism between the distributive lattice $\cal P(T)$ of all
partitions of diagonal lengths $T$, and the lattice corresponding to the direct product
of the lattices of partitions whose Ferrer's graphs are enclosed in boxes
$B_i(T),\mu\le i\le j$ (Theorem
\ref{thm4.19}).  In Section \ref{homgroups} we show that there is an additive
isomorphism over $\mathbb Z$ given by the hook code between the homology  $H^\ast
(\G_T)$, and
$H^\ast (\SGrass(T))$, that respects the ${\bf Z}_{2}$ action arising from
complementation or duality (Theorem \ref{thm5.1}), and we determine the Poincar\'{e}
polynomial of
$\G_T$ (Theorem \ref{thm5.2}).
%%%%%%%%%%%%%%%%%%%%%%%%%%%%%%%%%%%%%%%%%%%%%%%%%%%%%%%%%%
\subsection{Graded ideals of $R=k[x,y]$, and linear
systems on ${{\mathbb P}}^{1}$}\label{sec4A}
%%%%%%%%%%%%%%%%%%%%%%%%%%%%%%%%%%%%%%%%%%%%%%%%%%%%%%%%%%
%%%%%%%%%%%%%%%%%%%%%%%%%%%%%%%%%%%%%%%%%%%%%%%%%%%%%%%%%%
Let $R={k}[x,y]$ and  $A$ be an $R$-module,
%%%%%%%%%%%%%%%%%%%%%%%%%%%%%%%%%%%%%%%%%%%%%%%%%%%%%%%%%%
we let $A_{i}$ denote $M^{i}A/M^{i+1}A$, where $M$ is  
the maximal ideal $(x,y)$ of $R$, and
%%%%%%%%%%%%%%%%%%%%%%%%%%%%%%%%%%%%%%%%%%%%%%%%%%%%%%%%%%
we let 
$ H(A)=(t_{0},t_{1},\ldots ) ,
\,t_{i}=\dim_{{k}}A_{i}$
%%%%%%%%%%%%%%%%%%%%%%%%%%%%%%%%%%%%%%%%%%%%%%%%%%%%%%%%%%
denote the Hilbert function of $A$.
%%%%%%%%%%%%%%%%%%%%%%%%%%%%%%%%%%%%%%%%%%%%%%%%%%%%%%%%%%
If $I$ is an ideal of $R$,
%%%%%%%%%%%%%%%%%%%%%%%%%%%%%%%%%%%%%%%%%%%%%%%%%%%%%%%%%%
we let $d(I)=\min\{ i\ | \ I_{i}\neq 0 \}$ be the      
{\em order} of $I$, and if $I$ has finite colength,  
$j(I) = \max\{ i\ | \ I_{i}\neq R_{i} \}$ is the       
{\em socle degree} of $I$.
%%%%%%%%%%%%%%%%%%%%%%%%%%%%%%%%%%%%%%%%%%%%%%%%%%%%%%%%%%
The Hilbert function $T = H(A)$, $A=R/I$ of an Artin
quotient of $R$ satisfies,
\begin{equation}\label{e4.1}
T = (1, 2, \ldots , \mu , t_{\mu}, \ldots , t_{j}, 0), \ \
\mu\geq t_{\mu}\geq \cdots \geq t_{j}>0. 
\end{equation}
%%%%%%%%%%%%%%%%%%%%%%%%%%%%%%%%%%%%%%%%%%%%%%%%%%%%%%%%%%
If $T$ satisfies \eqref{e4.1}, we define the set
$\{\G_T\} = \{ graded \ ideals \ I\subset R \ | \
H(R/I)=T \} . $ 
%%%%%%%%%%%%%%%%%%%%%%%%%%%%%%%%%%%%%%%%%%%%%%%%%%%%%%%%%%
There is a natural inclusion $\iota :
\{\G_T\}\rightarrow \{\BGrass(T)\}$, of the set
$\{\G_T\}$ into the set of closed points of
$\BGrass(T)$, a product of ``Big" Grassmannians: 
$\BGrass(T) = \prod_{\mu\leq i\leq j}
\Grass(i+1-t_{i}, R_{i})$.  
%%%%%%%%%%%%%%%%%%%%%%%%%%%%%%%%%%%%%%%%%%%%%%%%%%%%%%%%%%
\begin{align}\iota:\, \{\G_T\}&\subset \{\BGrass(T)\} =  \prod_{\mu\leq i\leq j}
\{\Grass(i+1-t_{i}, R_{i})\} \notag\\\label{e4.2}
I&\rightarrow
(I_{\mu},I_{\mu +1},\ldots ,I_{j})  
\end{align}
%%%%%%%%%%%%%%%%%%%%%%%%%%%%%%%%%%%%%%%%%%%%%%%%%%%%%%%%%%
We give $\G_T$ the reduced subscheme structure 
induced from this inclusion $\iota$.
%%%%%%%%%%%%%%%%%%%%%%%%%%%%%%%%%%%%%%%%%%%%%%%%%%%%%%%%%%
Then $\G_T$ is a closed projective, nonsingular 
variety having a cover by affine spaces of the same
dimension ([I1], Theorem~2.9).
%%%%%%%%%%%%%%%%%%%%%%%%%%%%%%%%%%%%%%%%%%%%%%%%%%%%%%%%%%
%%%%%%%%%%%%%%%%%%%%%%%%%%%%%%%%%%%%%%%%%%%%%%%%%%%%%%%%%%

%%%%%%%%%%%%%%%%%%%%%%%%%%%%%%%%%%%%%%%%%%%%%%%%%%%%%%%%%%
%\newtheorem{exam41}{Example}[section]
\begin{example}\label{ex4.1}
												Let $T=(1,2,3,2,1)$, then
												$I_{2}=0,I_{3}=\langle f,g\rangle $ is 
            two-
												dimensional, $I_{4}$ has
												dimension four as vector subspace of
												$R_{4}$,
												while $I_{5}=R_{5}$.  The
												inclusion
												$\iota: \G_T\hookrightarrow \BGrass(T)$,
												is
          \begin{equation*}\iota :\G_T\subset \Grass(2,R_{3})\times
            \Grass(4,R_{4})=\Grass(2,4)\timesÊ
            {{\mathbb P}}^{4}, 
           \hspace{.2cm}	I\longmapsto (I_{3},I_{4}).\
\end{equation*} 
          In this example, $\G_T$  is the locus
of pairs	of vector spaces $(V_{3},V_{4})\in
											\Grass(2,R_{3})\times \Grass(4,R_{4})$  
										such that 
												$\langle xV_{3},yV_{3}\rangle \subset V_{4}$.

\end{example}
%%%%%%%%%%%%%%%%%%%%%%%%%%%%%%%%%%%%%%%%%%%%%%%%%%%%%%%%%%
%%%%%%%%%%%%%%%%%%%%%%%%%%%%%%%%%%%%%%%%%%%%%%%%%%%%%%%%%%
Given an ideal $I$ of $\G_T$  we let
$\iota(I)$  denote the sequence of vector spaces
$(I_\mu ,\ldots ,I_{j})$  and we let ${\cal L}(I)
={\cal L}(I_\mu ),\ldots ,{\cal L}(I_{j})$  denote the
corresponding sequence of linear systems on ${\check{\mathbb P}}^1$
 (see Definition~\ref{linearsystem}).
%%%%%%%%%%%%%%%%%%%%%%%%%%%%%%%%%%%%%%%%%%%%%%%%%%%%%%%%%%
%%%%%%%%%%%%%%%%%%%%%%%%%%%%%%%%%%%%%%%%%%%%%%%%%%%%%%%%%%

%%%%%%%%%%%%%%%%%%%%%%%%%%%%%%%%%%%%%%%%%%%%%%%%%%%%%%%%%%
%\newtheorem{lemm42}[exam41]{Lemma}
\begin{lemma}\label{lemma4.2}
Suppose that $T$ satisfies \eqref{e4.1}. 
The maps
$\iota : I\rightarrow \iota(I)$ and 
${\cal L} : I\rightarrow {\cal L}(I)$ 
give isomorphisms among
\begin{enumerate}           
\item[i.] The set of all ideals $I\in \G_T$.
\item[ii.] The set of all collections of vector spaces
	$(V_{\mu}, \ldots , V_{j})$ such that for each 
$i,\ \mu\leq i\leq j,  \ \ V_{i}\subset R_{i}, \ \dim_{k}V_{i}
=i+1-t_{i}$, and such that 
\begin{equation}\label{e4.3}
f\in V_{i},\,  L\in R_{1}\Longrightarrow
Lf\in 	V_{i+1} \,\   for\, \   \mu\leq i<j. 
\end{equation} 
\item[iii.] The set of sequences 
${\cal L} = ({\cal	L}_{\mu},\ldots ,{\cal L}_{j})$ of linear
systems on ${\check{\mathbb P}}^1$ such that for
$\mu\leq i\leq j$, $deg({\cal L}_{i})=i, \ \dim_{k}({\cal
				L}_{i})=i+1-t_{i}$, and satisfying, for $\mu\leq i<j$:
\end{enumerate} 
\begin{equation}\label{e4.4}
\text{ If }\ 
p_{f}=(\sum_{u}n_{u}p_{u}), \  p_{u}\in {\check{\mathbb P}}^1,\text{ satisfies
}\  p_{f}\in {\cal L}_{i} \ 
\text{ and }\ 
if \  p\in {\check{\mathbb P}}^1, \text{ then } \ (p+p_{f}) \in
{\cal L}_{i+1}. 
\end{equation}
\end{lemma}
%%%%%%%%%%%%%%%%%%%%%%%%%%%%%%%%%%%%%%%%%%%%%%%%%%%%%%%%%% 
\begin{proof}
Immediate from the definitions. 
\end{proof}                              
%%%%%%%%%%%%%%%%%%%%%%%%%%%%%%%%%%%%%%%%%%%%%%%%%%%%%%%%%%
%%%%%%%%%%%%%%%%%%%%%%%%%%%%%%%%%%%%%%%%%%%%%%%%%%%%%%%%%%
%\newtheorem{defi43}[exam41]{Definition} 
\begin{definition}\label{def4.3}
		The \em{ramification} $\QRAM(I,p)$ at $p$
		of an ideal
		$I\in \G_T$ in the direction $p\in {\check{\mathbb P}}^1$, given
by $p: L=0, L\in R_{1}$
		 is the sequence
\begin{equation}\label{e4.5}
\QRAM(I,p)=(\ldots , \QRAM(I_{i},p),
						\ldots ),\ \   \mu\leq
						i\leq j\ , 
\end{equation}
where $\QRAM (I_{i},p)$ is from
Definition 1.2. We denote the ramification in the
direction $p: x=0$ by $\QRAM(I)$.
 \end{definition}
%%%%%%%%%%%%%%%%%%%%%%%%%%%%%%%%%%%%%%%%%%%%%%%%%%%%%%%%%% 
Suppose that $\wp = (\wp_{i},\,  \mu\leq i\leq j)$ is a
sequence of partitions, and
that each $\wp_{i}$ is included in ${B}_{i}(T)$, a
$\dim (I_{i})\times \cod (I_{i})$ rectangle. This is the
condition for each $\wp_{i}$ to be possible as a
ramification partition by Lemma 1.9. Suppose also that
$p\in {\check{\mathbb P}}^1$ is a point.\\
%%%%%%%%%%%%%%%%%%%%%%%%%%%%%%%%%%%%%%%%%%%%%%%%%%%%%%%%%%

\noindent {\bf Question 4.3.1.} {\it What is the structure
of the subfamily ${\mathbb V}_{\wp , p}\subset \G_T$ of
graded ideals $I$  such that $\QRAM(I,p)=\wp$? \
What is its dimension?}\\
%%%%%%%%%%%%%%%%%%%%%%%%%%%%%%%%%%%%%%%%%%%%%%%%%%%%%%%%%%

\noindent For each $i$, the condition $\QRAM(I_{i},p)
= \wp_{i}$ \  is by Lemma 2.1B an open dense subset of  a
Schubert condition  on the big Grassmannian
$\Grass (i+1-t_{i},i+1)$. But these conditions for
different $i$ do not intersect properly in the product
$\BGrass(T)$.  First, in order to be compatible, there must
be a  monomial ideal $E\in \G_T$ such that
$\QRAM(E_{i}) = \wp_{i}$ for each $i$. Second, we need to study
the cell ${\mathbb V}(E,p)$. We answer Question 4.3.1 in
Lemmas \ref{lem4.3.1}, Proposition \ref{prop4.5}, and Theorem
\ref{thm4.9} below.\\
%%%%%%%%%%%%%%%%%%%%%%%%%%%%%%%%%%%%%%%%%%%%%%%%%%%%%%%%%%
\begin{lemma}\label{lem4.3.1} The family ${\mathbb V}_{\wp , p}
\subset \G_T$ of graded ideals satisfying $\QRAM (I,p) = \wp$
 is nonempty iff there is a monomial
ideal $E$ such that $\QRAM (E) = \wp$. Then ${\mathbb
V}_{\wp,p} = {\mathbb V}(E,P)$, the family of graded ideals
with initial ideal $E_{L}$ in the basis $(L_{p},C)$ for
$R$. 
\end{lemma}
%%%%%%%%%%%%%%%%%%%%%%%%%%%%%%%%%%%%%%%%%%%%%%%%%%%%%%%%%%
\begin{proof} By Lemma \ref{lemma1.9}, the family ${\mathbb
V}_{\wp , p}$ is the family of all graded ideals $I$ of
Hilbert function $H(R/I) = T$, such that for each $i, \
\mu\leq i\leq j, \ {Q}(I_{i},p) = ({B}(i+1-t_{i},t_{i}) -
\wp_{i})^{\wedge}$. By Lemma~\ref{lemma1.8}, each $I_{i}$ has the
unique initial monomial vector space
$E_{L}(i)$ in the basis $(L,C), L=L_{p}$ for $R$, such
that ${Q}(E_{L}(i),p) = {Q}(I_{i},p)$. The
initial monomial vector spaces $E_{L}(i)$ form an ideal
$E_{L}$. The family of graded ideals ${\mathbb V}_{\wp , p}$
is just ${\mathbb V}(E,p)$, the family having initial ideal
$E_{L}$. The corresponding monomial ideal $E$ in the
variables $(x,y)$ satisfies $\QRAM (E) = \wp$.
\end{proof}
\begin{definition}[The cell ${\mathbb V}(E,p)$ of $\G_T$]\label{def4.4}
If $E\in \G_T$ is an ideal generated by
monomials in $(L,C)$
then	$\mathbb{V}(E,p)$ is the
family of all graded ideals $I\in \G_T$ such that
$\In_{p}(I) = E$. 
\end{definition}
%%%%%%%%%%%%%%%%%%%%%%%%%%%%%%%%%%%%%%%%%%%%%%%%%%%%%%%%%%%
%%%%%%%%%%%%%%%%%%%%%%%%%%%%%%%%%%%%%%%%%%%%%%%%%%%%%%%%%%
%\newtheorem{propos45}[exam41]{Proposition}
\begin{proposition}\label{prop4.5}
														Let $E$ be a monomial ideal with 
														$H(R/E)=T$, 
and $p \in {\check{\mathbb P}}^1$, and suppose that $I$ is
an ideal with $H(R/I)=T$. The following are equivalent:
\begin{align}          
              i.\,& V\in \mathbb{V}(E,p)&\qquad\qquad\qquad\qquad \qquad
\qquad\qquad\qquad\notag\\
  															ii.\, &\QRAM(I,p)=\QRAM(E), \text { i.e. }\QRAM(I_{i},p) = 
																								\QRAM(E_{i}) \text{ for
																						each } i, \mu\leq i\leq j.&\label{e4.6a}\\
																		iii.\, & {Q}(I_{i},p) = 
													{Q}(E_{i}) \text{ for each } i, \ \mu\leq i\leq
j.&\notag
\end{align}
For each $i$, 
\begin{equation}\label{e4.6b}
\QRAM(I_{i},p) =
														({Q}(E_{i})^{c})^{\wedge},
\end{equation}
	the dual $({Q}(E_{i})^{c})^{\wedge}$ of the 
 complement ${Q}(E_{i})^{c}$ to ${Q}(E_{i})$  
 in the $(t_{i})\times (i+1-t_{i})$ rectangle
 ${B}_{i}(T)$.\\
\end{proposition}
%%%%%%%%%%%%%%%%%%%%%%%%%%%%%%%%%%%%%%%%%%%%%%%%%%%%%%%%%%
\begin{proof} The first formula \eqref{e4.6a} is a consequence of
\eqref{e2.2}, applied to each degree-$i$ piece $I_{i}$ and  
 $E_{i},\  \mu\leq i\leq j$. The second follows from
Lemma \ref{lemma1.9}. 
\end{proof}
%%%%%%%%%%%%%%%%%%%%%%%%%%%%%%%%%%%%%%%%%%%%%%%%%%%%%%%%%%
%\newtheorem{defi46}[exam41]{Definition}
		\begin{definition}[Partition ${ P}(E)$ of a 
					monomial ideal $E$]\label{def4.6}
    A monomial ideal $E\in \G_T$ has a
						complementary basis $E^{c}$  of monomials, 
					that we arrange in standard array form, with the
				 monomials $(1,x,\ldots ,x^{p_{0}(E)-1})$ in
					the cobasis forming the top row, and the
					monomials
     $(y^{i}, y^{i}x,\ldots , y^{i}x^{p_{i}(E)-1})$ 
     forming the $i$-th row from the top,
     $i=(0,1,\ldots , s)$. Then
     $E^{c}$  has the shape of a partition 
     ${ P}(E)=(p_{0}, p_{1},\ldots , p_{s} )$ 
					whose
					diagonal lengths are given by $T.$ Conversely,
					the partition $ P$  determines a unique
					monomial ideal $E = E({ P})$ satisfying
     ${P} = { P}(E({ P})).$ 
  \end{definition}
%%%%%%%%%%%%%%%%%%%%%%%%%%%%%%%%%%%%%%%%%%%%%%%%%%%%%%%%%%
%%%%%%%%%%%%%%%%% Inserting picture here %%%%%%%%%%%%%%%%%
\begin{figure}[hbtp]
\begin{center}
\leavevmode 
\begin{picture}(50,35)(25,-35)
\setlength{\unitlength}{.8mm}
\multiput(0,0)(0,-7){2}{\line(1,0){35}}
\put(0,-14){\line(1,0){14}}
\multiput(0,-21)(0,-7){2}{\line(1,0){7}}
\multiput(0,0)(7,0){2}{\line(0,-1){28}}
\put(14,0){\line(0,-1){14}}
\multiput(21,0)(7,0){3}{\line(0,-1){7}}
\put(2.5,-4.5){{\small \mbox{$1$}}}
\put(9.5,-4.5){{\small \mbox{$x$}}}
\put(16.5,-4.5){{\small \mbox{$x^2$}}}
\put(23.5,-4.5){{\small \mbox{$x^3$}}}
\put(30.5,-4.5){{\small \mbox{$x^4$}}}
\put(2.5,-11.5){{\small \mbox{$y$}}}
\put(9.5,-11.5){{\small \mbox{$yx$}}}
\put(2.5,-18.5){{\small \mbox{$y^2$}}}
\put(2.5,-25.5){{\small \mbox{$y^3$}}}
\end{picture}
\end{center}
\vspace{-.7cm}
\protect\caption{Partition ${ P}({\mathcal T})$ of
diagonal lengths $T=(1,2,3,2,1)$,\protect
\\
${\mathcal T}=(x^{5},x^{2}y,xy^{2},y^{4}).$  See Example
\ref{ex4.7}. \label{partT}}
\end{figure}
%\vspace{4cm}
%\special{picture Gr/fig(2) scaled 750}
%\vspace{1mm} %
%\underline{ {\bf Figure 2.}} Partition 
%${ P}({\cal T})$ of diagonal  lengths
%$T=(1, 2, 3, 2, 1), \ \ {\cal T} = (x^{5}, x^{2}y,
%xy^{2}, y^{4}).$ See Example \ref{ex4.7}.\\

%%%%%%%%%%%%%%%%%%%%%%%%%%%%%%%%%%%%%%%%%%%%%%%%%%%%%%%%%%
%%%%%%%%%%%%%%%%%%%%% Example n¼ 4.5 %%%%%%%%%%%%%%%%%%%%%
%\newtheorem{exam47}[exam41]{Example}
			\begin{example}\label{ex4.7}
   		Let $\cal T$ be the monomial ideal  
     ${\cal T} = (x^{5}, x^{2}y, xy^{2}, y^{4})$;
     then
     ${ P}({\cal T}) = (5, 2, 1, 1)$.  The
					diagonals of the shape
     ${ P}(C)$  are the those of
     ${\cal T}^{c}$: namely,
     $\langle 1\rangle , \langle x,y\rangle , \langle x^{2},xy,y^{2}\rangle ,\langle x^{3},y^{3}\rangle ,
\langle x^{4}\rangle ,$
					whose lengths are counted by the Hilbert 
     function $T(R/{\cal T})=(1, 2, 3, 2, 1).$   
			 (See Figure \ref{partT}.)
\end{example}
%%%%%%%%%%%%%%%%%%%%%%%%%%%%%%%%%%%%%%%%%%%%%%%%%%%%%%%%%
%%%%%%%%%%%%%%%%%%%%%%%%%%%%%%%%%%%%%%%%%%%%%%%%%%%%%%%%%
\subsection{Parameters for the cell ${\bf V}(E)$ 
of $\G_T$ }\label{sec4B}
%%%%%%%%%%%%%%%%%%%%%%%%%%%%%%%%%%%%%%%%%%%%%%%%%%%%%%%%%%
   In this section, we will further study the family
			$\mathbb{V}(E)$ of graded ideals having given monomial
			initial ideal $E$. We give $\mathbb{V}(E)$ the reduced
			subscheme structure coming from $\G_T$.  We show
			that the cell is an affine space of
			known dimension. Recall that we order the monomials of
			$R, \ \  1 < y < x < y^{2} < xy < x^{2} < y^{3} <
			\cdots	$
   by total degree, then $y\,$-degree. \\
%%%%%%%%%%%%%%%%%%%%%%%%%%%%%%%%%%%%%%%%%%%%%%%%%%%%%%%%%%

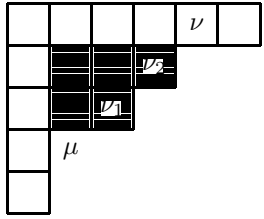
\begin{figure}[hbtp]
\begin{center}
\leavevmode 
\begin{picture}(50,40)(20,-40)
\setlength{\unitlength}{.8mm}
\thicklines
\multiput(0,0)(0,-7){2}{\line(1,0){42}}
\put(0,-14){\line(1,0){28}}
\put(0,-21){\line(1,0){21}}
\multiput(0,-28)(0,-7){2}{\line(1,0){7}}
\multiput(0,0)(7,0){2}{\line(0,-1){35}}
\multiput(14,0)(7,0){2}{\line(0,-1){21}}
\put(28,0){\line(0,-1){14}}
\multiput(35,0)(7,0){2}{\line(0,-1){7}}
\put(10.5,-24.5){\makebox(0,0){$\mu$}}
\put(17.5,-17.5){\makebox(0,0){$\nu_{1}$}}
\put(24.5,-10.5){\makebox(0,0){$\nu_{2}$}}
\put(31.5,-3.5){\makebox(0,0){$\nu$}}
\thinlines
\multiput(7.5,-7.5)(0,-.5){13}{\line(1,0){6}}
\multiput(14.5,-7.5)(0,-.5){13}{\line(1,0){6}}
\multiput(7.5,-14.5)(0,-.5){13}{\line(1,0){6}}
\multiput(21.5,-7.5)(0,-.5){4}{\line(1,0){6}}
\multiput(14.5,-14.5)(0,-.5){4}{\line(1,0){6}}
\multiput(14.5,-20.5)(0,.5){4}{\line(1,0){6}}
\multiput(21.5,-13.5)(0,.5){4}{\line(1,0){6}}
\multiput(14.5,-16.5)(0,-.5){5}{\line(1,0){1}}
\multiput(21.5,-9.5)(0,-.5){5}{\line(1,0){1}}
\multiput(20.5,-16.5)(0,-.5){5}{\line(-1,0){1.5}}
\multiput(27.5,-9.5)(0,-.5){5}{\line(-1,0){1.5}}
\end{picture}
\end{center}
\vspace{-1cm}
\protect\caption{Here ${P} = (6, 4,
3, 1, 1)$, $E=(x^{6}, yx^{4}, y^{2}x^{3}, y^{3}x, y^{5})$, 
$\mu = y^{3}x$, 
$f(\mu , I) = \mu - \alpha_{\mu\nu_{1}}\nu_{1} - 
\alpha_{\mu\nu_{2}}\nu_{2} - \alpha_{\mu\nu}\nu  \in I$,
and  $(\alpha_{\mu\nu_{1}},\alpha_{\mu\nu_{2}})$ are
parameters for $I$. The pair $(\mu,\nu_{2})\in {\cal
S}_{4}(E)$ corresponds to the shaded hook of ${ P}$.
See Example \ref{ex4.8.1} and Remark \ref{rem4.8.2}.\label{fig-parameters}}
\end{figure}

 We now define parameters for ideals $I$ in the
cell ${\mathbb V}(E)$.
%%%%%%%%%%%%%%%%%%%%%%%%%%%%%%%%%%%%%%%%%%%%%%%%%%%%%%%%%%
%%%%%%%%%%%%%%%%%%%% Definition n¼ 4.8 %%%%%%%%%%%%%%%%%%%
%\newtheorem{defi48}[exam41]{Definition}  \begin{defi48}
\begin{definition}\label{def4.8}
			Given a monomial ideal $E$ or
				partition ${P} = {P}(E)$ we let   
			${\cal S}_{i}(E)$ or  ${\cal S}_{i}({P})$
		denote the set of ordered pairs of
			monomials
\begin{align} 
			{\cal S}_{i}(E)&=\{ (\mu,\nu)\, | \,
			\mu\in E,
			(\mu:y)\not\in
   E,\  \mu\not\in E,\  x\mu\in E\notag\\   \label{e4.7}
			 & \text{ degree }\,\mu = \text{ degree }\nu = i,\text{ and } \ 
			\mu<\nu\}.	
\end{align}
We let ${\cal S}(E)$  or ${\cal S}(P)$
		denote the union 
\begin{equation}\label{e4.7a}
{\cal S}(E)=\bigcup_{\mu\leq i\leq j}{\cal S}_{i}(E). 
\end{equation}
		If $\mu$ is a monomial of $E$ and 
$I$ is an ideal with initial ideal $\In(I) = E$ we
		let
\begin{equation}\label{e4.8}
f(\mu,I)=\mu -\, \sum_{\nu
>\mu} \alpha_{\mu\nu}(I)\cdot\nu, 
		\ \,\alpha_{\mu\nu}\in k,\ \, \nu\not\in E, \ \
deg(\nu)=deg(\mu) 
\end{equation}
 denote the unique element of $I$  with
		leading term $\mu$.The parameters for the
graded ideals $I$ in the cell ${\mathbb V}(E)$ 
will be the set of coefficients
$\{\alpha_{\mu\nu}(I) \ | \ (\mu ,\nu)\in {\cal
S}(E)\}.$  (See Example \ref{ex4.8.1} and Figure \ref{fig-parameters}).
\end{definition}
\begin{example}\label{ex4.8.1}
If $T=(1,2,3,4,4,1)$ and the
initial ideal is $E=(x^{6}, yx^{4}, y^{2}x^{3}, y^{3}x, y^{5})$
then ${P}(E) = (6, 4, 3, 1, 1)$. If $\mu = y^{3}x$,
and $I\in {\mathbb V}(E)$, then \\ 
$f(\mu , I) = \mu - \alpha_{\mu\nu_{1}}\nu_{1} - 
\alpha_{\mu\nu_{2}}\nu_{2} - \alpha_{\mu\nu}\nu = 
\mu - \alpha_{\mu\nu_{1}}y^{2}x^{2} - 
\alpha_{\mu\nu_{2}}yx^{3} - \alpha_{\mu\nu}x^{4} \in I$,
and $(\alpha_{\mu\nu_{1}},\alpha_{\mu\nu_{2}})$ are among
the parameters for $I$. Each pair $(\mu ,\nu)\in {\cal
S}(E)$ corresponds to a difference-one hook of ${
P}$ with hand $\nu$ and foot $\mu : y$. (See Figure
\ref{fig-parameters}).
\end{example}
%%%%%%%%%%%% inserer ici la figure 3 %%%%%%%%%%%%%%%%%%%%
%
%\vspace{2.5cm}
%\special{picture params scaled 750}
%\vspace{1mm}

%%%%%%%%%%%%%%%%%%%%%%%%%%%%%%%%%%%%%%%%%%%%%%%%%%%%%%%%%
%
%\noindent \underline { Figure 3.} Here ${P} = (6, 4,
%3, 1, 1)$, $E=(x^{6}, yx^{4}, y^{2}x^{3}, y^{3}x, y^{5})$, 
%$\mu = y^{3}x$, 
%$f(\mu , I) = \mu - \alpha_{\mu\nu_{1}}\nu_{1} - 
%\alpha_{\mu\nu_{2}}\nu_{2} - \alpha_{\mu\nu}\nu  \in I$,
%and  $(\alpha_{\mu\nu_{1}},\alpha_{\mu\nu_{2}})$ are
%parameters for $I$. The pair $(\mu,\nu_{2})\in {\cal
%S}_{4}(E)$ corresponds to the shaded hook of ${ P}$.
%( See Example \ref{ex4.8.1} and Remark \ref{rem4.8.2}).
%%%%%%%%%%%%%%%%%%%%%%%%%%%%%%%%%%%%%%%%%%%%%%%%%%%%%%%%%%%%
%%%%%%%%%%%%%%%%%% remark 4.8.2. %%%%%%%%%%%%%%%%%%%%%%%%%
%%%%%%%%%%%%%%%%%%%%%%%%%%%%%%%%%%%%%%%%%%%%%%%%%%%%%%%%%%%
\begin{remark}\label{rem4.8.2} Note that ${\cal S}_{i}(E)$
consists of the pairs $(\mu,\nu)$ of monomials of degree
$i$, where $\nu$ is the right endpoint of a row of
$E^{c}$, where $\mu$ in $E$ is vertically  just below a
column of $E^{c}$, and $\mu$ is diagonally below $\nu$,
when we view the monomials of $E^{c}$ as  filling the
shape of the partition ${ P}(E)$. The degree-$i$
difference-one hooks of ${P}$ is the  set  
${\cal H}({P})_{i} = \{(\mu,\nu)\ | \  (y\mu,\nu) 
\in {\cal S}_{i}({P})\}$. (See Figure~\ref{fig-parameters} and
Definition \ref{def4.14}).
\end{remark}
%%%%%%%%%%%%%%%%%%%% Theorem n¼ 4.9 %%%%%%%%%%%%%%%%%%%%%
%\newtheorem{theo49}[exam41]{Theorem}
	\begin{theorem}\label{thm4.9}
			The cell ${\mathbb V}(E)$ of $\G_T$  is an affine space 
   ${\bf A}^{s(E)}$ of dimension 
   $s(E) = \#{\cal S}(E) = \#{\cal H}({P})$. The
			parameters of ${\mathbb V}(E)$ 
			are $\{\alpha_{\mu\nu}(I) \ | \  (\mu,\nu)\in {\cal
			S}(E)\}$.\\
\end{theorem}
%%%%%%%%%%%%%%%%%%%%%%%%%%%%%%%%%%%%%%%%%%%%%%%%%%%%%%%%%%
\noindent
{\bf Historical Remark.} Analogous results were well
			known to the Nice Hilbert scheme group, in particular
			J.~Brian\c{c}on, as early as 1972, for the family
			$\mathrm{Z_{T}}$ of
			all ideals (not just graded) defining quotient algebras
			$A=R/I$ of Hilbert function $T$. 
			J.~Brian\c{c}on studied  ``vertical strata"        
 		${\mathbb  Z}(E)$ of the punctual  Hilbert scheme
			$\Hilb^{n}R$, which when restricted to graded ideals 
			are identical with the cells ${\mathbb V}(E)$; his
			vertical strata involve ideals with different
			Hilbert functions, and are a key tool in his proof of the irreducibility of
the local punctual Hilbert scheme $\Hilb^n \mathbb C\{ x,y\}$ \cite{Br}.       
   J.~Yam\'{e}ogo gives a proof of Theorem \ref{thm4.9} in [Y3] and
			[Y4] using the map from $I$ to  $(I:x)$ and an
			induction,
			following the approach of Brian\c{c}on. Our proof here generalizes
			that given in [I1] for the special case of the
			``generic" or big cell $E(0)$ of maximum dimension, 
			where $E(0)_{i} = \langle y^{i},\ldots
    , y^{t_{i}}x^{i-t_{i}}\rangle $. 
			
			J. Yam\'{e}ogo showed that the ``vertical strata'' subvarieties $\mathbb Z(E)$ and
${\mathbb V}(E)$ are identical to the
			families of ideals collapsing to $E$  under a
			${\bf C}^{\ast}\,$-action (see \cite{Y3,Y4}),
			L.~G\"{o}ttsche proved that these latter are affine spaces, and he calculated the
dimension of
$\mathbb Z(E)$, and then of $\mathbb{V}(E)\subset \G_T$ from the dimension of the
fibre of $\mathrm {Z_T}$ over $\G_T$ \cite{G2}.  Given this
		identification, L.~G\"{o}ttsche's formula for $\dim \mathbb V(E)$
		(see Proposition \ref{prop4.15} below) was the first  covering
			explicitly 
			all the cells ${\mathbb V}(E)$ in the graded case. But his
			dimension formula is different than ours in Theorem \ref{thm4.9}!
		 We reconcile the two formulas
			using a combinatorial result from Part I \cite{IY2} (see Remark \ref{remark4.17}).
In our proof of Theorem \ref{thm4.9},
		 the parameters for $\mathbb{V}(E)$ give
		the geometric meaning of the ``hook code" of \S \ref{hookcode}.
%%%%%%%%%%%%%%%%%%%%%%%%%%%%%%%%%%%%%%%%%%%%%%%%%%%%%%%%%
%%%%%%%%%%%%%%%%%% Proof of Theorem 4.9. %%%%%%%%%%%%%%%%
	\begin{proof}[Proof of Theorem \ref{thm4.9}]                         
		The key to the proof is that
		since we are working in only two variables, the  minimal
		relations between standard generators of
		ideals  in ${\mathbb V}(E)$ have an echelon or upper
		triangular form. List the monomials
		$\beta_{0},\ldots ,\beta_{p}, \ p=p_{0}(E)$, in $E$ 
		just below the pattern:
		$\beta_{i} = x^{i}y^{q(i)} \in E $ and $(\beta_{i}:y)
		\in E^{c}$. \\
%%%%%%%%%%%%%%%%%%%%%%%%%%%%%%%%%%%%%%%%%%%%%%%%%%%%%%%%%
		The sequence $(q(0),q(1),\ldots ,q(p_{0}) = 0)$ is the 
		dual partition ${ P}^{\wedge}$ to ${ P}(E)$.
		Since $(x^{i})=(\beta_{i},\ldots ,\beta_{p})\oplus
		E^{c}\cap (x^{i})$, and $\beta_{i}$ is the leading form
		of $f(\beta_{i},I)$, it follows that for ideals     
		$I\in {\mathbb V}(E)$,
\begin{equation}\label{e4.9}
(x^{i})\cap I=
		(f(\beta_{i},I),
		f(\beta_{i+1},I),\ldots ,
		f(\beta_{p},I)).               
\end{equation}
%%%%%%%%%%%%%%%%%%%%%%%%%%%%%%%%%%%%%%%%%%%%%%%%%%%%%%%%%%
{\em Induction step : comparison to} $(I:x^{i})$.
%%%%%%%%%%%%%%%%%%%%%%%%%%%%%%%%%%%%%%%%%%%%%%%%%%%%%%%%%%
We will let $(E:x^{i})$ denote $(E\cap (x^{i}):(x^{i}))$,
the monomial  ideal with cobasis   $(E^{c}:x^{i}) =
(E^{c}\cap (x^{i}):x^{i})$. The shape ${ P}(E:x^{i})$
of the cobasis $(E^{c}:x^{i})$ is that of ${ P}$  with
the first $i$-columns omitted. Thus ${ P}(E:x^{i})$ is
the  partition dual to $(q(i),q(i+1),\ldots ,q(p_{0}))$.
We regard the monomials $\mu', \nu'$ of $(E:x^{i})$ and
$(E^{c}:x^{i})$ as those  of  $E^{c}$ and $E$ shifted   
left by $i$, and will thus for convenience set $\mu'=
(\mu:x^{i})$ and $(\nu'=\nu:x^{i})$ in  Claim B below. The
cell ${\mathbb V}(E:x^{i})$ of Claim B below lies on
$G_{T(i)}$ where $T(i)$ is the Hilbert function of
$R/(E:x^{i})$.\\
%%%%%%%%%%%%%%%%%%%%%%%%%%%%%%%%%%%%%%%%%%%%%%%%%%%%%%%%%%
{\em Claim A.} If $I\in{\mathbb V}(E)$, then for each $i \leq
p$ \  the polynomial $f(\beta_{i},I)$ is uniquely
determined by the set  of coefficients
\begin{equation}\label{e4.10a}
\{\alpha_{\mu\nu}(I) \ | \ \mu\geq \beta_{i}\ 
and\  (\mu,\nu)\in {\cal S}(E)\}.
\end{equation}                     
%%%%%%%%%%%%%%%%%%%%%%%%%%%%%%%%%%%%%%%%%%%%%%%%%%%%%%%%%%     
{\em Claim B.} For each $i\leq p$, given any set of
constants  $\{\psi_{\mu\nu}\in {k} \ | \  (\mu,\nu)\in
{\cal S}(E)$ \  and \  $\mu \geq \beta_{i}\}$ there is a
unique ideal $I_{\psi}^{(i)}\in {\mathbb V}(E:x^{i})$ such that
$\alpha_{\mu'\nu'}(I_{\psi}^{(i)})=\psi_{\mu\nu}.$                 
The ideal $x^{i}I_{\psi}^{(i)}$ has generators
$f(\beta_{i}),\ldots , f(\beta_{p})$ with initial terms
$\beta_{i},\ldots ,\beta_{p}$, respectively, and
satisfies  
\begin{equation}\label{e4.10b}
	x^{i}I^{(i)}\oplus E^{c}\cap
(x^{i})=(x^{i}). 
\end{equation}
%%%%%%%%%%%%%%%%%%%%%%%%%%%%%%%%%%%%%%%%%%%%%%%%%%%%%%%%%%
%%%%%%%%%%%%%% inserer ici la figure 4 %%%%%%%%%%%%%%%%%%%
\begin{figure}[hbtp]
\begin{center}
\leavevmode 
\begin{picture}(45,25)(45,-25)
\setlength{\unitlength}{1mm}
\thicklines
\multiput(0,0)(0,-4){2}{\line(1,0){36}}
\put(0,-8){\line(1,0){24}}
\put(0,-12){\line(1,0){18}}
\multiput(0,-16)(0,-4){2}{\line(1,0){6}}
\multiput(0,0)(6,0){2}{\line(0,-1){20}}
\multiput(12,0)(6,0){2}{\line(0,-1){12}}
\put(24,0){\line(0,-1){8}}
\multiput(30,0)(6,0){2}{\line(0,-1){4}}
\put(3,-22){\makebox(0,0){{\small \mbox{$\beta_{0}$}}}}
\put(9,-14){\makebox(0,0){{\small \mbox{$\beta_{1}$}}}}
\put(15,-14){\makebox(0,0){{\small \mbox{$\beta_{2}$}}}}
\put(21,-10){\makebox(0,0){{\small \mbox{$\beta_{3}$}}}}
\put(27,-6){\makebox(0,0){{\small \mbox{$\beta_{4}$}}}}
\put(33,-6){\makebox(0,0){{\small \mbox{$\beta_{5}$}}}}
\put(39,-2){\makebox(0,0){{\small \mbox{$\beta_{6}$}}}}
\thinlines
\multiput(12,-0.5)(0,-.5){7}{\line(1,0){6}}
\multiput(18,-0.5)(0,-.5){7}{\line(1,0){6}}
\multiput(24,-0.5)(0,-.5){7}{\line(1,0){6}}
\multiput(30,-0.5)(0,-.5){7}{\line(1,0){6}}
\multiput(12,-4.5)(0,-.5){7}{\line(1,0){6}}
\multiput(18,-4.5)(0,-.5){7}{\line(1,0){6}}
\multiput(12,-8.5)(0,-.5){7}{\line(1,0){6}}
\end{picture}
\end{center}
\vspace{-0.2cm}
\protect\caption{Partition ${ P}(E)=(6,4,3,1,1)$, and, 
shaded, ${ P}(E:x^{2})=(4,2,1)$. Here
$xf(\beta_{1},I)\subset (x^{2})\cap I=(f(\beta_{2},I),\ldots
,f(\beta_{5},I),\beta_{6})$. See Example
\ref{ex4.9.1}.}\label{partitionP(E)}
\end{figure}
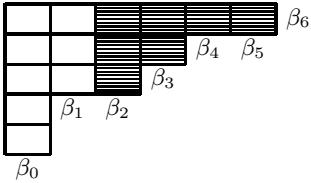
%
% \hspace{5.4cm} \underline{Figure 4.}  \ \ {\small
%Partition} ${P}(E) = (6, 4, 3, 1, 1)${\small ,} 
%
% \hspace{5.4cm} {\small  and, shaded,} 
% ${ P}(E : x^{2}) = (4, 2, 1)$. {\small Here}
%
% \hspace{5.4cm} 	$xf(\beta_{1},I) \subset
%(x^{2})\cap I =	(f(\beta_{2},I), \ldots , f(\beta_{5},I),
%\beta_{6})$.
%
% \hspace{5.4cm} {\small See Example \ref{ex4.9.1}.}\\
%												
%
%\special{picture partition(fig4) scaled 700} 
%
%
%
%\vspace{7mm}

%%%%%%%%%%%%%%%%%%%%%%%%%%%%%%%%%%%%%%%%%%%%%%%%%%%%%%%%%%

\noindent {\it Proof of Claims.} We will show these by
descending  induction on $i$, beginning with $i=p$. Since
$I\supset  x^{p}=\beta_{p}$, and the sets of A, B are
vacuous when  $i=p$, the claims are satisfied for $i=p$.
Suppose they are satisfied for $i+1$, and that $I\in {\mathbb
V}(E)$. Then $xf(\beta_{i},I)\in (x^{i+1})\cap
I.$\\                 The coefficients
$\{\alpha_{\mu\nu}(I)\, |\, \mu = \beta_{i}$ and 
$(\mu,\nu)\in {\cal S}(E)\}$ are precisely those of
$f(\beta_{i},I)$ on monomials $\nu$ such that $x\nu$ lands
{\em outside} $\, \langle E^{c}\rangle $ in $xf(\beta_{i},I)$. These can
be reduced in  a standard way by elements of 
$(f(\beta_{i+1},I),\ldots ,f(\beta_{p},I))$ to a
remainder  $g_{i}(I)$ in $\langle E^{c}\rangle _{i+1}$. The remaining
coefficients of $f(\beta_{i},I)$ are  on monomials $\nu$
such that $x\nu$ lands {\em inside}  $\, \langle E^{c}\rangle _{i+1}$:
these must be $-g_{i}(I)$, since their sum with  the
remainder $g_{i}(I)$ is in $I\cap \langle E^{c}\rangle $ so must be $0$.
This proves Claim A.\\                                   

Suppose Claim B is satisfied for $i+1$, and that the set 
$\{\psi_{\mu\nu}\in {k} | (\mu,\nu)\in {\cal S}(E)$
and    $\mu \geq \beta_{i}\}$ is specified, and let
$I^{(i+1)}$ denote  the unique ideal in ${\mathbb
V}(E:x^{i+1})$ satisfying Claim B. We have thus determined
an ideal $x^{i+1}I^{(i+1)}= (f(\beta_{i+1}),\ldots
,f(\beta_{p}))$ such that 
\begin{equation}\label{e4.11}
x^{i+1}I^{(i+1)}\oplus
E^{c}\cap (x^{i+1})=(x^{i+1}). 
\end{equation}
%%%%%%%%%%%%%%%%%%%%%%%%%%%%%%%%%%%%%%%%%%%%%%%%%%%%%%%%%%
By the proof of Claim A, there is a unique homogenous 
polynomial $f(\beta_{i})$ having the designated
coefficients $\{\alpha_{\mu\nu}=\psi_{\mu\nu}\in {k}
\  |\  (\mu,\nu)\in {\cal S}(E)\  and \ \mu = \beta_{i}\}$,
and satisfying $xf(\beta_{i})\in (f(\beta_{i+1}),\ldots
,f(\beta_{p}))$, so we have a  relation $r_{i}$,
\begin{equation}\label{e4.12}
r_{i} : xf(\beta_{i}) = h_{i+1}f(\beta_{i+1}) +
\ldots +  h_{p}f(\beta_{p}).
\end{equation}
Since $\beta_{i}$ has leading term of $x\,$-degree $i$,
any  other  relation $r$ between $f(\beta_{i})$ and
$f(\beta_{i+1}),\ldots , f(\beta_{p})$ involves a
multiple of $xf(\beta_{i})$, so can be reduced by the
relation $r_{i}$  to one involving only\linebreak
$f(\beta_{i+1}),\ldots ,f(\beta_{p})$.             
But \eqref{e4.11}  together with the form of $f(\beta_{i})$ now
show  that                           
\begin{equation}\label{e4.13}
(f(\beta_{i}),\ldots ,f(\beta_{p}))\oplus
E^{c}\cap (x^{i})=(x^{i}). 
\end{equation}                        
We set $I^{(i)} = (f(\beta_{i}):x^{i},\ldots
,f(\beta_{p}):x^{i})$.  By \eqref{e4.13} and the form of
$f(\beta_{i})$, the ideal $I^{(i)}$ is  in ${\mathbb
V}(E:x^{i})$. This completes the proof of Claim B.      
Taking $i = 0$, and $I = I^{0} $ \  we obtain Theorem~\ref{thm4.9}. 
\end{proof}
%%%%%%%%%%%%%%%%%%%%%%%%%%%%%%%%%%%%%%%%%%%%%%%%%%%%%%%%%
%%%%%%%%%%%%% mettre ici l'exemple 4.9.1 %%%%%%%%%%%%%%%%

\begin{example}[Induction step]\label{ex4.9.1}
Consider ${ P}(E) = (6, 4, 3, 1, 1), \beta_{2} = 
xy^{3}$, and  $f(\beta_{2},I) = \mu -
\alpha_{\mu\nu_{1}}y^{2}x^{2} - \alpha_{\mu\nu_{2}}yx^{3}
- \alpha_{\mu\nu}x^{4}$. If $(x^{2}\cap I)$ has been
chosen, the new information in $f(\beta_{2},I)$ is the
coefficients $(\alpha_{\mu\nu_{1}}, \alpha_{\mu\nu_{2}})$
on those monomials $\nu_{1} , \nu_{2}$ such that
$x\nu_{1}, x\nu_{2}$ lie {\em outside} of ${ P}(E)$:
these may be chosen arbitrarily (Claim B). In the usual
standard basis way, $xf(\beta_{2},I)$ may be reduced using
$(x^{2}\cap I = (f(\beta_{2},I), \ldots , \beta_{6})$ to a
linear combination of cobasis monomials in ${
P}(E)\cap(M^{5})$, which must be zero. The coefficient
$\alpha_{\mu\nu}$ of $f$ on a monomial $\nu = x^{4}$ such
that $x\nu$ lies {\em inside} of ${ P}(E)$ is thus
determined by $(x^{2})\cap I$ and by
$(\alpha_{\mu\nu_{1}}, \alpha_{\mu\nu_{2}})$ (Claim A).
(See Figure \ref{partitionP(E)}).
\end{example}
%%%%%%%%%%%%%%%%%%%%%%%%%%%%%%%%%%%%%%%%%%%%%%%%%%%%%%%%%
%\newtheorem{defi410}[exam41]{Definition}
\begin{definition}[Small Grassmannians and the big cell of
$\G_T$]\label{def4.10} Let $T$ satisfy \eqref{e4.1}. Then we let
$\delta_{i} = \delta_{i}(T) =t_{i-1}-t_{i}$.  We let
$\SGrass(T)$  denote the product of small 
Grassmannians,
\begin{equation}\label{e4.14}
\SGrass(T) =\prod_{\mu (T)\leq i\leq
j(T)}\Grass(\delta_{i+1},1+\delta_{i}+\delta_{i+1}).
\end{equation}
 We let $E=E_{0}$ denote the big
cell of
$\G_T$                                               
determined by the unique partition           
${P}_{0}(E)= (p_{1},p_{2},\ldots ,p_{r})$ of
diagonal lengths $T$ having  distinct parts:
${P}_{0}(T)$ is maximal in the partial order on
${\cal P}(T)$  defined by the inclusion of the
partitions ${Q}(E({P}))$  (See also Part I,
\cite[Definition 2.18B]{IY2}). For $i\geq d$ and $E=E_{0}$ 
we define the vector  spaces $U_{i} = (E^{c}_{i+1}:x)$, 
$V_{i} = (x\langle E^{c}\rangle _{i-1}+y\langle E^{c}\rangle _{i-1})/U_{i}$, 
and $W_{i} = \langle E^{c}\rangle _{i}/U_{i}$ . Then 
$\dim_{k}W_{i} = \delta_{i+1}$ , and $\dim_{k}V_{i} =
t_{i-1}+1-t_{i+1} =  \delta_{i}+\delta_{i+1}+1$.
Given $I\in {\mathbb V}(E)$ we
define
\begin{equation}\label{e4.15a}
             {\cal G}_{i}(I) =
(I_{i}\cap ( x\langle E^{c}\rangle _{i-1}+y\langle E^{c}\rangle _{i-1} ) +
U_{i})/U_{i}, 
\end{equation}  
a complement to
$W_{i}$ {\em in} $V_{i}$ . We
define                                             
${\cal G} : {\mathbb V}(E_{0})\rightarrow  \SGrass(T)$ 
by
\begin{equation}\label{e4.15b}
{\cal G}(I) = ({\cal G}_{\mu (T)}(I),\ldots ,{\cal
G}_{j(T)}(I)). 
\end{equation}
\end{definition}
%%%%%%%%%%%%%%%%%%%%%%%%%%%%%%%%%%%%%%%%%%%%%%%%%%%%%%%%%
%%%%%%%%%%%%%%%%%% Proposition n¼ 4.11 %%%%%%%%%%%%%%%%%%%
%\newtheorem{propos411}[exam41]{Proposition}
	\begin{proposition}\label{prop4.11}
			There is a birational map ${\cal G} : \G_T 
			\rightarrow \SGrass(T)$, that is defined on the big 
    cell ${\mathbb V}(E_{0})$.
\end{proposition}
%%%%%%%%%%%%%%%%%%%%%%%%%%%%%%%%%%%%%%%%%%%%%%%%%%%%%%%%%
%%%%%%%%%%%%%% Proof of Proposition n¼ 4.11 %%%%%%%%%%%%%%
\begin{proof} The proof of Theorem \ref{thm4.9} shows that ${\mathbb
V}(E_{0})$ is an  open dense affine cell in $\G_T$, and
gives its parameters; these parameters can be identified
with the usual parameters for the corresponding open
dense cell ${\cal G}({\mathbb V}(E_{0}))$ of $\SGrass(T)$.  
\end{proof}
%%%%%%%%%%%%%%%%%%%%%%%%%%%%%%%%%%%%%%%%%%%%%%%%%%%%%%%%%
%%%%%%%%%%%%%%%%%%% Section n¼ 4C %%%%%%%%%%%%%%%%%%%%%%%
\subsection{Dimension of the cell ${\bf V}(E)$ and the
hooks of $P(E)$ }\label{dimcell}
%%%%%%%%%%%%%%%%%%%%%%%%%%%%%%%%%%%%%%%%%%%%%%%%%%%%%%%%%
In this section we obtain a new formula for the dimension
of the cells $\mathbb{V}(E)$, and compare  with a
previous dimension result of L.~G\"{o}ttsche. Since
G\"{o}ttsche's result is based on his dimension formula
for cells of nongraded ideals, we introduce the family
${\mathbf  Z}_T $ of all algebra quotients $A=R/I$ of $R$ with
Hilbert function $H(A)=T$.
%%%%%%%%%%%%%%%%%%%%%%%%%%%%%%%%%%%%%%%%%%%%%%%%%%%%%%%%%%
%%%%%%%%%%%%%%%%%%%%%%%%%%%%%%%%%%%%%%%%%%%%%%%%%%%%%%%%%%
%\newtheorem{defi412}[exam41]{Definition}
\begin{definition}\label{def4.12}
Let ${\mathbb  Z}(E) \subset {\mathbf  Z}_T $  parametrize the family
of all ideals of $R$ having initial monomial ideal $E$. 
We let
\begin{align}\label{e4.16a}
 {\cal W}^{+}(E)&=\{(\mu,\nu) \  | \ 
			\mu\in E,\  
				(\mu:y)\not\in E,\  \nu\not\in E,\  x\nu\in E, \ 
		 \text{	and }\,\   degree\,  \mu< \ degree\, \nu\}.\\\label{e4.16b}
			{\cal W}^{-}(E)&=\{(\mu,\nu) \  | \ 
\mu\in E, \  
    (\mu:x)\not\in
				E, \  \nu\not\in E, \  y\nu\in  E, \ \text{ and degree }
			\mu<\text{ degree }  \nu\}.\\\notag            
{\cal W}(E)&={\cal W}^{+}(E)\cup {\cal
W}^{-}(E)
\end{align}
and we let $w(E)$ denote  $\#{\cal
W}(E)$.
\end{definition}  
%%%%%%%%%%%%%%%%%%%%%%%%%%%%%%%%%%%%%%%%%%%%%%%%%%%%%%%%%%%
We let
%%%%%%%%%%%%%%%%%%%%%%%%%%%%%%%%%%%%%%%%%%%%%%%%%%%%%%%%%%%
$h:  {\mathbf  Z}_T  \rightarrow \G_T$ denote the morphism 
defined on closed points by $h(A) = Gr_{m}(A)$ ,    
the associated graded algebra $A^{\star}$ of $A$ 
with respect to  the maximal ideal $ M = (x,y)$  
of $R$.  
Recall that $n(T) = \sum t_{i}$, that
the  initial  degree of ideals in $\G_T$ is  $\mu (T) =
\min\{i \  | \  t_{i} \leq i\}$, the socle degree of
$R/I$ is   $j(T) = \max \{i \  | \  t_{i} \neq 0\}$ , and we
let $\delta_{i} = t_{i-1}-t_{i}$. In [I1], Theorems 2.11 
and 3.14, the first author showed the first of the
following two Propositions. The second readily follows
from similar arguments. 
%%%%%%%%%%%%%%%%%%%%%%%%%%%%%%%%%%%%%%%%%%%%%%%%%%%%%%%%%
%%%%%%%%%%%%%%%%%% Proposition n¼ 4.13 %%%%%%%%%%%%%%%%%%
%\newtheorem{propos413}[exam41]{Proposition}
	\begin{proposition}[A. Iarrobino,{[I1]}]\label{prop4.13A}
 			If $char({k})=0$ or is greater than $j(T)$, then
				the
				morphism $h :  {\mathbf  Z}_T  \rightarrow \G_T$ is a
				locally trivial fibration whose fiber is an affine 
				space ${\bf A}^{f(T)}$, of dimension 
		\begin{equation}\label{e4.17a}
f(T)\ =\  n(T)\ -\
				\sum_{i\geq \mu (T)} \delta_{i}(\delta_{i}+1)/2 \  -\
				\sum_{i\geq \mu (T)} (\delta_{i}+1)\delta_{i+1}.       
		\end{equation}
\end{proposition}
%%%%%%%%%%%%%%%%%%%%%%%%%%%%%%%%%%%%%%%%%%%%%%%%%%%%%%%%%
Here\  $\dim \G_T\ =\ \sum_{i\geq
\mu (T)} (\delta_{i}+1)\delta_{i+1}$\ ,  and   $\dim
{\mathbf  Z}_T  \ =\  n(T)\ -\ \sum_{i\geq \mu (T)}
\delta_{i}(\delta_{i}+1)/2$\ .\\
%%%%%%%%%%%%%%%%%%%%%%%%%%%%%%%%%%%%%%%%%%%%%%%%%%%%%%%%%
\begin{proposition}\label{prop4.13B}. The fiber $F_{{\mathbb Z}/{\mathbb V}}(E)$ of
the cell
${\mathbb Z}(E)$ over ${\mathbb V}(E)$ satisfies,
\begin{equation}\label{e4.17b}
F_{{\mathbb Z}/{\mathbb V}}(E)\cong {\bf A}^{w(E)}, \text { an affine space }.
\end{equation}
\end{proposition}
%%%%%%%%%%%%%%%%%%%%%%%%%%%%%%%%%%%%%%%%%%%%%%%%%%%%%%%%%
\begin{proof}[Comment on proof]  The proofs of
Propositions
	\ref{prop4.13A} and \ref{prop4.13B} are  analogous to that of 
		Theorem~\ref{thm4.9}. However Proposition \ref{prop4.13A} requires in
		addition a proof that in a generic basis $X, Y$ for
		$R_{1}$, a given ideal $I$ lies in the maximal cell
		$E_{0}$: this requires the  assumption on
		$\cha ({k})$.		 
A simple reconciliation of the
dimension formulas in Proposition \ref{prop4.13A} and Proposition
\ref{prop4.13B} follows from our combinatorial results in Part I
\cite{IY2}: see  Remark \ref{remark4.17} below, which gives an
alternative proof of Proposition \ref{prop4.13B}. 
\end{proof}\par       
		L.~G\"{o}ttsche gives a different formula for the
		dimension  of the cell $Z(E)$, from which he 
		derives a formula for $\dim({\mathbb V}(E))$ by using
		Proposition \ref{prop4.13B}. In order to show the G\"ottsche
formula from ours, we need to define
		the set ${\cal H}^{a}(E)_{i}$ of  degree-$i$
difference-$a$  hooks of $E^{c}$.
%%%%%%%%%%%%%%%%%%%% Proposition n¼ 4.14 %%%%%%%%%%%%%%%%
%%%%%%%%%%%%%%%%%%%%%%%%%%%%%%%%%%%%%%%%%%%%%%%%%%%%%%%%%%
%\newtheorem{defi414}[exam41]{Definition}
\begin{definition}[Difference-$a$ hooks]\label{def4.14}
Suppose ${P}={P}(E)$ and $a\geq 0$; 
we denote by ${\cal
		H}^{a}({P})$ or ${\cal H}^{a}(E)$ the set
 ${\cal H}^{a}({P}) =
\bigcup {\cal H}^{a}({P})_{i}$ where
\begin{align}
	{\cal H}^{a}(E)_{i}= & \{(\mu,\nu)\  |\ 
	\mu,\nu\not\in E, \ y\mu\in E, x\nu\in E,\notag \\\label{e4.18}
  & \text{ degree }\  x^{a}\muÊ= \text{ degree }\,
	\nu,\ =i, \  
\text{	and }\   x^{a}\mu \leq \nu\}. 
\end{align}
\end{definition}
%%%%%%%%%%%%%%%%%%%%%%%%%%%%%%%%%%%%%%%%%%%%%%%%%%%%%%%%%%
These pairs of monomials are the endpoints of
difference-$a$ hooks of  the shape ${ P}(E)$ whose arm
length is $a$ units greater than the leg length. When
$a=0$,  we call such hooks {\em balanced}. When $a=1$, we
write ${\cal H}(E)_{i},\ {\cal H}(E)$, for ${\cal
H}^{1}(E)_{i}, \ {\cal H}^{1}(E)$.\\
%%%%%%%%%%%%%%%%%%%%%%%%%%%%%%%%%%%%%%%%%%%%%%%%%%%%%%%%%%

\begin{figure}[hbtp]
\begin{center}
\leavevmode 
\begin{picture}(25,15)(25,-15)
\setlength{\unitlength}{1mm}
\thicklines
\multiput(0,0)(0,-4){2}{\line(1,0){24}}
\multiput(0,-8)(0,-4){2}{\line(1,0){6}}
\multiput(0,0)(6,0){2}{\line(0,-1){12}}
\multiput(12,0)(6,0){3}{\line(0,-1){4}}
\put(3,-10){\makebox(0,0){{\small \mbox{$\mu$}}}}
\put(21,-2){\makebox(0,0){{\small \mbox{$\nu$}}}}
\end{picture}
\end{center}
\vspace{-.5cm}
\protect\caption{Hook $(\mu,\nu)$ of arm length four, leg length
three.\label{hook-fig}}
\end{figure}
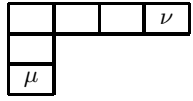

\noindent {\em Figure \ref{hook-fig}:} A {\em hook} in the shape
of ${P}$ is a subshape as shown in Figure \ref{hook-fig}.
consisting of an {\em arm} whos {\em hand} or right
endpoint $\nu$ is the roght endpoint of a row of $
P$, and a {\em leg} whose {\em foot} $\mu$ or lowermost
entry is the endpoint of a column of $ P$. The hook
shown has an {\em arm} of length four, and a {\em leg}
of length three, so has arm-leg {\em difference} equal
one.\\
%%%%%%%%%%%%%%%%%%%%%%%%%%%%%%%%%%%%%%%%%%%%%%%%%%%%%%%%%%
%%%%%%%%%%%%%%%%%%%% hook figure %%%%%%%%%%%%%%%%%%%%%%%%%
%\reline
%\vspace{2.3cm} 
%
%\special{picture hook(fig5) scaled 800}
%\vspace{1mm}
%\underline{Figure 5.} Hook $(\mu,\mu)$ of arm length
%four, leg length three.\\ 
%%%%%%%%%%%%%%%%%%%%%%%%%%%%%%%%%%%%%%%%%%%%%%%%%%%%%%%%%%%

%%%%%%%%%%%%%%%%%%%%%%%%%%%%%%%%%%%%%%%%%%%%%%%%%%%%%%%%%
\noindent We let $E^{\wedge}$  denote the monomial ideal
$\sigma(E)$ obtained by  switching the variables $x, y$ in
$E$. The associated  partition ${ P}(E^{\wedge})$
satisfies     ${P}(E^{\wedge})=(
P(E))^{\wedge}$, the dual of ${ P}(E)$, obtained by
switching rows and columns in  the shape of ${
P}(E)$. L.~ G\"{o}ttsche showed, using somewhat
different language,
%%%%%%%%%%%%%%%%%%%%%%%%%%%%%%%%%%%%%%%%%%%%%%%%%%%%%%%%%
%%%%%%%%%%%%%%%%% Proposition n¼ 4.16 %%%%%%%%%%%%%%%%%%%
%\newtheorem{propos415}[exam41]{Proposition}
	\begin{proposition}[L. G\"{o}ttsche \cite{G2,G4}]\label{prop4.15}		
 		The dimension $z(E)$ of the cell ${\mathbf  Z}(E)$ of
			$ {\mathbf  Z}_T $ satisfies
\begin{equation}\label{e4.19a}
z(E)\  =\  n(T)\  - \ \#{\cal S}(E^{\wedge})\ 
						-\  \#{\cal H}^{0}(E). 
\end{equation}
The dimension $v(E)$ of the  cell
			${\mathbb V}(E)$ of $\G_T$ satisfies
\begin{equation}\label{e4.19b}
v(E)\  =\  n(T)\  -\  f(T)\  -\ 
						\#{\cal S}(E^{\wedge})\  -\   \#{\cal
						H}^{0}(E). 
\end{equation}
\end{proposition}
%%%%%%%%%%%%%%%%%%%%%%%%%%%%%%%%%%%%%%%%%%%%%%%%%%%%%%%%%

\noindent We will now show Proposition \ref{prop4.15} as a
consequence of  Theorem \ref{thm4.9}. We need two results from Part
I \cite{IY2}     \     Let  $\#{\cal H}^{0}(E)_{i}, \ {\cal
S}_{i}(E)$, respectively, denote the subsets of monomial
pairs in ${\cal H}^{0}(E)$, ${\cal S}(E)$  having degree
$i$. 
%%%%%%%%%%%%%%%%%%%%%%%%%%%%%%%%%%%%%%%%%%%%%%%%%%%%%%%%%
%%%%%%%%%%%%%%%%%%%%%%%%%%%%%%%%%%%%%%%%%%%%%%%%%%%%%%%%%
%%%%%%%%%%%%%%%% Lemma 4.17 former 4.16A %%%%%%%%%%%%%%%%
\begin{lemma}\label{lemma4.16A} If $i\geq \mu (T)-1$,
									the number of difference-zero degree-$i$ hooks  $\#({\cal H}^{0}(E)_{i})$
satisfies
\begin{equation}\label{e4.20}
         \#{\cal H}^{0}(E)_{i} = 
									\left( 
																		\begin{array}{c}
																																		\delta_{i}+1\\
																																							2
																		\end{array}
								\right) 
\end{equation}
\end{lemma}
%%%%%%%%%%%%%%%%%%%%%%%%%%%%%%%%%%%%%%%%%%%%%%%%%%%%%%%%%
\begin{lemma}\label{lemma4.16B} 
		The sizes of ${\cal S}_{i}(E)$ and of 	                
		${\cal S}_{i}(E^{\wedge})$ are related by
\begin{equation}\label{e4.21}	
\#{\cal S}_{i}(E^{\wedge}) =
		(\delta_{i}+1)\delta_{i+1} - \#{\cal S}_{i}(E)
\end{equation}
\end{lemma}
%%%%%%%%%%%%%%%%%%%%%%%%%%%%%%%%%%%%%%%%%%%%%%%%%%%%%%%%%
\begin{proof} These are Theorem 1.17, and a
consequence of Theorem 1.8, respectively of Part I \cite{IY2}.
\end{proof}
%%%%%%%%%%%%%%%%%%%%%%%%%%%%%%%%%%%%%%%%%%%%%%%%%%%%%%%%%
\begin{proof}[Proof of Proposition 4.15]
 Since the
morphism          $h : {\mathbf  Z}_T  \rightarrow \G_T$ satisfies 
$h^{-1}({\mathbb V}(E))  =  {\mathbf  Z}(E)$, Theorem \ref{thm4.9}, Propositions
\ref{prop4.13A} and \ref{prop4.13B} imply
that                         
\begin{equation*}
\begin{array}{cl}
 								z(E)	&= \#{\cal S}(E) + f(T)\\
              &= \dim\G_T - \#{\cal S}(E^{\wedge}) +
																		f(T)\\
		    	       &= \dim\G_T - \#{\cal S}(E^{\wedge}) + 
                 (n(T) - \#{\cal H}^{0}(E) -  \dim\G_T)\\
		   	        &= n(T) - \#{\cal S}(E^{\wedge}) -   
                 \#{\cal H}^{0}(E)\, ,
												\end{array}
\end{equation*}
which is G\"{o}ttsche's formula \eqref{e4.19a}.
\end{proof}\par
\begin{remark}\label{remark4.17} 
We can now reconcile the formulas \eqref{e4.17b} and
\eqref{e4.17a}
		for the  fibre dimension of the cell                		
		${\mathbf  Z}(E)$  over ${\mathbb V}(E)$, as a direct
		consequence of Lemmas \ref{lemma4.16A} and \ref{lemma4.16B}. The
size  
		$\#{\cal W}(E)$  of the  fibre counts the size of
		$\{{\cal H}^{a}(E)\ | \ a\neq 0,\ 1, \ -1 \}$,  the
		number of hooks of $E$ having hook-difference any 
		integer but $0$, $1$ and $-1$ . There is
		exactly one hook  having corner $\mu$  for each
		monomial $\mu$ of $E^{c}$ so there are
		$n(T)$ hooks in all. Thus                  			
		$\#{\cal W}(E)$ satisfies
\begin{equation}\label{e4.22}
\begin{tabular}{cl}
		$\#{\cal W}(E)$ &$= n(T) - \#{\cal H}^{0}(E) -  
																					\#{\cal H}^{1}(E) - 										
																						\#{\cal H}^{-1}(E)$\\
                    &$= n(T) -
                      \sum_{i\geq \mu (T)}
																							\delta_{i}(\delta_{i}+1)/2 -  
																							 \sum_{i\geq \mu (T)}
																									(\delta_{i}+1)\delta_{i+1}\,
																										,$\\
\end{tabular}
\end{equation}
 by Lemmas \ref{lemma4.16A},\ref{lemma4.16B},  since $\#{\cal H}^{1}(E) =
\#{\cal S}(E)$, and $\#{\cal H}^{-1}(E) =  \#{\cal
S}(E^{\wedge})$.
\end{remark}

%%%%%%%%%%%%%%%%%%%%%%%%%%%%%%%%%%%%%%%%%%%%%%%%%%%%%%%%%%
%%%%%%%%%%%%%%%%%%%%%%%%%%%%%%%%%%%%%%%%%%%%%%%%%%%%%%%%%%
\subsection{The cells of $\G_T$ and the hook code}\label{hookcode}
%%%%%%%%%%%%%%%%%%%%%%%%%%%%%%%%%%%%%%%%%%%%%%%%%%%%%%%%%%

We recall from Part I \cite{IY2} the ``hook code" for the cells
${\mathbb V}(E)$ of $\G_T$.  The hook code is important,
because it gives a simple way of understanding the Betti
numbers of $\G_T$, and because the homology ring
structure on $H^{\ast}(\G_T)$ in the cases it is
known, is most simply stated in terms of the hook code.
 We show here that the hook code arises in natural way
form the local parameters of the cell ${\mathbb V}(E)$ of
$\G_T$.

Recall that when $T=(1, 2, \ldots , \mu,t_\mu, \ldots
,t_{j},0)$, satisfies the condition $\mu\geq t_\mu\geq
\ldots \geq t_{j}>0$ of \eqref{e4.1}, we let $\mu (T) =\mu$ and
$j(T)=j$. Let ${ P}$ be a partition whose shape has
diagonal lengths $T$. Then $\mu (T)$ is the index of the
 first diagonal of the shape of $P$ - counting from zero -
in which there are spaces, and $j(T)$ is the index of the
last nonempty diagonal.\\
%%%%%%%%%%%%%%%%%%%%%%%%%%%%%%%%%%%%%%%%%%%%%%%%%%%%%%%%%%

\begin{definition}[Partially ordered
sets ${\cal P}(T)$  and ${\cal Q}(T)$].\label{def4.18A}
 Recall that
${\cal P}(T)$ is the set of partitions $ P$ having
diagonal lengths $T$. We denote by ${\cal B}(T)$ the
sequence ${\cal B}(T) = ({\cal B}_{\mu}(T),\ldots \ldots,
{\cal B}_{j}(T))$ where
\begin{align}\notag              
 {\cal B}_{i}(T) = &\text { the partition of 
																								rectangular shape 	having }
																								t_{i}-t_{i+1}\\\label{e4.23}
		   												 				&\text{ parts each of length } 
																							(1+t_{i-1}-t_{i})\text { for }
																	 					\mu\leq i\leq j.
	\end{align}
We denote by ${\cal Q}(T)$ the set of  sequences of
partitions ${\cal Q}_{i}$, with ${\cal Q}_{i}$ {\em
included in}  ${\cal B}_{i}(T)$.  
\begin{equation}\label{e2.24}
{\cal Q}(T) = \{{\cal Q}_{\mu},\ldots ,{\cal
Q}_{i},\ldots ,{\cal Q}_{j})\, \  | \, \  {\cal Q}_{i}
\subset {\cal B}_{i}(T)\}\, . 
\end{equation}
%%%%%%%%%%%%%%%%%%%%%%%%%%%%%%%%%%%%%%%%%%%%%%%%%%%%%%%%%%
%%%%%%%%%%%%%%%%%%%%%%%%%%%%%%%%%%%%%%%%%%%%%%%%%%%%%%%%%%
\noindent If ${\cal Q} = ({\cal Q}_\mu,\ldots
,{\cal Q}_{j}) \in {\cal Q}(T)$, then 
$({\cal Q}_{i})^{c} = {\cal B}_{i}(T)-{\cal Q}_{i}$ is the
complement  of  ${\cal Q}_{i}$ in ${\cal B}_{i}(T)$. We
also let  $\#{\cal Q} = \sum_{i} \#({\cal Q}_{i})$   be
the total  length of the  partitions in the sequence
${\cal Q}$. \\ We partially order ${\cal Q}(T)$ by {\em
inclusion}  of the component partitions, and we partially
order ${\cal P} (T)$ by inclusion of each of the diagonal
partitions ${\cal Q}(E_{i})$, $\mu (T)\leq i\leq j(T)$,
$E=E({ P})$. Both ${\cal Q}(T)$ and ${\cal P} (T)$ are
distributive lattices (see Part I, \cite{IY2}).
\end{definition}
%%%%%%%%%%%%%%%%%%%%%%%%%%%%%%%%%%%%%%%%%%%%%%%%%%%%%%%%%%
%%%%%%%%%%%%%%%%%%%%%%%%%%%%%%%%%%%%%%%%%%%%%%%%%%%%%%%%%%

\begin{definition}[Hook code]\label{def4.18B} Suppose
${ P} = { P}(E)\in {\cal P}(T)$ is a
partition of diagonal lengths $T$  , and  $\mu\leq
i\leq j$. We let ${\cal Q}_{i}({ P})$ be the
partition of the integer $\#{\cal H}^{1}(E)_{i}$ 
determined by the distribution of degree-i difference-one
hooks $(\mu ,\nu)$ according to the hand monomial
$\nu$.  We order the $\delta_{i+1} = t_{i} -
t_{i+1}$ monomials of $(E^{c})_{i}$ such
that $x\nu \in E_{i+1}$ in the reverse alphabetic
order. Then
\begin{equation}\label{e4.25a}
 {\cal Q}_{i}({P}) = (q_{i1}, \ldots
,q_{i\delta_{i+1}}), q_{ik} = \#\{ (\mu , \nu) \in {\cal
H}^{1}(E)_{i} \ | \ \nu=\nu_{k} \}. 
\end{equation}
\noindent The partition ${\cal Q}_{i}({P})$ is considered to be inside the box $B_i(T)$,
so it has
$\delta_{i+1}$ parts, some possibly zero. The {\em hook code} ${\cal D}({ P})$ of ${
P}$ is the sequence 
\begin{equation}\label{e4.25b}
{\cal D}({ P}) = ({\cal Q}_\mu (P), \ldots , {\cal
Q}_{j}({ P})) \subset {\cal B}(T).
\end{equation}
The length $\ell ({\cal D}({ P})) = \ell ({\cal
Q}_\mu({ P})) + \ldots + \ell ({\cal Q}_{j}({
P}))$, and is $\#{\cal H}^{1}(E)$. 
\end{definition}
%%%%%%%%%%%%%%%%%%%%%%%%%%%%%%%%%%%%%%%%%%%%%%%%%%%%%%%%%%
%%%%%%%%%%%%%%%%%%%%%%%%%%%%%%%%%%%%%%%%%%%%%%%%%%%%%%%%%%
\begin{theorem}\label{thm4.19} 
 The coding map $\cal{ P} \rightarrow {\cal D}({ P})$ 
	gives an isomorphism of distributive lattices \linebreak ${\cal
D}: {\cal P}(T)\rightarrow
	{\cal Q}(T)$, such that
	\begin{enumerate}
 \item[i.] The dimension and codimension of the cell                
            ${\mathbb V}(E)$ satisfies                  
\begin{align}\label{edimTcell}
 \dim {\mathbb V}(E) &=  \ell ({\cal D}(P(E)))\\
\cod \mathbb V(E)&=\ell ({\cal D}({ P(E)}^{\wedge})) = \#{\cal
H}^{-1}(E).\label{ecodTcell}
\end{align}
\item[ii.] ${\cal Q}({ P}^{\wedge}) = ({\cal Q}({
P}))^{c}$.\\
\end{enumerate} 
\end{theorem}
%%%%%%%%%%%%%%%%%%%%%%%%%%%%%%%%%%%%%%%%%%%%%%%%%%%%%%%%%%
\begin{proof}[Proof of Theorem \ref{thm4.19}] The dimension
formula in    Part i. of the Theorem is immediate from
Theorem \ref{thm4.9}  and the definition of $\cal D$. Part ii. of
the Theorem is the case $a = 1$ of Theorem 1.8 of Part I \cite{IY2}. 
That $\cal D$ is an isomorphism of distributive lattices
is Theorem 2.28 of \cite{IY2}.
\end{proof}
%%%%%%%%%%%%%%%%%%%%%%%%%%%%%%%%%%%%%%%%%%%%%%%%%%%%%%%%%%
%%%%%%%%%%%%%%%%%%%%%%%%%%%%%%%%%%%%%%%%%%%%%%%%%%%%%%%%%%              
{\bf Remark.} J. Yam\'{e}ogo in [Y3] outlines a different
proof of Theorem \ref{thm4.19}. He uses a parametrized family of
linkings to reduce from the cell $\mathbb{V}(E)$ to the cell
$\mathbb{V}((x^{a},y^{b}):E)$ whose partition ${P}'$ is
the complement of ${P}(E)$ in a $b\times a$ rectangle.\\
%%%%%%%%%%%%%%%%%%%%%%%%%%%%%%%%%%%%%%%%%%%%%%%%%%%%%%%%%% 
%%%%%%%%%%%%%%%%%%%%%%%%%%%%%%%%%%%%%%%%%%%%%%%%%%%%%%%%%%
{\bf Warning.} The hook partion ${\cal Q}_{i}(E)$ is
quite different from the  partition ${ Q}(E_{i})$
that we defined above. Neither determines the other,
except in special cases such as $i=\mu (T)=
j(T)$.\footnote{Here are two such special cases, where
${\cal Q}_{i}$ determines ${ Q}_{i}$. When $T=(1,2,
\ldots , \mu (T), t,0)$ then $i=j(T)=\mu (T), \G_T =
\Grass(i+1-t,i+1)$ as in Section \ref{sec2}; and ${\cal Q}_{i}
={ Q}_{i}$. When $T=(1,2, \ldots , \mu (T), a, \ldots
,a)$ then ${\cal Q}_{j(T)}$ determines ${ Q}_\mu,
{ Q}_{\mu+1}, \ldots ,{ Q}_{j}$ and any ${
Q}_{i}$ determines ${\cal Q}_{j(T)}$.} \ However, the
sequence of partitions $({\cal Q}_{\mu}(E), \ldots {\cal
Q}_{j}(E))$ determines $({ Q}(E_{\mu}), \ldots ,{
Q}(E_{j}))$ and vice versa. The connection between the
sequences is explained in \S 2C of Part I: see Theorem
2.23ff and Example 2.31 of \cite{IY2}.\\

%%%%%%%%%%%%%%%%%%%%%%%%%%%%%%%%%%%%%%%%%%%%%%%%%%%%%%%%%%
%%%%%%%%%%%%%% Insert Table 2 figure here %%%%%%%%%%%%%%%%
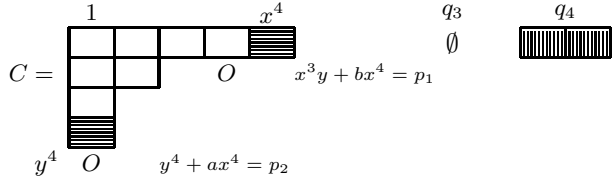
\begin{figure}[th]
\begin{center}
\leavevmode 
\begin{picture}(80,100)(6,-100)
\setlength{\unitlength}{1mm}
\thicklines
\multiput(6,-4)(0,-4){2}{\line(1,0){30}}
\multiput(6,-56)(0,-4){2}{\line(1,0){30}}
\put(6,-12){\line(1,0){12}}
\multiput(6,-16)(0,-4){2}{\line(1,0){6}}
\multiput(6,-64)(0,-4){2}{\line(1,0){12}}
\multiput(66,-4)(0,-4){2}{\line(1,0){12}}
\multiput(66,-56)(0,-4){2}{\line(1,0){12}}
\multiput(54,-56)(0,-4){2}{\line(1,0){6}}
\multiput(6,-4)(6,0){2}{\line(0,-1){16}}
\put(18,-4){\line(0,-1){8}}
\multiput(24,-4)(6,0){3}{\line(0,-1){4}}
\multiput(6,-56)(6,0){3}{\line(0,-1){12}}
\multiput(24,-56)(6,0){3}{\line(0,-1){4}}
\multiput(66,-4)(6,0){3}{\line(0,-1){4}}
\multiput(66,-56)(6,0){3}{\line(0,-1){4}}
\multiput(54,-56)(6,0){2}{\line(0,-1){4}}
\thinlines
\multiput(30,-4.5)(0,-.5){7}{\line(1,0){6}}
\multiput(6,-16.5)(0,-.5){7}{\line(1,0){6}}
\multiput(12,-64.5)(0,-.5){7}{\line(1,0){6}}
\multiput(30,-56.5)(0,-.5){7}{\line(1,0){6}}
\multiput(66.5,-8)(.5,0){11}{\line(0,1){3.5}}
\multiput(72.5,-8)(.5,0){11}{\line(0,1){3.5}}
\multiput(54.5,-60)(.5,0){11}{\line(0,1){3.5}}
\multiput(66.5,-60)(.5,0){11}{\line(0,1){3.5}}
\multiput(72.5,-60)(.5,0){11}{\line(0,1){3.5}} 
\put(57,-6){\makebox(0,0){$\emptyset$}}
\put(9,-2){\makebox(0,0){{\small \mbox{$1$}}}}
\put(33,-2){\makebox(0,0){{\small \mbox{$x^{4}$}}}}
\put(57,-2){\makebox(0,0){{\small \mbox{$q_{3}$}}}}
\put(72,-2){\makebox(0,0){{\small \mbox{$q_{4}$}}}}
\put(1,-10){\makebox(0,0){{\small \mbox{$C=$}}}}
\put(1,-62){\makebox(0,0){{\small \mbox{$D=$}}}}
\put(27,-10){\makebox(0,0){{\small \mbox{$O$}}}}
\put(36,-11){{\scriptsize \mbox{$x^{3}y+bx^{4} = p_{1}$}}}
\put(3,-22){\makebox(0,0){{\small \mbox{$y^{4}$}}}}
\put(9,-22){\makebox(0,0){{\small \mbox{$O$}}}}
\put(18,-23){{\scriptsize \mbox{$y^{4}+ax^{4} = p_{2}$}}}
\put(27,-62){\makebox(0,0){{\small \mbox{$O$}}}}
\put(36,-63){{\scriptsize \mbox{$x^{3}y+a_{1}x^{4} = p_{1}$}}}
\put(9,-70){\makebox(0,0){{\small \mbox{$O$}}}}
\put(15,-70){\makebox(0,0){{\small \mbox{$O$}}}}
\put(24,-71){{\scriptsize \mbox{$xy^{3}+a_{2}x^{4} = p_{2}$}}}
\put(8,-75){{\scriptsize \mbox{$y^{3}+bxy^{2}+\cdots = p_{3}$}}}
\put(30,-45){\shortstack[l]{{\small \mbox{The forms $p_{1}$,
$p_{2}$  in $I$ determine}}
\\
{\small \mbox{$f=x^2y+bx^3, \ yf=x^2y^2+byx^3$}}
\\
{\small \mbox{$g=xy^2-b^2x^3, \ yg=xy^3-b^2yx^3$ in $I$}}
\\
\vspace{1.2mm}
\\
{Cell ${\Bbb V}(C), \ \ C=(x^{5},x^{2}y,xy^{2},y^{4})$}}}
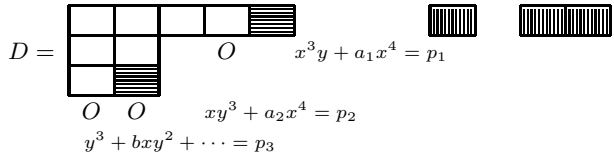
\put(30,-92){\shortstack[l]{{\small \mbox{The coefficients
$a_{1}$, $a_{2}$, $b$}}
\\
{\small \mbox{determine $I$ in the cell ${\mathbb V}(D)$}}
\\
\vspace{1.2mm}
\\
{ Cell ${\mathbb V}(D), \ \ D=(x^{5},x^{2}y,y^{3})$}
}}
\end{picture}
\end{center}
\vspace{-0.5cm}
\renewcommand{\figurename}{Table}
\protect\caption{Two cells of $\G_T $, $T=(1,2,3,2,1)$ and their
codes. See Example \ref{ex4.20}.}\label{twocells}
\end{figure}

%\vspace{9.5cm}
%
%\special{picture table2(fig) scaled 800}     
%
%%%%%%%%%%%%%%%%%%%%%%%%%%%%%%%%%%%%%%%%%%%%%%%%%%%%%%%%%%%
%%%%%%%%%%%%%%%%%%%%%%%%%%%%%%%%%%%%%%%%%%%%%%%%%%%%%%%%%%%
%\medskip
%%%%%%%%%%%%%%%%%%%%%%%%%%%%%%%%%%%%%%%%%%%%%%%%%%%%%%%%%%
\noindent
%%%%%%%%%%%%%%%%%%%%%%%%%%%%%%%%%%%%%%%%%%%%%%%%%%%%%%%%%%   
\begin{example}\label{ex4.20} Table \ref{twocells} depicts several
cells in
$\G_T$, $T = (1, 2, 3, 2, 1)$ with their parameters and
corresponding pruning code.\par
The first cell, labeled $C$, corresponds to the partition
${P}(C) = (5, 2, 1, 1)$, and is important for a
counterexample of J.~Yam\'eogo (see (Ci) in \S \ref{sec4F} below). In
Table \ref{twocells} the shape or Ferrer's graph of ${P}$ is shaded, and
corresponds  to the cobasis $(1, x, x^{2}, x^{3}, x^{4}; y, yx; y^{2};
y^{3})$. 
%%%%%%%%%%%%%%%%%%%%%%%%%%%%%%%%%%%%%%%%%%%%%%%%%%%%%%%%%%
The cell $\mathbb{V}(C)$ consists of all ideals $I$ in $R$
having the monomial initial ideal $C = (x^{5}, x^{2}y,
xy^{2}, y^{4})$.  If $I$  is in the cell,  the two
key elements $p_{1} = x^{3}y+bx^{4}$ and $p_{2} =
y^{4}+ax^{4}$ of $I$ determine the $I$ and are listed at
right. They have the form  $(\mu + c\nu)$, where $c\in k$,
where  $\nu$ Ð darkly shaded Ð is a monomial at the {\em
right endpoint} of a row , and $\mu$ Ð indicated by a
space with ``O" inside Ð is a monomial {\em just below the
shape, and  lying on the same diagonal as} $\nu$ {\em but
lower than} $\nu$. There are two such pairs $(\mu,\nu)$
with $\nu = x^{4}$; the two pairs count as the unique part
of the code partition $q_{4} = (2)$. There are no such
pairs with $\nu = y^{3}$, and we have written at right
$q_{3} = (\emptyset)$, the empty partition. \par
The second cell in Table \ref{twocells} corresponds to
the partition ${P} = (5, 2, 2)$. The ideal $I$ is in
the cell if its monomial initial ideal is $(x^{5}, yx^{2},
y^{3})$; if so, the ideal is completely determined by the
coefficients $a_{1}, a_{2}, b$ of the forms $p_{1}, p_{2},
p_{3}$.  At right we have written the code partition
$q_{4} = (2)$ corresponding to the pairs of monomials\\ 
\begin{tabular}{ll}
	$(\mu,\nu)$ & = (initial monomial, cobasis monomial) \\
		&  = $(x^{3}y, x^{4}), (xy^{3}, x^{4})$ of the
coefficients $a_{1}, a_{2}$.
\end{tabular}\\
%%%%%%%%%%%%%%%%%%%%%%%%%%%%%%%%%%%%%%%%%%%%%%%%%%%%%%%%%%
The code partition $q_{3} = (1)$ of the second cell
corresponds to the pair $(y^{3}, y^{2}x)$ of the
coefficient $b$ in $p_{3}$.\\
\end{example}
\noindent
%%%%%%%%%%%%%%%%%%%%%%%%%%%%%%%%%%%%%%%%%%%%%%%%%%%%%%%%%%Ê
%%%%%%%%%%%%%%%%%%%%%%%%%%%%%%%%%%%%%%%%%%%%%%%%%%%%%%%%%%%
%\underline{Table 2. Two cells of $\G_T$, $T = (1,
%2, 3, 2, 1)$ and their codes}\\  

%%%%%%%%%%%%%%%%%%%%%%%%%%%%%%%%%%%%%%%%%%%%%%%%%%%%%%%%%%

%%%%%%%%%%%%%%%%%%%%%%%%%%%%%%%%%%%%%%%%%%%%%%%%%%%%%%%%%%
%%%%%%%%%%%%%%%%%%% Sections 5-6-7 %%%%%%%%%%%%%%%%%%%%%%
%%%%%%%%%%%%%%%%%%%%%%%%%%%%%%%%%%%%%%%%%%%%%%%%%%%%%%%%%%
\subsection{The homology groups $H^{\ast}(\G_T)$}\label{homgroups}
We apply Theorem \ref{thm4.19} to determining the Betti numbers of
the variety $\G_T$.
%%%%%%%%%%%%%%%%%%%%%%%%%%%%%%%%%%%%%%%%%%%%%%%%%%%%%%%%%%
First, recall that the generating function $B(a,b;q) = 
\sum_{n} p(a,b;n)q^{n}$ for the number $p(a,b;n)$ of
partitions of $n$ into at most $b$ parts each less or equal to $a$ is
the $q$-binomial coefficient,%
%%%%%%%%%%%%%%%%%%%%%%%%%%%%%%%%%%%%%%%%%%%%%%%%%%%%%%%%%%
\begin{equation}\label{e5.1}
 \sum_{n} p(a,b,n)q^{n} = \left[ \begin{array}{c} a+b\\
a \end{array} \right] =
\frac{(q)_{a+b}}{(q)_{a}(q)_{b}},\, \ where\, \ (q)_{a} =
(1-q^{a})(1-q^{a-1})\cdots (1-q)\, . 
\end{equation} 
%%%%%%%%%%%%%%%%%%%%%%%%%%%%%%%%%%%%%%%%%%%%%%%%%%%%%%%%%% 
The polynomial $B(a,b;q^{2})$ is also the Poincar\'{e}
polynomial  $B(\Grass(a,a+b))$ of the Grassmannian of
$a$-subspaces of a complex vector space of dimension
$a+b$.
%%%%%%%%%%%%%%%%%%%%%%%%%%%%%%%%%%%%%%%%%%%%%%%%%%%%%%%%%%

A partition $ P$ of diagonal lengths $T$
corresponds to a unique monomial ideal $E({ P})$,
defining a quotient $R/E({ P})$ of Hilbert function
$T$; and conversely, a monomial ideal  $E$ determines a
unique partition ${ P} = { P}(E)$
%%%%%%%%%%%%%%%%%%%%%%%%%%%%%%%%%%%%%%%%%%%%%%%%%%%%%%%%%% 
(see Definition \ref{def4.4}).
%%%%%%%%%%%%%%%%%%%%%%%%%%%%%%%%%%%%%%%%%%%%%%%%%%%%%%%%%%

We let $b(T,h)$ denote the number of partitions ${P}$
of $n = \sum t_{i}$ having diagonal lengths $T$ and for
which the corresponding cell $\mathbb{V}(E({P}))$ in
$\G_T$ has dimension $h$, or, equivalently, for which 
${P}$ has $h$ difference-one hooks (Theorem \ref{thm4.9}).
%%%%%%%%%%%%%%%%%%%%%%%%%%%%%%%%%%%%%%%%%%%%%%%%%%%%%%%%%% 

In the following result, we consider $\G_T$ as complex
variety,, because of our use of the theorem of Bialynicki-Birula. The results concerning
cellular decomposition of $\G_T$ are valid in characteristic zero, or $\cha k= p>j$, and
the results below will count the cells of given dimension for such fields $k$.
%%%%%%%%%%%%%%%%%%%%%%%%%%%%%%%%%%%%%%%%%%%%%%%%%%%%%%%%%% 
We denote by $H^{i}(\G_T)$ the homology group in
codimension $i$. Note that $\Grass(\delta_{i+1},1+\delta_i+\delta_{i+1})=\Grass(t_{i}-t_{i+1},
1+t_{i-1}-t_{i})$. Over the complexes the cells occur in even (real) dimensions, and we
denote by $B(T,q^2)$ the Poincar\'{e} polynomial of $\G_T$, $B(T,q^2)=\sum \dim
H^i(\G_T)q^i$.
%%%%%%%%%%%%%%%%%%%%%%%%%%%%%%%%%%%%%%%%%%%%%%%%%%%%%%%%%%
%\newtheorem{theo51}{Theorem}[section]
\begin{theorem}[G.~Gotzmann in codimension one
\cite{gotz2}]\label{thm5.1} Suppose that $k = {\bf C}$. 
The isomorphism $\cal D$ of 
Theorem \ref{thm4.19} induces an additive isomorphism of homology groups
$\tau: H^{\ast}(\G_T)\longrightarrow   
H^{\ast}(\SGrass(T))$, over $\mathbb Z$, where $\SGrass(T)$ is the product of small
Grassmannians (Definition
\ref{def4.10}):

\begin{equation}\label{e5.2}
\tau : \, H^{\ast}(\G_T) 
\cong_{add} 
\prod_{\mu (T)\leq i\leq j(T)}
H^{\ast}(\Grass(\delta_{i+1}, 1+\delta_i+\delta_{i+1}))
\end{equation}
The homomorphism $\tau$ respects the $\mathbb Z_2$ duality action.
\end{theorem}
%%%%%%%%%%%%%%%%%%%%%%%%%%%%%%%%%%%%%%%%%%%%%%%%%%%%%%%%%%
%\newtheorem{theo52}[theo51]{Theorem}
\begin{theorem}\label{thm5.2} 
When $k=\mathbb C$ the Poincar\'{e} polynomial 
$B(\G_T) = B(T;q^{2})$ satisfies
\begin{equation}\label{e5.3}
B(T; q^{2}) = 
\prod_{\mu\leq i\leq j}
B(\Grass(\delta_i, 1+\delta_i+\delta_{i+1})=\prod_{\mu\le i\le
j}B(1+\delta_i,\delta_{i+1},q^2).
\end{equation}
When $\cha k=0$, or $\cha k>j(T)$ the coefficent of $q^{2u}$ in \eqref{e5.3} counts the
number of cells $\mathbb V(E)$ having $k$-dimension $2u$; and the coefficient counts
also the number of partitions of diagonal lengths $T$ having $u$ hooks of difference
one. \par
The total number of cells $\mathbb V(E)\subset \G_T$ is
$b(T) = B(T;1)$, which is the following product of binomial
coefficients
\begin{equation}\label{e5.4}
b(T) = 
\prod_{\mu (T)\leq i\leq j(T)}
\left (
\begin{array}{c}
1+t_{i-1}-t_{i+1}\\
t_{i}-t_{i+1}
\end{array}
\right )
\end{equation}
The product $b(T)$ also counts the number of monomial
ideals $E$ in ${k}[x,y]$  defining a quotient algebra of
Hilbert function $T$; and when $k=\mathbb C$ it is also the Euler
characteristic  $\chi(\G_T)$.\\
\end{theorem}
%%%%%%%%%%%%%%%%%%%%%%%%%%%%%%%%%%%%%%%%%%%%%%%%%%%%%%%%%%
\begin{proof}
%%%%%%%%%%%%%%%%%%%%%%%%%%%%%%%%%%%%%%%%%%%%%%%%%%%%%%%%%%
By Theorem \ref{thm4.9} the cells $\{ \mathbb{V}(E)$ such that $ E$ is a 
monomial ideal with $H(R/E) = T\}$ 
%%%%%%%%%%%%%%%%%%%%%%%%%%%%%%%%%%%%%%%%%%%%%%%%%%%%%%%%%%
form a cellular decomposition of $\G_T$, in any
characteristic $p > j$, or in characteristic $0$.
%%%%%%%%%%%%%%%%%%%%%%%%%%%%%%%%%%%%%%%%%%%%%%%%%%%%%%%%%%
 When $k = {\mathbb C}$, it follows from \cite{B} that the classes of the
cells form a ${\mathbb Z}$-basis of the homology
$H^{\ast}(\G_T)$.
%%%%%%%%%%%%%%%%%%%%%%%%%%%%%%%%%%%%%%%%%%%%%%%%%%%%%%%%%%
Theorems~\ref{thm5.1} and \ref{thm5.2} now follow from Theorem \ref{thm4.19} and
the well-known homology of $\Grass (a,b)$.
\end{proof}

%%%%%%%%%%%%%%%%%%%%%%%%%%%%%%%%%%%%%%%%%%%%%%%%%%%%%%%%%%
%%%%%%%%%%%%%%%%%% section 6 now 4F %%%%%%%%%%%%%%%%%%%%%%%%%%%%%
\subsection{Do the classical results for the
									Grassmannian extend to $\G_T$?}\label{sec4F}
%%%%%%%%%%%%%%%%%%%%%%%%%%%%%%%%%%%%%%%%%%%%%%%%%%%%%%%%%%
In Section \ref{classicalresults} we listed the classical results for
the Grassmannian $\Grass (d,R_{j})$. Here we state which extend
to $\G_T$.\\
%%%%%%%%%%%%%%%%%%%%%%%%%%%%%%%%%%%%%%%%%%%%%%%%%%%%%%%%%%
%%%%%%%%%%%%%%%%%%%%%%%%%%%%%%%%%%%%%%%%%%%%%%%%%%%%%%%%%%

\par {\bf A.} $\G_T$ has a cellular decomposition
whose  cells $\mathbb V(E)$ correspond $1-1$ with initial monomial ideals
$E$, such that the Hilbert function $H(R/E)=T$, and
$I\in {\mathbb V}(E)\rightarrow \In (I)=E\ .$ (Theorem \ref{thm4.9}).
%%%%%%%%%%%%%%%%%%%%%%%%%%%%%%%%%%%%%%%%%%%%%%%%%%%%%%%%%%
The cells ${\mathbb V}(E)$ correspond $1-1$ with the 
partitions  ${ P}(E)$\  of\   $n=\sum t_{i}$    
having diagonal lengths $T$. 
%%%%%%%%%%%%%%%%%%%%%%%%%%%%%%%%%%%%%%%%%%%%%%%%%%%%%%%%%% 
The shape of the partition ${ P}(E)$ can be identified 
with the standard cobasis $E^{c}$ for the ideal $E$. (See Proposition \ref{prop4.5},
Theorem \ref{thm4.9})\\
%%%%%%%%%%%%%%%%%%%%%%%%%%%%%%%%%%%%%%%%%%%%%%%%%%%%%%%%%%

\par {\bf B.} The set of cells ${\mathbb V}(E)$ of
$\G_T$  correspond $1-1$ to the set of all
sequences ${\cal Q} = ({\cal Q}_\mu, \ldots , {\cal
Q}_{j})$ where each ${\cal Q}_{i}$ is an ordinary
partition having no more than $(t_{i}-t_{i+1})$ nonzero
parts, each of size no greater than $(1+t_{i-1}-t_{i})$.
Thus, each ${\cal Q}_{i}({ P})$ is a subpartition of a
$(t_{i}-t_{i+1})\times (1+t_{i-1}-t_{i})$ rectangle ${\cal
B}_{i}(T)$. (See Theorem \ref{thm4.19}.)\\
%%%%%%%%%%%%%%%%%%%%%%%%%%%%%%%%%%%%%%%%%%%%%%%%%%%%%%%%%%
%%%%%%%%%%%%%%%%%%%%%%%%%%%%%%%%%%%%%%%%%%%%%%%%%%%%%%%%%%

\par {\bf C.}
%%%%%%%%%%%%%%%%%%%%%%%%%%%%%%%%%%%%%%%%%%%%%%%%%%%%%%%%%%
The dimension of ${\mathbb V}(E) $ satisfies $\dim \mathbb V(E)=(\ell \cal Q
(P(E))=\#\cal H (P(E))$, the
 total length of the  partitions comprising the code $\mathcal Q(P(E))$, or the
total number of difference one hooks of $P(E)$.
%%%%%%%%%%%%%%%%%%%%%%%%%%%%%%%%%%%%%%%%%%%%%%%%%%%%%%%%%
	The code of the dual partition ${ P}^{\wedge}$ is the
complement  of ${\cal Q}({ P})$:
%%%%%%%%%%%%%%%%%%%%%%%%%%%%%%%%%%%%%%%%%%%%%%%%%%%%%%%%%%
${\cal Q}({ P}^{\wedge}) = ({\cal Q}(E)^{c})$ (Theorem \ref{thm4.19}).\\
%%%%%%%%%%%%%%%%%%%%%%%%%%%%%%%%%%%%%%%%%%%%%%%%%%%%%%%%%%

\par {\bf Ci.}
%%%%%%%%%%%%%%%%%%%%%%%%%%%%%%%%%%%%%%%%%%%%%%%%%%%%%%%%%%
Dimensional properness fails when          
%%%%%%%%%%%%%%%%%%%%%%%%%%%%%%%%%%%%%%%%%%%%%%%%%%%%%%%%%%
$T=(1, 2, 3, 2, 1)$  and the cell is $C$ of codimension 
$2$ in $\G_T$  from Table \ref{twocells}. 
%%%%%%%%%%%%%%%%%%%%%%%%%%%%%%%%%%%%%%%%%%%%%%%%%%%%%%%%%%
The family of ideals
%%%%%%%%%%%%%%%%%%%%%%%%%%%%%%%%%%%%%%%%%%%%%%%%%%%%%%%%%%
$I_{a} = (y^{4}+ax^{4}, yx^{2}, y^{2}x), a \neq 0 $ 
%%%%%%%%%%%%%%%%%%%%%%%%%%%%%%%%%%%%%%%%%%%%%%%%%%%%%%%%%%
has initial ideal $(y^{4}, yx^{2}, y^{2}x)$ corresponding 
%%%%%%%%%%%%%%%%%%%%%%%%%%%%%%%%%%%%%%%%%%%%%%%%%%%%%%%%%%
to the   cell $C$ in the basis $(x,y)$ where $y<x$, and
has the corresponding initial ideal
$(x^{4}, xy^{2}, x^{2}y)$ in the basis $(y,x),\, \  y > x$.
%%%%%%%%%%%%%%%%%%%%%%%%%%%%%%%%%%%%%%%%%%%%%%%%%%%%%%%%%%
Thus, the $1$ dimensional family $I_{a}$  
satisfies
%%%%%%%%%%%%%%%%%%%%%%%%%%%%%%%%%%%%%%%%%%%%%%%%%%%%%%%%%%
$$I_{a}\subset \overline{{\mathbb V}(C,x)}\,\cap
\,\overline{{\mathbb V}(C,y)}\ ,$$                          
%%%%%%%%%%%%%%%%%%%%%%%%%%%%%%%%%%%%%%%%%%%%%%%%%%%%%%%%%%
so the intersection is not proper (Example of
J.~ Yam\'{e}ogo \cite{Y1}).\\
 When the
intersection  is proper, the homology class of intersection $\overline{{\mathbb
V}(E,p)}\,\cap
\,\overline{{\mathbb V}(E',p')}$ of several ramification
conditions at distinct points can be calculated from the still unknown homology ring
of $\G_T$. See  Theorem \ref{thm7.0}, Example
\ref{ex7.4} and the Problem at the end of  Section \ref{sec5B}. In Section \ref{sec5A}
we list the cases where the ring structure $H^\ast(\G_T)$ is known.
%%%%%%%%%%%%%%%%%%%%%%%%%%%%%%%%%%%%%%%%%%%%%%%%%%%%%%%%%%
%%%%%%%%%%%%%%%%%%%%%%%%%%%%%%%%%%%%%%%%%%%%%%%%%%%%%%%%%%
\par {\bf Cii.}
%%%%%%%%%%%%%%%%%%%%%%%%%%%%%%%%%%%%%%%%%%%%%%%%%%%%%%%%%%
The analog of the Eisenbud-Harris result concerning specialization of intersections on
$\Grass(d,R_j)$ (see \S \ref{classicalresults} (Cii)) is still open.\\
%%%%%%%%%%%%%%%%%%%%%%%%%%%%%%%%%%%%%%%%%%%%%%%%%%%%%%%%%%
%%%%%%%%%%%%%%%%%%%%%%%%%%%%%%%%%%%%%%%%%%%%%%%%%%%%%%%%%%

\par {\bf D.} The Poincar\'{e} duality is not exact,
when expressed in terms of  the cells. J.~Yam\'{e}ogo shows 
that the  intersection $\overline{C}\cdot \overline{C} =\,
[ {\mathrm{ point}}]\, = \overline{C}\cdot \overline{C^{\wedge}}$
when  $T = (1,2,3,2,1)$.
%%%%%%%%%%%%%%%%%%%%%%%%%%%%%%%%%%%%%%%%%%%%%%%%%%%%%%%%%%
However, he has also shown that the duality  is always exact
in codimension one (See \cite{Y6}).\\
%%%%%%%%%%%%%%%%%%%%%%%%%%%%%%%%%%%%%%%%%%%%%%%%%%%%%%%%%%
%%%%%%%%%%%%%%%%%%%%%%%%%%%%%%%%%%%%%%%%%%%%%%%%%%%%%%%%%%

\par {\bf E.}\, J.~Yam\'{e}ogo has shown that the
closure  $\overline{C}$ is  not a union of cells.
%%%%%%%%%%%%%%%%%%%%%%%%%%%%%%%%%%%%%%%%%%%%%%%%%%%%%%%%%%
He proves in [Y-5]  the following result:  
\begin{uthm} 
If $U_{c}$  is a cell in $\G_T$
of dimension $c$ then  there are cells $\overline{U_{+}}\supset
																				\overline{U}\supset \overline{U_{-}}$\ ,
																				where  $\dim U_{+}=c+1$
																				 and $\dim U_{-} = c-1\ .$
\end{uthm}
%%%%%%%%%%%%%%%%%%%%%%%%%%%%%%%%%%%%%%%%%%%%%%%%%%%%%%%%%%

%%%%%%%%%%%%%%%%%%%%%%%%%%%%%%%%%%%%%%%%%%%%%%%%%%%%%%%%%%
\noindent His method is to use a parametrized family of
linkings   to reduce the question to cells of lower
dimension.
%%%%%%%%%%%%%%%%%%%%%%%%%%%%%%%%%%%%%%%%%%%%%%%%%%%%%%%%%

%%%%%%%%%%%%%%%%%%%%%%%%%%%%%%%%%%%%%%%%%%%%%%%%%%%%%%%%%
%%%%%%%%%%%%%%%% Chain conjecture %%%%%%%%%%%%%%%%%%%%%%% 
\begin{conjecture}[Chain conjecture]
The cells $U_{c}$, $U_{e}$ of dimensions $c,\, e$ 
in $\G_T$  satisfy $\overline{U_{c}}\supset U_{e}$  
iff there is  a chain of  cells
$\, \overline{U_{c}}\supset \overline{U_{c+1}}\supset 
\cdots \supset \overline{U_{e}}$ such that $U_{i+1}$ has
codimension $1$ in $\overline{U_{i}}$.
\end{conjecture}  
%%%%%%%%%%%%%%%%%%%%%%%%%%%%%%%%%%%%%%%%%%%%%%%%%%%%%%%%%%
%%%%%%%%%%%%%%%%%%%%%%%%%%%%%%%%%%%%%%%%%%%%%%%%%%%%%%%%%%
%%%%%%%%%%%%%%%%%%%%%% section n¼ 7 %%%%%%%%%%%%%%%%%%%%%%
\section{Homology ring of $\G_T$ and the intersection of
cells}\label{sec5}
  In Section \ref{sec5A} we note the relevence of the homology ring $H^\ast(\G_T)$ in
computing the intersection of ramification loci at different points $p\in
\check{\mathbb P}^1$ (Theorem \ref{thm7.0}), and we summarize what is known
about the homology ring. In Section \ref{sec5B} we determine the homology ring
for $\G_T$ when $T=T(\mu ,j)=(1,2,\ldots \mu, \ldots \mu, 1=t_j)$ (Theorem
\ref{thm7.3}). In Section \ref{sec5C} we show that this $\G_T, T=T(\mu ,j)$ is a
desingularisation of the $\mu$-secant variety $\Sec (\mu ,j)$ to the degree-$j$
rational normal curve, and we identify the homology classes of the pullbacks to
$H^\ast (G_T)$ of the higher singular loci of $\Sec (\mu ,j)$ (Theorem \ref{thm8.0}).
%%%%%%%%%%%%%%%%%%%%%%%%%%%%%%%%%%%%%%%%%%%%%%%%%%%%%%%%%%
%%%%%%%%%%%%%%%%%%% Section n¼ 7A %%%%%%%%%%%%%%%%%%%%%%%%
\subsection{Intersection of ramification conditions}\label{sec5A}
%%%%%%%%%%%%%%%%%%%%%%%%%%%%%%%%%%%%%%%%%%%%%%%%%%%%%%%%%%
Knowledge of the homology ring $ H^{\ast}(\G_T)$
would allow us to understand  the homology class of the
intersection of ramification conditions 
%%%%%%%%%%%%%%%%%%%%%%%%%%%%%%%%%%%%%%%%%%%%%%%%%%%%%%%%%%
$\overline{{\mathbb V}(E;x_{0})}\,\cap \,\overline{{\mathbb
V}(E',x_{1})}$ on ideals of linear systems at different
points $x_{0},\, x_{1}$  of the curve ${{\mathbb P}}^{1}$.       
%%%%%%%%%%%%%%%%%%%%%%%%%%%%%%%%%%%%%%%%%%%%%%%%%%%%%%%%%%
We give an example of this in \S \ref{sec5B}  (Example \ref{ex7.4}).
Although by Theorem \ref{thm5.1} the homology of  $\G_T$ is
additively isomorphic to that of a product $\SGrass(T)$ of
small Grassmanians, the ring structure is  not in general
isomorphic to that of the product. Nevertheless there are
some simplifications. We first formulate the basic
result, motivating our study of the homology ring of
$\G_T$. We suppose that ${k} = {\mathbb C}$ in this
section.\\
%%%%%%%%%%%%%%%%%%%%%%%%%%%%%%%%%%%%%%%%%%%%%%%%%%%%%%%%%%
\begin{theorem}\label{thm7.0} Let $p_1,p_2,\ldots , p_s$ be $s$ points of
${\check{\mathbb P}}^1$, let ${\mathbb V}(E_1,p_1), \ldots , {\mathbb V}(E_s,p_s)
\subset \G_T$ be the cells associated to $\QRAM (V,p_1) = \QRAM (E_1),\ldots
\QRAM (V,p_s) = \QRAM (E_s)$, and suppose that the intersection $\Z =
\overline{{\mathbb V}(E,p)}\cap \overline {{\mathbb V}(E',p')}$ is proper. Then
the homology class of $\Z$ is the intersection product of
the classes of $\overline {{\mathbb V}(E_1,p_1)},\ldots  ,
\overline {{\mathbb V}(E_s,p_s)}$. The homology class of the
closure $\overline {{\mathbb V}(E,p)}$ is independent of
$p$. 
\end{theorem}
%%%%%%%%%%%%%%%%%%%%%%%%%%%%%%%%%%%%%%%%%%%%%%%%%%%%%%%%%%
\begin{proof} The independence of the homology
class of  $\overline {{\mathbb V}(E,p)}$ is result of the
$\mathrm{PGL}(1)$ action on $\G_T$; if $g(p)=p'$ then $g$ takes 
$\overline {{\mathbb V}(E,p)}$ to $\overline {{\mathbb
V}(E,p)}$. 
The rest is a consequence of homology theory on the 
nonsingular projective algebraic variety $\G_T$.
\end{proof}\par
\noindent {\bf Remark:} Despite the counterexample of
Yam\'eogo (see (Ci) of \S \ref{sec4F}), the intersection $\Z$ of two
such conditions often is proper, and there is an actual
ramification calculus in those cases. \\

\noindent There are two major simplifications of the
problem of determining the homology ring structure
$H^{\ast}(\G_T)$.\\

%%%%%%%%%%%%%%%%%%%%%%%%%%%%%%%%%%%%%%%%%%%%%%%%%%%%%%%%%%
\par{\bf i. Decomposition into a product of elementary}
$G_{T_{k}}.$  
%%%%%%%%%%%%%%%%%%%%%%%%%%%%%%%%%%%%%%%%%%%%%%%%%%%%%%%%%%
We may assume that $T$ has no consecutive constant  values
less than $\mu (T)$. We say that $T$ is {\em elementary} if
$$t_{i}\neq t_{i+1}\, \  for\, \  t_{i}<\mu (T)\, .$$
%%%%%%%%%%%%%%%%%%%%%%%%%%%%%%%%%%%%%%%%%%%%%%%%%%%%%%%%%%
%%%%%%%%%%%%%%%%%%%%%%% Lemma n¼ 7.1 %%%%%%%%%%%%%%%%%%%%%
%\newtheorem{lemm71}{Lemma}[section] \begin{lemm71}
\begin{lemma}\label{lemma7.1} There is
a decomposition of $\G_T$ as a product  
\begin{equation}\label{e7.1}
\G_T =
\prod_{k}Ê{\mathrm{G_{T_{k}}}} \, \  for\, \  T_{k}\, \  elementary. 
\end{equation}
\end{lemma}\noindent
For details see \cite{Y3,IY1}.\\
%%%%%%%%%%%%%%%%%%%%%%%%%%%%%%%%%%%%%%%%%%%%%%%%%%%%%%%%%%
\par
%%%%%%%%%%%%%%%%%%%%%%%%%%%%%%%%%%%%%%%%%%%%%%%%%%%%%%%%%%
{\bf ii. Surjectivity of $\iota^{\star}$.}
%%%%%%%%%%%%%%%%%%%%%%%%%%%%%%%%%%%%%%%%%%%%%%%%%%%%%%%%%%
A result of A. King and C. Walter \cite{KW} shows that the inclusion
%%%%%%%%%%%%%%%%%%%%%%%%%%%%%%%%%%%%%%%%%%%%%%%%%%%%%%%%%%
$\iota : \G_T \rightarrow \BGrass(T)$ of \eqref{e4.2} into a product
of big Grassmannian varieties satisfies
%%%%%%%%%%%%%%%%%%%%%%%%%%%%%%%%%%%%%%%%%%%%%%%%%%%%%%%%%%
\begin{uthm} {\cite{KW}}  The homomorphism $\iota^{\star} : 
H^{\star}(\BGrass(T)) \rightarrow H^{\star}(\G_T)$ of homology rings is a
surjection (and $i_\star :H^\star (\G_T)\rightarrow H^\star ( \BGrass(T))$ is an
inclusion).
\end{uthm} 
%%%%%%%%%%%%%%%%%%%%%%%%%%%%%%%%%%%%%%%%%%%%%%%%%%%%%%%%%%

%%%%%%%%%%%%%%%%%%%%%%%%%%%%%%%%%%%%%%%%%%%%%%%%%%%%%%%%%%
\noindent {\bf Remark.} {\sc status of problem of
determining} $H^{\ast}(\G_T)$.\\
%%%%%%%%%%%%%%%%%%%%%%%%%%%%%%%%%%%%%%%%%%%%%%%%%%%%%%%%%%
{\em Cases where the ring structure} $H^{\star}(\G_T)$
{\em is known:}

   A. For  $T=(1,2,3,2,1)$ (see J. Yam\'{e}ogo, cite{Y3,Y6}).

   B. For  $T = T(\mu,j) = (1,2,\ldots ,\mu ,\mu ,\ldots ,\mu ,1)$ of
socle degree j. 
(see  \S \ref{sec5B} below, and a sequel article).
  
   C. We have determined the class $\iota_{\star}(\G_T)$ in
$\BGrass(T)$ for $T = (1,2,\ldots ,\mu ,a,b,0)$, using a
vector bundle argument suggested by  G. Ellingsrud (See
\cite{IY1}, and a sequel article under preparation).\\
%%%%%%%%%%%%%%%%%%%%%%%%%%%%%%%%%%%%%%%%%%%%%%%%%%%%%%%%%%

\noindent By \eqref{e7.1} these results determine
$H^{\ast}(\G_T)$ for any $T$ composed of the elementary
Hilbert functions in  A. and B.\\

%%%%%%%%%%%%%%%%%%%%%%%%%%%%%%%%%%%%%%%%%%%%%%%%%%%%%%%%%%

\noindent {\em Sections of a vector bundle.} J.~Yam\'{e}ogo
has  shown that if we set $T' = (t_{0},\ldots ,t_{j-1})$,
%%%%%%%%%%%%%%%%%%%%%%%%%%%%%%%%%%%%%%%%%%%%%%%%%%%%%%%%%%
then $\G_T$ is the zeroes of a section of a known
vector  bundle ${\cal E}_{T}$ on $\mathrm{G_{T'}}\times
\Grass(j+1-t_{j},j+1)$, where the codimension of $\G_T$ in
the product is equal  to the rank of ${\cal E}_{T}$.
[Y6].
%%%%%%%%%%%%%%%%%%%%%%%%%%%%%%%%%%%%%%%%%%%%%%%%%%%%%%%%%%
This allows computation of the class of [$\G_T$] in
$H^{\ast}(\mathrm{G_{T'}})\otimes H^{\ast}(\Grass (j+1-t_{j},j+1))$,
%%%%%%%%%%%%%%%%%%%%%%%%%%%%%%%%%%%%%%%%%%%%%%%%%%%%%%%%%%
and with the theorem of A.~King and C.~Walter \cite{KW} in principle allows 
calculation of $H^{\ast}(\G_T)$ inductively, since then
%%%%%%%%%%%%%%%%%%%%%%%%%%%%%%%%%%%%%%%%%%%%%%%%%%%%%%%%%%
$$H^{\ast}(\G_T) \cong
[\G_T]\,\cdot\,(H^{\ast}(\mathrm{G_{T'}})
\otimes
H^{\ast}(\Grass(j+1-t_{j},j+1)),$$ 
%%%%%%%%%%%%%%%%%%%%%%%%%%%%%%%%%%%%%%%%%%%%%%%%%%%%%%%%%%
the ideal of [$\G_T$] in the homology ring of the
product.  However, it is a nontrivial further step to
identify the ring structure $H^{\ast}(\G_T)$ in terms 
of the classes $b_{E}  = [\mathbb{V}(E)]$ of the cells. It is
this further step, open in general, that is needed to convert Theorem \ref{thm7.0}
to an effective calculus of ramification conditions. 
%%%%%%%%%%%%%%%%%%%%%%%%%%%%%%%%%%%%%%%%%%%%%%%%%%%%%%%%%
%%%%%%%%%%%%%%%%%%% Section n¼ 7B %%%%%%%%%%%%%%%%%%%%%   
\subsection{The ring $H^{\ast}(\G_T)$ when $T=T(\mu,j)$}\label{sec5B}
%%%%%%%%%%%%%%%%%%%%%%%%%%%%%%%%%%%%%%%%%%%%%%%%%%%%%%%%%% 
  In this section we determine  the 
ring $H^{\ast}(\G_T)$ for the special case $T=T(\mu,j)$:
\begin{equation}
T(\mu,j)=(1,2,\l\dots
\mu,\ldots , \mu,1=T_j).
\end{equation}
 We then  use the calculation to illustrate
Theorem \ref{thm7.0}, determining the number of ideals of linear systems over
${\check{\mathbb P}}^1$ satisfying an intersection of ramification
conditions (Example \ref{ex7.4}). The ring structure of
$H^{\ast}(\G_T)$, $T=T(\mu,j)$ involves binomial
coefficients,  which are related to the degrees of ideals
of  determinantal minors of circulant matrices
	$$\left( 
\begin{array}{ccc}
a_{0} &\cdots &a_\mu\\
a_{1} &\cdots &a_{\mu +1}\\
\vdots &\vdots &\vdots\\
a_{j-\mu} &\cdots &a_{j}
\end{array} 
\right) \ .$$
%%%%%%%%%%%%%%%%%%%%%%%%%%%%%%%%%%%%%%%%%%%%%%%%%%%%%%%%%%
The variety $\G_T(\mu,j)$ has dimension $2\mu-1$, and has 
the  structure of a projective ${{\mathbb P}}^{\mu-1}$ bundle 
over ${{\mathbb P}}^{\mu}$.                                     
An ideal $I$ having Hilbert function $T$ is determined 
by  $I_\mu=\langle f\rangle $ and by a codimension one subspace $I_{j}$
of $R_{j}$, since $I_{i}=(f)\cap R_{i}$ for $\mu <i<j$. 
Thus, we have an inclusion of projective varieties:\\
\begin{align}
\iota : \, \G_T  &\subset \prod_{i=\mu,j} (\text{Large Grassmanians})=
\Grass(1,\mu +1)\times \Grass(j,j+1)\notag\\				 
				 I &\rightarrow (I_\mu,I_{j}) \subset {{\mathbb P}}^{\mu}\times
{{\mathbb P}}^{j}.\label{epispecial}
\end{align}

%%%%%%%%%%%%%%%%%%%%%%%%%%%%%%%%%%%%%%%%%%%%%%%%%%%%%%%%%  
The projection $pr_{1}: \G_T \rightarrow {{\mathbb P}}^{\mu}$ 
has fibre over  $\langle f\rangle \in R_{\mu}$ the projective space
${{\mathbb P}}(R_{j}/(R_{j-\mu}f))\cong {{\mathbb P}}^{\mu-1}$.   
By this fact, or by Theorem \ref{thm5.1}, the additive structure  
of $H^{\ast}(\G_T)$ is\\
%%%%%%%%%%%%%%%%%%%%%%%%%%%%%%%%%%%%%%%%%%%%%%%%%%%%%%%%%
$$\begin{array}{ll}
			H^{\ast}(\G_T) & \cong_{add}\,
						 H^{\ast}(\Grass(\mu -1,\mu))\otimes
					  H^{\ast}(\Grass(1,\mu))\\
						 & \cong_{add}\,   H^{\ast}({{\mathbb P}}^{\mu -1})\otimes
H^{\ast}({{\mathbb P}}^{\mu })\,.
\end{array}$$
By \eqref{e4.23} the block ${\cal B}_{j-1}(T)$ is $(\mu -1)\times 1$,
and ${\cal B}_{j}(T)$  is  $(1\times \mu )$.               
We let $[a,b]$ denote the cohomology class corresponding 
to the cell
%%%%%%%%%%%%%%%%%%%%%%%%%%%%%%%%%%%%%%%%%%%%%%%%%%%%%%%%%
$$E[a,b] = {\cal D}^{-1}((\mu -1-a)^{\wedge}\, ,(\mu -b),\,  0
\leq a \leq \mu -1,\,  0 \leq b \leq \mu.$$
%%%%%%%%%%%%%%%%%%%%%%%%%%%%%%%%%%%%%%%%%%%%%%%%%%%%%%%%%%
of codimension $a+b$ in $\G_T$. We let $[a,b] =0$ when
$a > \mu -1$ or $b > \mu$.  Recall that  $H^{\ast}(\G_T)$
denotes  the homology ring  of $\G_T$ graded by codimension
of the homology class. If $Y$ is subvariety we denote by
$[Y]$ its homology class in the relevant homology ring.\\
%%%%%%%%%%%%%%%%%%%%%%%%%%%%%%%%%%%%%%%%%%%%%%%%%%%%%%%%%%
%%%%%%%%%%%%%%%%% Lemma n¼ .2 %%%%%%%%%%%%%%%%%%%%%%%%%%%
%\newtheorem{lemm72}[lemm71]{Theorem}
\begin{lemma}\label{lemma7.2}
Suppose $T=T(\mu,j)$.  
Let $\zeta$ be the cohomology class of a hyperplane 
section of ${{\mathbb P}}^{\mu}$, and $\eta$ the class of a
hyperplane section of ${{\mathbb P}}^{j}$. 
Then the  class of $\iota_{\ast}(\G_T)$ satisfies
\begin{equation}\label{e7.1}
[\iota_{\ast}(\G_T)]\,=
\,(\zeta +\eta)^{j+1-\mu}\,.
\end{equation}
Furthermore, $\iota_{\ast}$ is injective, and
$\iota^{\ast}$ is surjective  on homology or 
cohomology and they satisfy
\begin{equation}\label{e7.1a}
\iota^{\ast}(\zeta^{u}\eta^{v})\,=  
										\left\{  
										\begin{array}{l}
										[u,v]\, \   for\, \ u+v<\mu, \\
										\sum_{i=0}^{i=j+1-\mu}
										(_{i}^{j+1-\mu})[u+i-1,v-i+1]\, \ for\, \  u+v\geq
										\mu.
									\end{array}
									\right.
\end{equation}
\begin{equation}\label{e7.1b}
\iota_{\ast}[a,b]\,=\, 
			\left\{
				\begin{array}{l}
				[\zeta^{a+1}\eta^{b+j-\mu}]\, \  for\, \ 
   a+b\geq \mu,\\
			\zeta^{a}\eta^{b}[\iota_{\ast}\G_T]\, \   for\, \ 
    a+b < \mu.
			\end{array}
				\right.
\end{equation}
\end{lemma}
%%%%%%%%%%%%%%%%%%%%%%%%%%%%%%%%%%%%%%%%%%%%%%%%%%%%%%%%%

\noindent {\bf Comment.} The $\zeta^{i}$ term in \eqref{e7.1} is
zero unless $i\leq \mu$, as $\zeta^{\mu +1}=0$. We originally proved
Lemma~\ref{lemma7.2} and Theorem \ref{thm7.3} by intersecting
$\G_T$ with cycles of $H^\ast({{\mathbb P}}^{\mu}\times {{\mathbb P}}^{j})$
having complementary dimension, obtaining the degree of
certain determinantal minors of circulant matrices, which
we then calculated. We are grateful to G.~Ellingsrud, who
suggested the approach here. A similar method applies to
Hilbert functions of the types $T=(1,2, \ldots , \mu ,
\ldots ,\mu ,a,0)$ or $T = ( 1,2, \ldots ,\mu ,a,b,0)$.
That the formula \eqref{e7.1} depends on $j$ shows that the homology ring
structure is not that of $\SGrass(T)$ which here depends only on $\mu$.

\begin{proof}[Proof of Lemma \ref{lemma7.2}] On ${{\mathbb P}}^{\mu}\cong
\Grass(1,R_{t})$ we let $-\zeta = {\cal O}(-\zeta)\cong
{\cal O}(1)$ denote the universal rank one subbundle, and
on ${{\mathbb P}}^{j}$ we let ${\cal O}(\eta) \cong {\cal
O}(-1)$ denote the universal rank one quotient. Let
$p_{1}$, $p_{2}$ denote the projections from ${{\mathbb P}} =
{{\mathbb P}}^{\mu}\times {{\mathbb P}}^{j} \rightarrow {{\mathbb P}}^{\mu} \ or
\ {{\mathbb P}}^{j}$, respectively. Then we have on 
${{\mathbb P}}^{\mu}\times {{\mathbb P}}^{j}$ a composite map
$$\alpha : \ p_{1}^{\ast}(R^{j-\mu}\otimes {\cal O}(-\zeta))
\rightarrow p_{2}^{\ast}(R^{j})\rightarrow
p_{2}^{\ast}({\cal O}(\eta))$$
whose vanishing at $p\in {{\mathbb P}}^{\mu}\times {{\mathbb P}}^{j}$ is
a condition for $p\in \G_T$. The vanishing of $\alpha$ is
equivalent to that of $\alpha^{\vee}$
$$\alpha_{\vee} : P_{2}^{\ast}({\cal O}(-\eta))
\rightarrow p_{1}^{\ast}(R^{\vee j-\mu}\otimes {\cal
O}(\zeta)),$$
or that of $\alpha ' = \alpha^{\vee}\otimes {\cal O}(\eta)$
\begin{equation}\label{e7.2}
\alpha ' : {\cal O}_{{{\mathbb P}}} \rightarrow
(p_{1}^{\ast}(R^{\vee j-\mu}\otimes{\cal O}(\zeta))\otimes
{\cal O}(\eta ) = {\cal O}_{{{\mathbb P}}}^{j+1-\mu}(\zeta + \eta
) = {\cal E}. 
\end{equation}
Since the maps are linear on fibers of ${{\mathbb P}}$ over ${{\mathbb P}}^{\mu}$, and
since the codimension of $\G_T$ in ${{\mathbb P}}$  satisfies 
\begin{equation*}
[\mu+j-(2\mu+1)] =
j+1-\mu = rk({\cal E}),
\end{equation*}
 and since
$\G_T$ is nonsingular, the class $c_{j+1-\mu}({\cal E}) =
[V(\alpha ')]$ is reduced and equal to $[\G_T]$. This
class is $(\zeta + \eta )^{j+1-\mu}$, proving \eqref{e7.1}.\\
We omit the proofs of the remaining formulas of Lemma \ref{lemma7.2}, which are
now straightforward.
\end{proof}

%%%%%%%%%%%%%%%%%%%%%% Theorem n¼ 7.3 %%%%%%%%%%%%%%%%%%%
%\newtheorem{theo73}[lemm71]{Theorem} \begin{theo73}
\begin{theorem}\label{thm7.3}
																When $T=T(\mu,j)$ the product 
																$[a,b]\cdot [c,e]$ satisfies
$$ [a,b]\cdot [c,e]=
\left\{
\begin{array}{ll}
[a+c,b+e] &\left\{ 
															   \begin{array}{c}
																		{\scriptstyle if\, \ a+b+c+e<\mu\, \ or}\\
																		{\scriptstyle a+b<\mu\, \ and\, \ c+e\geq
																			\mu,}\\
																		{\scriptstyle  or\, \ vice-versa}
																		\end{array} \right.\\
\sum_{i=0}^{i=j+1-\mu}(_{i}^{j+1-\mu})
								[a+c+i-1,b+e-i+1]
        &\left\{ 
																		\begin{array}{c}
																		{\scriptstyle if\, \ a+b<\mu,\, \ 
																				c+e<\mu}\\
																		{\scriptstyle but\, \ a+b+c+e \geq
																				\mu}
																	 \end{array}\right.\\
0 						&{\quad \, \scriptstyle if\, \ a+b\geq \mu\, \ and\, \
										c+e\geq \mu}
\end{array} \right. 
$$
\end{theorem}
%%%%%%%%%%%%%%%%%%%%%%%%%%%%%%%%%%%%%%%%%%%%%%%%%%%%%%%%

\begin{proof} The proof of Theorem \ref{thm7.3} is a
straightforward calculation, given Lemma \ref{lemma7.2}.
\end{proof}\par
%%%%%%%%%%%%%%%%%%%%%%%%%%%%%%%%%%%%%%%%%%%%%%%%%%%%%%%%

\noindent {\bf Remark.} When the classes
$[a,b]$ and  $[c,e]$ are in complementary dimensions, so
$(a+b+c+e) = 2\mu -1$,
%%%%%%%%%%%%%%%%%%%%%%%%%%%%%%%%%%%%%%%%%%%%%%%%%%%%%%%%%%
then $a+b<\mu$ implies $c+e \geq \mu$, and vice-versa. Thus the
product 
$$[a,b]\cdot [c,e]=[a+c,b+e]=[a+c,2\mu -1-(a+c)],$$
which is nonzero only when $a+c=\mu -1$, and $b+e=\mu$, namely
when $E[a,b]=(E[c,e])^{\vee}$, the class of the dual
partition. 
%%%%%%%%%%%%%%%%%%%%%%%%%%%%%%%%%%%%%%%%%%%%%%%%%%%%%%%%%%
It follows that when $T=T(\mu,j)$ the classes of cells do in
fact give the  exact duality in complementary dimensions
for $H^{\ast}(\G_T)$, unlike what is true in general.\\
%%%%%%%%%%%%%%%%%%%%%%%%%%%%%%%%%%%%%%%%%%%%%%%%%%%%%%%%%%
%%%%%%%%%%%%%%%%%%%%%%%%%%%%%%%%%%%%%%%%%%%%%%%%%%%%%%%%%%
%%%%%%%%%%%%%%%%%%% Example n¼ 7.4 %%%%%%%%%%%%%%%%%%%%%%%

%\newtheorem{exam74}[lemm71]{Example}
\begin{example}[Determining ideals with  given ramification]\label{ex7.4}
Let $\mu=3$ and $j=6$, so that\
$T=(1,2,3,3,3,3,1)$ . We wish to find the
number of ideals having the ramification  
$[1,1]$ at $x$, $[0,2]$ at $y$,
and $[1,0]$ at $(x+y)$
classes of total codimension
$2+2+1 = 5 = \dim (\G_T)$. 
By Theorem~\ref{thm7.3}, the product 
$[1,1]\cdot [0,2]$ 
satisfies 
$$ [1,1]\cdot [0,2]\, =\, 
\sum_{0\leq i
\leq 4}(_{i}^{4})\cdot	[i,4-i]\, .
$$
 By the exact duality, the intersection of this 
product  with $[1,0]$ picks out the coefficient for
$i=1$, namely four. There are thus four ideals -
counting multiplicities, having this ramification at 
$x,y, x+y$.  What are these four ideals?
%%%%%%%%%%%%%%%%%%%%%%%%%%%%%%%%%%%%%%%%%%%%%%%%%%%%%%%%%%
Each of the four ideals has ramification at $x = 0$ 
given by the the powers of $x$ in the monomial ideal  
$E_{1}= (y^{2}x,y^{6},x^{6})$ with cobasis given by the
partition  $(6,6,1,1,1,1)$ of code
%%%%%%%%%%%%%%%%%%%%%%%%%%%%%%%%%%%%%%%%%%%%%%%%%%%%%%%%%%
${\cal D}(E_{1}) = ({\cal Q}_{5}=(1,0);
{\cal Q}_{6}=(2))$. 
%%%%%%%%%%%%%%%%%%%%%%%%%%%%%%%%%%%%%%%%%%%%%%%%%%%%%%%%%%
%%%%%%%%%%%%%%%%%%%%%%%%%%%%%%%%%%%%%%%%%%%%%%%%%%%%%%%%%%
The monomial ideal $E_{2}$ in coordinates $(y,C)$ with
code ${\cal D}(E_{2}) = ({\cal Q}_{5} =
(1,1);{\cal Q}_{6} = 1)$ at $y$ is,  $E_{2} =
(y^{6}, Cy^{5}, C^{3})$, of partition $(6,5,5)$.
%%%%%%%%%%%%%%%%%%%%%%%%%%%%%%%%%%%%%%%%%%%%%%%%%%%%%%%%%%
Finally  the monomial ideal $E_{3}$ in coordinates
$(x+y, C)$ for the diagram ${\cal D}(E_{3}) =
({\cal Q}_{5} = 1; {\cal Q}_{6} = 3)$ at $x+y$ is that of 
the ideal $((x+y)^{7}, C(x+y)^{5}, C^{2}(x+y), C^{6})$, 
with cobasis given by the partition $(7,5,1,1,1,1)$.
%%%%%%%%%%%%%%%%%%%%%%%%%%%%%%%%%%%%%%%%%%%%%%%%%%%%%%%%%%
It follows that  the generator  $f = f_{a}$ satisfies
	$$f=x(x+y)(x+ay)\, \  where\, \  a\, \  is\, \  not\, \ 
	zero\, .$$
%%%%%%%%%%%%%%%%%%%%%%%%%%%%%%%%%%%%%%%%%%%%%%%%%%%%%%%%%
We know that $I_{6}$ contains $y^{6}, xy^{5}$, and $x^{6}$,
as well as  $R_{3}f$, but has dimension only $6$.       
It follows that there  is a dependency in the columns of
the $4\times 4$ matrix $M$ obtained by writing as rows
$x^{3}f, x^{2}yf, xy^{2}f, y^{3}f$  of \      
$(R_{6}f/\langle y^{6}, xy^{5}, x^{6}\rangle )$ in the  basis  $\langle x^{5}y,
x^{4}y^{2}, x^{3}y^{3}, x^{2}y^{4}\rangle $ for \  $(R_{6}/\langle y^{6},
xy^{5}, x^{6}\rangle )$:
%%%%%%%%%%%%%%%%%%%%%%%%%%%%%%%%%%%%%%%%%%%%%%%%%%%%%%%%% 
$$M\, =\, \left( \begin{array}{cccc}
																	1+a &a &0 &0\\
																	1 &1+a &a &0\\
																	0 &1 &1+a &a\\
																	0 &0 &1 &1+a
																	\end{array}
											\right)\, .$$
%%%%%%%%%%%%%%%%%%%%%%%%%%%%%%%%%%%%%%%%%%%%%%%%%%%%%%%%%%
The four roots $a_{1},\ldots ,a_{4}$ of $\det M = 0$ give
the four  ideals
\begin{equation}\label{e7.3}
I_{u} = (f_{a}, y^{6}, xy^{5},x^{6}),\, \  a =
a_{u},\,
\   u = 1,\ldots ,4\,. 
\end{equation}
having the given ramification at $x, y$, and $x+y$.
\end{example}
%%%%%%%%%%%%%%%%%%%%%%%%%%%%%%%%%%%%%%%%%%%%%%%%%%%%%%%%%%
%%%%%%%%%%%%%%%%% Proposition n¼ 7.5 %%%%%%%%%%%%%%%%%%%%%
%\newtheorem{propos75}[lemm71]{Theorem}
\begin{proposition}\label{prop7.5}
If $I\in \G_T$,
the sum of the codimensions or lengths  of the
ramification   conditions it satisfies at all points of
${\check{\mathbb P}}^1$ is  
\begin{equation}\label{e7.4}
\begin{array}{ll}
\sum_{p\in {\check{\mathbb P}}^1}\ell (\QRAM(I,p))
									&= \dim(\BGrass(T) )\\
  							&= \sum_{\mu (T)\leq i\leq j(T)} (t_{i})(i+1-t_{i})
\end{array}
\end{equation}
\end{proposition}
%%%%%%%%%%%%%%%%%%%%%%%%%%%%%%%%%%%%%%%%%%%%%%%%%%%%%%%%%%
\begin{proof} Immediate from Lemma \ref{lemma1.4} applied
to each $I_{i}$.
\end{proof}
%%%%%%%%%%%%%%%%%%%%%%%%%%%%%%%%%%%%%%%%%%%%%%%%%%%%%%%%%%
%\newtheorem{exam76}[lemm71]{Example}
\begin{example}[Overdetermination of ideals by
ramification]\label{ex7.6}
Each of the four solution ideals $I$ for
Example \ref{ex7.4}  has non-generic ramification at
$(x+a_{u}y)$  because $x+a_{u}y$ is a factor of
$f$.
%%%%%%%%%%%%%%%%%%%%%%%%%%%%%%%%%%%%%%%%%%%%%%%%%%%%%%%%%%
Each solution ideal $I(u), u = 1,\ldots ,4$ has
ramification also at three other points $L_{uv} = 0, v =
1, 2, 3$ satisfying $L_{uv}^{6}\in I(u)_{6}$ (see Lemma \ref{lemma1.4}).
\end{example}
%%%%%%%%%%%%%%%%%%%%%%%%%%%%%%%%%%%%%%%%%%%%%%%%%%%%%%%%%%
We define a Wronskian morphism $W : \G_T \rightarrow {{\mathbb P}}$, ${{\mathbb P}} = \prod {{\mathbb P}}^{N_{i}}\, , N_{i} =
t_{i}(i+1-t_{i})$ as the \linebreak composition of $\iota : \G_T
\rightarrow \BGrass(T)$ and of the product $w$ of the
Wronskisn maps \linebreak
$w_{i} : \Grass(i+1-t_{i}, i+1) \rightarrow {{\mathbb P}}^{N_{i}}$ (Definition
\ref{def1.4}). Each $w_{i}$ is a finite cover of degree given by Proposition
\ref{prop2.2}. Hence we have
%%%%%%%%%%%%%%%%%%%%%%%%%%%%%%%%%%%%%%%%%%%%%%%%%%%%%%%%%%
\begin{proposition}\label{prop7.7} The Wronskian morphism 
$W : \G_T \rightarrow {{\mathbb P}}$ is a finite cover of its
image, of degree $\prod_{i=\mu}^j\deg w_i$.
\end{proposition}
%%%%%%%%%%%%%%%%%%%%%%%%%%%%%%%%%%%%%%%%%%%%%%%%%%%%%%%%%%
\noindent {\bf Remark.} {\sc relations among the
ramification conditions}. It follows from Propositions
\ref{prop7.5}  and \ref{prop7.7} that in general there are relations between
the set of ramifications $\{\QRAM(I,L_{p}) \ | \ p
\in {{\mathbb P}}^{1}\}$ in different directions $p$ of an ideal
$I\in \G_T$.
%%%%%%%%%%%%%%%%%%%%%%%%%%%%%%%%%%%%%%%%%%%%%%%%%%%%%%%%%%
This contrasts with the case of a single linear system,
where the total codimension of the ramification conditions
adds up to the dimension $N_{i}$ of $\Grass (d_{i},R_{i})$
(Lemma \ref{lemma1.4}). For ideals the total codimension of the
ramification conditions adds up to the dimension of the
product $\BGrass(T)$ of big Grassmannians, in which $\G_T$
is embedded; this is much larger than the dimension of
$\G_T$. We expect that
there is a kind of algebraic variety of relations among the
ramification conditions for ideals in $\G_T$.
%%%%%%%%%%%%%%%%%%%%%%%%%%%%%%%%%%%%%%%%%%%%%%%%%%%%%%%%%%
This can be seen in the special case of Example \ref{ex7.6}, where
evidently the ``extra" ramification points $x+a_{i}y = 0$,
and at $L_{uv} = 0$, for $I(u)$ satisfy algebraic
relations. What are the equations defining the image
$W(\G_T)$ in ${{\mathbb P}}$? \\\medskip\par\noindent
%%%%%%%%%%%%%%%%%%%%%%%%%%%%%%%%%%%%%%%%%%%%%%%%%%%%%%%%%%
%%%%%%%%%%%%%%%%%%%%%%%% section n¼ 7C %%%%%%%%%%%%%%%%%%%
{\bf Problem}\, {\sc A Generalized Schubert calculus?}
%%%%%%%%%%%%%%%%%%%%%%%%%%%%%%%%%%%%%%%%%%%%%%%%%%%%%%%%%%
The intersection of the classes of cells in\linebreak
$H^{\ast}(\G_T),\  T = T(\mu,j)$ can be expressed simply,
using the codes $[a,b]$ of cells; the intersection 
numbers involve binomial coefficients.  
%%%%%%%%%%%%%%%%%%%%%%%%%%%%%%%%%%%%%%%%%%%%%%%%%%%%%%%%%%
When $d=j$, then  $T = (1,2,\ldots ,\mu,a,0)$, and  $\G_T$ is
a Grassmann variety $\Grass (\mu+1-a,R_\mu)$. In this case the
cells are the Schubert cells, and the intersection numbers
are given by  the Littlewood-Richardson rule. 
%%%%%%%%%%%%%%%%%%%%%%%%%%%%%%%%%%%%%%%%%%%%%%%%%%%%%%%%%%
When ${k} = {\mathbb C}$, $H^{\ast}(\G_T)$ is
additively a ${\bf Z}$-module with basis the classes of
the cells ${\mathbb V}(E({P}))$, as $E({P})$ runs
through the monomial ideals $E({ P})$ attached to
partitions $P$  of diagonal lengths $T$. The question
of finding the ring structure on $H^{\ast}(\G_T)$ in terms
of this basis  asks for a generalization of (at least) two
rules. First, the Schubert calculus and its
Littlewood-Richardson rule when $\mu =j$; second, the
intersection numbers given above for $T=T(d,j)$
involving binomial coefficients (Theorem \ref{thm7.3}). The
intersection numbers there are most simply  expressed in
terms of the hook code of the cells.

%%%%%%%%%%%%%%%%%%%%%%%%%%%%%%%%%%%%%%%%%%%%%%%%%%%%%%%%%%

The problem of determining the homology ring of $\G_T$, generalizing
the special case $\mu=j$ which is Schubert calculus, is, given two partitions $P,Q$ of
diagonal lengths  $T$, to find the intersection numbers
$\alpha(P, Q:S)$, such that the homology
classes    $[P],\ [Q]$ of the cells
$\mathbb{V}(P)$ and $\mathbb{V}(Q)$ satisfy
\begin{equation}\label{e7.5}
[P]\,\cdot \,[Q] =
\sum_{S}\alpha(P,Q:S)[S] ,
\end{equation}
%%%%%%%%%%%%%%%%%%%%%%%%%%%%%%%%%%%%%%%%%%%%%%%%%%%%%%%%%%
Here $S$ runs through the partitions of diagonal
lengths $T$, having codimension $\cod\mathbb{V}(S) =
\cod\mathbb{V}(P) + \cod\mathbb{V}(Q)$. 
%%%%%%%%%%%%%%%%%%%%%%%%%%%%%%%%%%%%%%%%%%%%%%%%%%%%%%%%%%
The related geometric problem is to connect these numbers
naturally to the geometry of the rational normal curve. 
%%%%%%%%%%%%%%%%%%%%%%%%%%%%%%%%%%%%%%%%%%%%%%%%%%%%%%%%%%
%%%%%%%%%%%%%%%%%%%%%%%%%%%%%%%%%%%%%%%%%%%%%%%%%%%%%%%%%%
\subsection{Desingularization of the secant bundle $V(\mu,j)$
to the rational normal curve}\label{sec5C}
Let $k[A]$ be the polynomial ring $k[A]=k[a_0,\ldots ,a_j]$,
suppose $2\mu <j+1$ and 
denote by $V(\mu ,j)$ the projective subvariety of 
${\mathbb P}^j=\mathrm{Proj}(k[A])$ defined by $I(\mu,j-\mu ,A)$, the ideal of $
(\mu +1)\times (\mu +1)$ minors of the Hankel matrix
\begin{equation}\label{e7.5}
\mathrm{HANKEL}(\mu ,j-\mu ,A)=\left(
\begin{array}{cccc}
a_0&a_1&\ldots &a_{j-\mu }
\\
a_1&a_2&\ldots &a_{j+1-\mu }
\\ 
\vdots & \vdots&\ldots &\vdots
\\ 
a_\mu &a_{\mu +1}&\ldots &a_{j}
\end{array}
\right) . 
\end{equation}
It is well known that $V(\mu ,j)$ is irreducible, and that 
$I(\mu ,j-\mu ,A)$ is also the ideal of $(\mu +1)\times (\mu +1)$ minors of 
$\mathrm{HANKEL}(\mu ',j-\mu ',A)$ provided $\mu\leq \min(\mu ',j-\mu ')$ (see
\cite{GP,watan,eis,geram},\cite[\S 1.3]{IK}).
It is also well known that if $\cha k=0$ or $\cha k>j$, 
then $V(\mu ,j)$ is the $\mu $-secant locus $\Sec (\mu ,j)$ to the degree-$j$ 
rational normal curve $X(j)$ 
 parametrizing projective ${{\mathbb P}}^{\mu -1}$-secants to 
$X(j)$.
We explain this connection briefly. Suppose that  
$$F_a = a_{0}x^j+{j\choose 1}a_1x^{j-1}y+\ldots +
{j\choose i}a_ix^{j-i}y^i+\ldots +a_{j}y^j \ \  (a_i \in k)$$
is an element of $R_j$, 
 $R = k[x,y]$. $\Sec (\mu ,j)$ is the closure of the family of 
forms $F$ that can be written
as the sum of $F=L_1^j+\ldots L_\mu ^j$ of $\mu $ $j$-th powers
 of linear forms $L_i \in R_1$.
 Then 
\begin{lemma}\label{lemma7.8}
$F_a \in \Sec (\mu ,j)$ iff $a \in V(\mu ,j)$. 
\end{lemma} 
We now explain the connection with the variety $\G_T, T = T(\mu ,j)$.
 Following the classical theory of apolarity, let $R$ act on $R$ as 
higher order partial differential operators:
$h\circ f = h(\partial / \partial x, \partial /\partial y) \circ f$.
 Then if $F \in R_j$,
the ideal $I_F = \Ann(f)$ is a Gorenstein ideal - hence complete
 intersection of $R$, and
we have  $I_F = (g,h)$, with $\mu =\mu _F= \deg (g) \leq \deg(h)=j+2-\mu $
 (see \cite[Theorem 1.54]{IK}). We have
\begin{lemma}\label{lemma7.9}
$\mu _F \leq \mu  $ iff $F \in \Sec (\mu ,j)$. 
\end{lemma}
Let $T = T(\mu ,j)$. We define a morphism ${\alpha } : G_T \rightarrow 
\Sec (\mu ,j)$ by $\alpha (I) = F_I, 
F_I = (I_j)^{\perp} 
\in R_j$. Thus, $\alpha$ is identical to (the dual of) $pr_2:
 G_T \rightarrow {{\mathbb P}}^j$,
 in the notation of \S \ref{sec5B}. Since, as is easy to show,
 $I_\mu  = g$ satisfies $g\circ F_I=0$, we have $\mu _{F_I}
\leq \mu $, implying by the Lemma that $F_I \in \Sec (\mu ,j)$. \par
We stratify $V(\mu ,j)$ by rank: $V(\mu ,j) \supset V(\mu -1,j) \supset
 \ldots \supset V(1,j)$.
It is easy to see that the strata
satisfy $\dim (V(i,j)) = 2i-1$ if $i \leq (j+1)/2$, or $\dim (V(i,j))
 = 2i-2$ if $i = (j+2)/2$. If $2\mu <j+1$  
and $i>0$ the singular locus of $V(i,j)$ is $V(i-1,j)$; if $i=\mu $
and $2\mu =j+1$, or $2\mu =j+2$ then $V(\mu ,j) = {{\mathbb P}}^{j}$ or
$\Grass(2,R_{j})$, respectively, and is non singular. 
\par 
Recall that a desingularization is semismall iff the fiber 
over a codimension $c$
stratum has dimension at most $2c$ (see \S 6 of \cite{nakaj3} or 
\cite{bor-macpher}). Below
 $\iota: G_T \rightarrow {{\mathbb P}}^\mu  \times {{\mathbb P}}^j$
is the inclusion of \eqref{epispecial}. We have
\begin{theorem}\label{thm8.0}
If $2\mu <j+1$, the morphism $\alpha$ makes $G_T$ a semismall 
desingularization of $\Sec (\mu ,j)$, whose fiber over a point 
$p\in V(i,j)-V(i-1,j)$, for $ 1 \leq i \leq \mu $ is a
 projective space ${{\mathbb P}}^{\mu -i}$. The
class in the homology ring $H^{\ast}(G_T)$ of ${\alpha}^{-1}(V(i,j))$
 satisfies 
\begin{equation}\label{e7.6a}
|{\alpha}^{-1}(V(\mu -1,j))| = 
{\iota}^{\ast}(\mu \cdot
\eta  -(j+2-2\mu )\zeta)=\mu [0,1]-(j+2-2\mu )[1,0],
\end{equation}
 and
\begin{equation}\label{e7.6b}
|{\alpha}^{-1}(V(i,j))|= \ \mbox{coeff of} \ 
t^{\mu -i}
\ 
\mbox{in} \ 
\iota^{\ast}\lbrace (1-\zeta t)^{j-\mu -i+1}(1+\eta t)^{i+1}) \rbrace 
\end{equation}
\begin{equation}\label{e7.6c} \ \  = \sum_{u+v=\mu -i}
{{j-\mu -i+1}\choose{u}}{{i+1}\choose{v}}[u,v]. 
\end{equation}
\end{theorem}
\begin{proof} First, if the point $p\in V(i,j)-V(i-1,j)$ 
corresponds to $F$, then $I_F = (g',h')$,
with $\deg(g') = i$, and the fibre of $G_T$ over $p$ is all
 pairs $(I_\mu ,I_j)$ with $I_j=(I_F)_j$
and $g'|I_\mu $; thus, $I_\mu  = <g"g'>,g" \in R_{\mu -i}$  and the
 fiber over $p$ is 
${{\mathbb P}}(R_{\mu -i}) = {{\mathbb P}}^{\mu -i}$. 
\par
Consider $\G_T \subset {{\mathbb P}}^\mu  \times {{\mathbb P}}^j$, let $pr_{1}$ be 
the projection of $\G_T $ on ${{\mathbb P}}^\mu $ and $pr_{2}$ its projection 
on ${{\mathbb P}}^j$. Let ${\mathcal S}_{\mu }$ denote the pull-back on
$\G_T $ of the tautological sub-bundle of ${{\mathbb P}}^\mu $ and 
${\mathcal Q}_{j}$ the pull-back on $\G_T $ of the tautological 
quotient-bundle of ${{\mathbb P}}^j$.
For $1\leq i\leq \mu $, we have an injection ${\mathcal S}_{\mu }\otimes
R_{j-\mu -i}\hookrightarrow R_{j-i}$ of vector bundles on $\G_T $ 
that on fibres maps $f\otimes h$ to the product $fh$. We identify 
${\mathcal S}_{\mu }\otimes R_{j-\mu -i}$ to a sub-bundle of $R_{j-i}$
and let $E_{i} = R_{j-i}/({\mathcal S}_{\mu }\otimes R_{j-\mu -i})$. 
Now consider the homomorphism of vectors bundles on $\G_T $, \  
$\phi_{i} : E_{i}\rightarrow {\mathcal Q}_{j}^{\oplus i+1}$, that
on fibres maps the class $\overline{f}$ to $(x^{i}\cdot\overline{f},
\ldots, x^{i-l}\cdot\overline{f},\ldots ,y^{i}\cdot\overline{f})$. It is easy 
to see that $\alpha^{-1}(V(i,j))$ is the locus where $\phi_{i}$ has rank
less or equal $i$. Since $\alpha^{-1}(V(i,j))$ has the right
codimension ($\mu -i$) in $\G_T $, we have 
$[\alpha^{-1}(V(i,j))]= c_{\mu -i}\left({\mathcal Q}_{j}^
{\oplus i+1}\ - \ E\right)$, where 
$$\begin{array}{ll}
c_t\left({\mathcal Q}_{j}^
{\oplus i+1}\ - \ E\right)& = c_t({\mathcal Q}_{j}^{\oplus i+1})
/c_t(E)=c_t(S_\mu )^{j-\mu -i+1}c_t({\mathcal Q}_{j})^{i+1}
\\
&=(1-\zeta t)^{j-\mu -i+1}(1+\eta t)^{i+1}
\end{array}.$$ 
This proves the equation \eqref{e7.6b}. The first equation \eqref{e7.6a}
is just  a particular case of \eqref{e7.6b} when $i=\mu -1$.
\end{proof}

%%%%%%%%%%%%%%%%%%%%%%%%%%%%%%%%%%%%%%%%%%%%%%%%%%%%%%%%%%
\addcontentsline{toc}{section}{\sc references}
%%%%%%%%%%%%%%%% References %%%%%%%%%%%%%%%%%%%%%%%%%%%%%%
\bibliographystyle{amsalpha}

%%%%%%%%%%%%%%%%%%%%%%%%%%%%%%%%%%%%%%%%%%%%%%%%%%%%%%%%%%

\vspace{2cm}
  
%%%%%%%%%%%%%%%%%%%%%%%%%%%%%%%%%%%%%%%%%%%%%%%%%%%%%%%%%%%
%\begin{center}
%\begin{minipage}{10cm}                                  
%A. Iarrobino\\ Mathematics Department, 526B NI\\
%Northeastern University\\ Boston, MA 02115\\
%e-mail: iarrobin@northeastern.edu \\                
%      \\        
%             %%%%%%%%%%%%%%%%%%%%%%%%%%%%%%%           
%Joachim Yam\'{e}ogo\\                                
%IMSP, Universit\'{e} de Nice\\ Parc Valrose\\
%06034 Nice\\
%e-mail: yameogo@math.unice.fr
%\end{minipage} 
%
%\end{center}
%%%%%%%%%%%%%%%%%%%%%%%%%%%%%%%%%%%%%%%%%%%%%%%%%%%%%%%%%%%

\begin{center}
\begin{tabular}{l}                                 
Anthony A. Iarrobino,
\\ 
Mathematics Department, 567 LA
\\
Northeastern University Boston, MA 02115, USA
\\
e-mail: iarrobin@neu.edu 
\\
\\
Joachim Yam\'eogo,
\\
Laboratoire J.-A. Dieudonn\'e, UMR CNRS 6621 
\\                              
Universit\'e de Nice-Sophia Antipolis,
\\
F-06108 Nice cedex 02, France 
\\
e-mail: yameogo@math.unice.fr
\end{tabular} 
\end{center}
\end{document}